\documentclass{elsart}

\usepackage{amssymb,amsmath,graphicx}
\usepackage{amsfonts}
\newtheorem{corollary}{Corollary}
\newtheorem{example}{Example}
\newtheorem{definition}{Definition}
\newtheorem{lemma}{Lemma}
\newtheorem{proposition}{Proposition}
\newtheorem{remark}{Remark}
\newtheorem{theorem}{Theorem}
\newcounter{QQ}
\newtheorem{cla}[QQ]{Composition Lemma For Associative Algebras}
\newtheorem{cll}[QQ]{Composition Lemma For Lie Algebras}
\newtheorem{cbl}[QQ]{Special Bracketing Lemma}
\newtheorem{Proposition}[QQ]{Proposition GS}

\newcounter{AD}
\newtheorem{Definition}[AD]{Definition}

\newcommand{\pp}{\noindent {\em Proof. }}

\newcommand{\bee}[1]{\begin{equation}\label{#1}}
\newcommand{\ene}{\end{equation}}
\newcommand{\beq}[1]{\begin{eqnarray}\label{#1}}
\newcommand{\eqe}{\end{eqnarray}}
\newcommand{\bea}{\begin{eqnarray*}}
\newcommand{\eqa}{\end{eqnarray*}}
\newcommand{\cd}{\cdots}
\newcommand{\cir}{\,{\scriptstyle\circ}\,}
\newcommand{\ld}{\ldots}
\newcommand{\ra}{\rightarrow}
\newcommand{\vp}{\varphi}
\newcommand{\ve}{\varepsilon}

\renewcommand{\o}{\otimes}

\newcommand{\ca}{\mathcal{A}}
\newcommand{\cb}{\mathcal{B}}
\newcommand{\cc}{\mathcal{C}}
\newcommand{\cf}{\mathcal{F}}

\newcommand{\cl}{\mathcal{L}}
\newcommand{\cm}{\mathcal{M}}
\newcommand{\cn}{\mathcal{N}}
\newcommand{\cs}{\mathcal{S}}

\newcommand{\cu}{\mathcal{U}}
\newcommand{\cw}{\mathcal{W}}
\newcommand{\N}{\mathbb{N}}

\newcommand{\dist}{\mathrm{dist}}
\newcommand{\ds}[2]{\mathrm{dist}_{#1}^{#2}}

\newcommand{\dg}{\mathrm{dg}\,}

\renewcommand{\th}{^\mathrm{th}}

\newcommand{\Sp}[1]{\mathrm{Span}\,\left\{#1\right\}}
\newcommand{\da}{d_{\alpha}}

\newcommand{\wX}{\widetilde{X}}

\newcommand{\wq}{\widetilde{q}}
\renewcommand{\wp}{\widetilde{p}}

\newcommand{\ba}{\overline{a}}

\newcommand{\bc}{\overline{c}}
\newcommand{\bd}{\overline{d}}
\renewcommand{\bf}{\overline{f}}
\newcommand{\bg}{\overline{g}}

\newcommand{\bh}{\overline{h}}

\newcommand{\bs}{\overline{s}}

\newcommand{\bS}{\overline{S}}

\newcommand{\epf}{\hfill$\Box$}

\newcommand{\fd}{finite - dimensional }
\newcommand{\fg}{finitely generated }

\newcommand{\fa}{free associative algebra }

\newcommand{\fas}{free associative algebras }

\newcommand{\fac}{free associative algebra, }

\newcommand{\wax}{\mathcal{W}(X)}
\newcommand{\fx}[1]{\mathcal{A}(#1)}
\newcommand{\lax}{\mathcal{L}(X)}

\newcommand{\mS}{\mathcal{S}}

\newcommand{\mw}{\mathcal{W}^{\lambda}}

\newcommand{\mv}{\mathcal{V}^{\lambda}}

\renewcommand{\alg}[1]{\mathrm{alg}\,{#1}}

\begin{document}

\begin{frontmatter}

\title{Filtrations and Distortion in Infinite-Dimensional Algebras}

\author[yb,ybao,yb1]{Yuri Bahturin} and \author[ao,ybao,ao1]{Alexander Olshanskii}
\address[yb]{Department of Mathematics and Statistics, Memorial University of Newfoundland, St. John's, NL, A1C 5S7, \textsc{Canada}}
\address[ao]{Department of Mathematics,
1326 Stevenson Center,
Vanderbilt University, Nashville, TN 37240, \textsc{USA}}
\address[ybao]{Department of Algebra, Faculty of Mathematics and Mechanics, 119899 Moscow, \textsc{Russia}}
\thanks[yb1]{Partially supported by NSERC grant \# 227060-09}
\thanks[ao1]{Partially supported by NSF grant, DMS 0700811, and Russian Fund for Fundamental Research, grant 08-01-00573}

\begin{abstract} A tame filtration of an algebra is defined by the growth of its terms, which has to be majorated by an exponential function. A particular case is the degree filtration used in the definition of the growth of \fg algebras. The notion of tame filtration is useful in the study of possible distortion of degrees of elements when one algebra is embedded as a subalgebra in another. A geometric analogue is the distortion of the (Riemannian) metric of a (Lie) subgroup when compared to the metric induced from the ambient (Lie) group. The distortion of a subalgebra in an algebra also reflects the degree of complexity of the membership problem for the elements of this algebra in this subalgebra. One of our goals here is to investigate, mostly in the case of associative or Lie algebras, if a tame filtration of an algebra can be induced from the degree filtration of a larger algebra.
\end{abstract}

\begin{keyword}
\end{keyword}
\end{frontmatter}

\tableofcontents
\addtocounter{section}{-1}
\section{Introduction}\label{sDF}

Let $A$ be a linear algebra over a field $F$. An ascending filtration $\alpha=\{ A_n\}$ on $A$ is  a sequence of subspaces $A_0\subset A_1\subset\ld\subset A_n\subset\ld$ such that $A=\cup_{n=0}^\infty A_n$ and $A_kA_l\subset A_{k+l}$ for all $k,l=0,1,2,\ld$ Given $a\in A$, the $\alpha$-\emph{degree} of $a$, denoted by $\deg_\alpha a$, is defined as the least $n$ such that $a\in A_n$. If $B$ is a subalgebra of $A$ with a filtration $\alpha$, as above, then we call the filtration $\beta=\{ B\cap A_n\}$ the \emph{restriction} of $\alpha$ to $B$ and we write $\beta=\alpha\cap B$. In the case of monoids (that is, semigroups with 1) the terms of filtrations are simply subsets, with all other conditions being the same. In the case of groups the terms of filtrations must be closed under inverses. 

In this paper we will consider \emph{finitary} filtrations, that is, filtrations all of whose terms are finite-dimensional. If we deal with monoids or groups, the terms of finitary filtrations are simply finite.

A generic example is as follows. Let $A$ be a unital associative algebra generated by a finite set $X$. Then a finitary filtration $\alpha=\{ A_n\}$ arises on $A$, if one sets $A_0=\Sp{1}$ and $A_{n}=A_{n-1}+\Sp{X^n}$, for any $n>1$. The $\alpha$-degree of $a\in A$ in this case is an ``ordinary'' degree with respect to the generating set $X$, that is the least degree of a polynomial in $X$ equal $a$. We write $\deg_\alpha a=\deg_Xa$. Such filtration $\alpha$ is called the \emph{degree filtration} defined by the generating set $X$. 

If $B$ is a subalgebra of an algebra $A$ such that $A$ is generated by a finite set $X$ and $B$ by a finite set $Y$, $\alpha=\{ A_n\}$ and $\beta=\{ B_n\}$ are respective degree filtrations then there is $t$ such that $Y\subset A_t$ and then it immediately follows that for any $n\ge 0$ we  have $B_n\subset A_{tn}$. If $\beta'=\alpha\cap B$ then $\beta$ is majorated by $\beta'$ in the sense of the following definition.

\begin{definition}\label{dEF}
Given two filtrations $\beta=\{ B_n\}$ and $\beta'=\{ B_n'\}$ on the same algebra $B$, we say that $\beta$ is \emph{majorated} by $\beta'$ if there is an integer $t>0$ such that $B_n\subset B'_{tn}$, for all $n\ge 0$. We then write $\beta\preceq\beta'$. If $\beta\preceq\beta'$ and $\beta'\preceq\beta$ then we say that $\beta$ and $\beta'$ are \emph{equivalent} and write $\beta\sim\beta'$.
\end{definition}

Setting $B=A$ in the argument just before the definition, we can say that \emph{any two degree filtration on the same algebra are equivalent}.

In terms of degrees, $\beta\preceq\beta'$ if and only if, there is natural $t$ such that $\deg_{\beta'} b\le t\deg_{\beta} b$ and $\beta\sim\beta'$ if both $\deg_{\beta'} b\le t\deg_{\beta} b$ and $\deg_{\beta} b\le t\deg_{\beta'} b$, for any $b\in B$. If $B$ is a subalgebra of $A$, $\beta$ defined by a finite generating set $Y$ of $B$, $\alpha$ by a finite generating set $X$ of $a$, $\beta'=\beta\cap B$ then inequality  $\beta\preceq\beta'$ is equivalent to the existence of $t\in\N$ such that $\deg_{X} b\le t\deg_{Y} b$, for any $b\in B$.

In a similar way one can speak about the degree filtrations and their equivalence in other classes of algebras, say, Lie, Jordan, etc. 

In the case of monoids  the degree filtration on a monoid $M$ with generating set $X$ is defined by setting $M_0=\{ 1\}$ and, for $n>1$, $M_n=\{ x_1\cdots x_m\,|\, x_i\in X,\, m\le n\}$. But if $G$ is a group with generating set $X$, then one sets $G_0=\{ 1\}$ and $G_n=\{ y_1\cdots y_m\,|\, y_i\in X\cup X^{-1},\, m\le n\}$. In the case of general universal algebras see Subsection \ref{sDA}.

The rate of growth of a degree filtration defined by a finite generating set on an algebra with a finite set of operations or a linear algebra with a finite set of multilinear operations is bounded from above by an exponential function. In other words, for such filtration $\{ A_n\}$,

\medskip

(D) \emph{there is $c>0$ such that $\dim A_{n}<c^n$, for all $n=1,2,\ld$}

\medskip

For example, in the case of associative algebras, if $\# X=d<\infty$ then $\dim A_{n}\le 1+d+...+d^n   \le (d+1)^n$. The same estimate works in the case of monoids and similar estimates are true for groups and linear algebras with one binary operation, such as Lie, Jordan, alternative. etc. For general case see Subsection \ref{sDA}.
These observations lead us to the main definition of the paper. 

\begin{definition}\label{dDFil}
A filtration $\alpha=\{ A_n\}$ on an algebra $A$ satisfying \emph{(D)} is called a \emph{tame filtration}. 
\end{definition}

As noted above, any degree filtration is a tame filtration. Moreover, if $B$ is a subalgebra of an algebra $A$ and $\alpha$ is a tame filtration of $A$ then the restriction $\alpha\cap B$ is a tame filtration of $B$. In particular, \emph{the restriction of a degree filtration of an algebra $A$ to a subalgebra $B$ is always a tame filtration of $B$}. This simple fact will be used without further comments in several theorems of this paper designed to deal with the following natural questions.

\begin{enumerate}
\item[(1)] \emph{Is it true that every tame filtration of an algebra $B$ is equivalent (or equal) to a filtration restricted from the degree filtration of a \fg algebra $A$ where $B$ is embedded as a subalgebra}?
\item[(2)] \emph{If the answer to the previous question is ``yes'', can one choose $A$ finitely presented? If not, indicate conditions ensuring that the answer is still ``yes''.}
\end{enumerate}

It follows from our examples in Section \ref{sDA}, that the set of pairwise inequivalent tame filtrations on any  infinite-dimensional algebra $B$ is uncountable. At the same time, the set of finitely generated subalgebras of finitely presented algebras over a fixed countable field is countable. Thus, the answer to the second question is generally ``no''. (Similar argument works in the case of infinite groups, monoids, etc.) However, there is an extensive class of so called \emph{constructive} tame filtrations, see definition in Section \ref{sHEUA}, for which the answer is ``yes''.

Now our answers to the above questions are as follows. Note that an algebra $A$ is called \emph{unital} if $A$ has identity element $1$; a subalgebra $B$ of such $A$ is called \emph{unital} if $1\in B$.

\emph{\begin{enumerate}
\item[$\mathrm{(i)}$] In the case of associative and Lie algebras,  we prove that a filtration $\beta$ on an algebra $B$  is  a tame filtration if and only if $\beta\sim\alpha\cap B$ where $\alpha$ is a degree filtration on a \emph{2-generator} algebra $A$, in which $B$ is embedded as a subalgebra. 
\item[$\mathrm{(ii)}$] Under the same condition as in \emph{(i)}, we prove that a filtration $\beta$ on an algebra $B$  is  a tame filtration if and only if $\beta=\alpha\cap B$ where $\alpha$ is a degree filtration on a \emph{finitely generated} algebra $A$, in which $B$ is embedded as a subalgebra. 
\item[$\mathrm{(iii)}$] In the case of associative algebras over any field, which is finitely generated over prime subfield, a constructive filtration $\beta$ on a (unital) countable algebra $B$ is a tame filtration if and only if $\beta\sim\alpha\cap B$ where $\alpha$ is a degree filtration on a (unital) finitely presented algebra $A$, in which $B$ is embedded as a (unital) subalgebra. 
\end{enumerate}}

The results of (i) and (iii) have their analogues in Group Theory, see \cite{AO} and \cite{OSD}, where the author answers questions asked by M. Gromov \cite{Gr}.

If $B$ is a \fg subalgebra of a \fg algebra $A$ then we say then $B$ is  embedded in $A$ \emph{without distortion} (or that $B$ is an \emph{undistorted subalgebra} of $A$) if a degree filtration of $B$ is equivalent to the restriction to $B$ of a degree filtration of $A$. As will be shown in Section \ref{sBP}, then any degree filtration of $B$ is equivalent to the restriction of any degree filtration of $A$. Thus the property of being undistorted does not depend on the particular choice of finite generating sets in $A$ and $B$. All this is obviously true in the case of monoids or groups.

Some results about distorted and undistorted subalgebras are as follows.

\emph{\begin{enumerate}
\item[$\mathrm{(iv)}$] In the case of \emph{commutative} associative algebras, all \fg subalgebras of \fg algebras are undistorted.
\item[$\mathrm{(v)}$] In the case of \emph{noncommutative} associative algebras, free associative algebras of finite rank and also free algebras  of finite rank in certain proper varieties of associative algebras have distorted \fg subalgebras. 
\item[$\mathrm{(vi)}$]All \fg subalgebras in free \emph{Lie algebras} of finite rank are undistorted. 
\item[$\mathrm{(vii)}$] Any finitely generated associative, respectively, Lie algebra can be embedded without distortion in a simple 2-generator associative, respectively, Lie algebra.
\end{enumerate}}

An important tool while comparing two finitary filtrations on the same algebra is given by the following.

\begin{definition}\label{ddf}
Let $\alpha=\{ A_n\}$ and $\alpha'=\{ A'_n\}$ be two finitary filtrations on the same algebra $A$. For each $n\ge 0$ we set

\[
\dist_\alpha^{\alpha'}(n)=\min\{m\,|\, A_n\subset A'_{m}\}.
\]

The function $\dist_\alpha^{\alpha'}$ thus obtained is called the \emph{distortion function of $\alpha'$ with respect to $\alpha$.}
\end{definition}

In terms of distortion functions, two finitary filtrations $\alpha$ and $\alpha'$ of an algebra $A$ are equivalent if and only if both $\dist_\alpha^{\alpha'}$ and $\dist_{\alpha'}^\alpha$ are bounded from above by a linear function. 
The notion of undistorted embedding in terms of distortion functions (see Proposition \ref{BOr1}) is equivalent to the following: given two degree filtrations, $\alpha$ in $A$ and $\beta$ in $B$, the function $\dist_{ \alpha\cap B}^{\beta}$ is bounded from above by a linear function. Again, all the above is true in the case of monoids, groups or even more general universal algebras.

Also the following result will be proven in a far greater generality than we state here (see Subsection \ref{sDA}).
\emph{\begin{enumerate}
\item[$\mathrm{(vii)}$] Let $g:\N\to\N$ be any function. On any \fg infinite-dimensional associative or Lie algebra over a field there exist tame filtrations $\alpha$ and $\alpha'$ such that $\dist_\alpha^{\alpha'}(n)>g(n)$, for any $n\ge 1$.
\end{enumerate}}

Given two functions $f,g:\N\to\N$ we say that $f$ is \emph{majorated} by $g$ and write $f\preceq g$ if there is natural $t$ such that  $f(n)\le tg(tn)$, for all $n\in\N$. If both $f\preceq g$ and $g\preceq f$ then we say that $f$ is \emph{equivalent} to $g$ and write $f\sim g$. The equivalence classes acquire partially ordering if one sets $[f]\le[g]$ once $f\preceq g$.  

Now suppose we are given two pairs of equivalent filtrations: $\alpha\sim\beta$ and $\alpha'\sim\beta'$. Then it is easy to show (see below Claim (5) of Proposition \ref{pGPDF}) that $\dist_{\alpha}^{\alpha'}\sim\dist_{\beta}^{\beta'}$. This allows us, given an embedding $\vp:B\to A$, where $A$ and $B$ are \fg algebras, to define the distortion $\dist(\vp)$ of $\vp$ as follows. Take any two degree filtrations $\alpha$ in $A$, $\beta$ in $B$ and consider the image $\vp(\beta)$ under $\vp$. Then $\dist(\vp)$ is defined as the equivalence class of $\dist_{\alpha\cap\vp(B)}^{\vp(\beta)}$ under the above equivalence relation. In particular, if $B$ is a subalgebra of $A$, we define $\dist_A^B=\dist(\vp)$ where $\vp$ is the natural embedding of $B$ in $A$. The definition of this asymptotic invariant is analogous to the definition introduced by Gromov \cite{Gr} in order to measure the distortion of the lengths of elements of groups that may occur when one group is embedded in another.

In our simple Proposition \ref{BOpGA} we note that $\dist_G^H=\dist_{F[G]}^{F[H]}$, where $H$ is a \fg subgroup of a \fg group $G$, and $F[H]$, $F[G]$ are the group algebras of $H$ and $G$, with natural embedding of $F[H]$ in $F[G]$. Similar is the situation in the case of semigroups. In the case of Lie algebras and their universal enveloping algebras the situation is somewhat different. For a subalgebra $M$ of a Lie algebra $L$, let $U(M)$ be the associative subalgebra generated by $M$ in the universal enveloping algebra $U(L)$ for $L$. It is well-known that $U(M)$ is isomorphic to the universal enveloping algebra for $M$. Given a function $f:\N\to\N$ we define its \emph{superadditive closure} as the least function $\bf:\N\to\N$ such that $f(n)\le\bf(n)$, for all $n\in\N$, and satisfying $\bf(m)+\bf(n)\le \bf(m+n)$.

\emph{\begin{enumerate}
\item[$\mathrm{(ix)}$] If $M$ is a \fg subalgebra in a \fg Lie algebra $L$ then $\dist_L^M\le \dist_{U(L)}^{U(M)}$. More precisely, if $f$ represents $\dist_L^M$ then $\dist_{U(L)}^{U(M)}$ is represented by $\bf$. There are examples of pairs $M\subset L$ for which the above inequality is strict.
\end{enumerate}}

As we mentioned earlier, any filtration $\alpha$ on an algebra $A$ allows one to define the degree function $\deg_{\alpha}:A\to\N$. If $\alpha$ is a tame filtration, we call $\deg_{\alpha}$ \emph{tame}. Tame degrees have all ``natural'' properties of ordinary degrees  with respect to a finite generating set of $A$. Still, in Subsection \ref{ssEDF} we show that a particular tame degree can be a pretty wild function, even in the case of very basic associative or Lie algebras.

\section{General properties of distortion}\label{sGPD}

In this section we consider tame filtrations on universal algebras. In Subsection \ref{BOfua} we discuss filtrations on universal algebras. Using this generality, in Subsection \ref{sBP} we establish some basic facts about distortion functions and embeddings with or without distortion. In Subsection \ref{sOPUE} we discuss connections with some algorithmic decidability problems. Finally, in Subsection \ref{sDA} we provide with a construction that exhibits uncountably many pairwise inequivalent tame filtrations on virtually any infinite (-dimensional) algebra.

\subsection{Filtrations on universal algebras}\label{BOfua}

Given a signature $\Omega$, we will be considering general universal algebras with operations in $\Omega$, we called them $\Omega$-algebras, and linear $\Omega$-algebras, which are vector spaces over a field $F$, with multilinear operations in $\Omega$. 

Let us indicate some problems which arise when we try to generalize the notion of a filtration $\{ A_n\}$, $A_0\subset A_1\subset\ldots,$ on a universal algebra $A$, which may have nullary or unary operations. Suppose, following the pattern exhibited in Introduction, we require that $A_{i_1}\ld A_{i_n} \omega \subset A_{i_1+...+i_n}$, where $\omega$ is an arbitrary operation of arbitrary arity $n$. Then the condition $A_1\lambda\subset A_1$, for a unary operation $\lambda$ can quickly make $A_1$ infinite (or infinite-dimensional) even if $A$ is finitely-generated and $\lambda$ is the only operation on $A$. Another kind of confusion may arise if we require $A_0\ld A_0\omega \subset A_0$, for any $n$-ary operation $\omega$. These problems can be avoided if we give the following.

\begin{definition}\label{dBOunialg}
If $A$ is an $\Omega$-algebra, respectively a linear $\Omega$-algebra, then an ascending chain of subsets, respectively, subspaces,  $\alpha=\{ A_n\}$, $n=0,1,2,\ld$, is a \emph{filtration} provided that $A=\cup_{n=0}^\infty A_n$ and for any $m$-ary operation $\omega$ and any $a_1\in A_{i_1},\ld,a_m\in A_{i_m}$ one has $a_1\ld a_m\omega\in A_{i_1+\ld+i_m+1}$. 
\end{definition}

As we can see, in this approach we count the number of operations rather than the number of elements involved in these operations. In the case of rings, semigroups or groups (in the latter case each term $T$ of a filtration is usually assumed symmetric: $T^{-1}=T$), a filtration $\{ A_n\}$ in the sense of Definition \ref{dBOunialg} produces a filtration $\{ B_n\}$ in the ``common'' sense if one sets $B_{n+1}=A_n$, for all $n=0,1,\ld$. Two filtrations are equivalent in the ``new'' sense if and only if they are equivalent in the ``old'' sense. Speaking about the notions of the equivalence of two filtrations or the distortion functions, they are the same as given in Definitions \ref{dEF} and \ref{ddf}. The notion of the tame filtration remains the same as in Introduction, except that in the case where we do not have the structure of vector space on $A$ (like in groups, monoids, etc.) we simply have to require the existence of an integer $c$ such that $\# A_n\le c^n$, for all $n \ge 1$. Here and in what follows, given a set $M$, we denote the cardinality of $M$ by $\# M$.

Let us now define the \emph{degree filtration} of an $\Omega$-algebra, with \emph{finite signature} $\Omega$ as follows. First we write $\Omega=\Omega_0\cup\Omega_1\ld\cup\Omega_q$ where the operations in $\Omega_m$ have arity $m$. We identify $\Omega_0$ with the subset of elements of $A$ selected by these operations. If $X$ is a generating set for $A$ then the degree filtration $\{ A_n\}$, $n=0,1,2,\ld$, is defined as follows: $A_0=X\cup\Omega_0$ and if $n>0$ then
 
\[
A_n = \bigcup_{m = 1}^q \bigcup_{\omega\in\Omega_m} \bigcup_{i_1+...+i_m+1\le n} A_{i_1}...A_{i_m}\omega.
\]

In the case of linear $\Omega$-algebras one has to define $A_0$ as the linear span of $X\cup\Omega_0$ and replace all unions in the previous formula by vector space sums. It follows from this definition that the degree filtration is a filtration and if $X$ is finite then all sets $A_i$ are finite (-dimensional). Any two degree filtrations $\alpha=\{ A_n\}$ and $\alpha'=\{ A_n'\}$ on a finitely generated algebra $A$ are equivalent. Indeed, if $\alpha$, respectively, $\alpha'$ is defined by a finite generating set $X$, respectively, $X'$ then there is $k$ such that $X'\subset A_k$. Let us set $t=kq+1$. Then, using induction, for each $\omega\in\Omega_m$, $m\le q$, and $n=i_1+\ld+i_m+1$, we will have 

\[
A^{\prime}_{i_1}\ld A^{\prime}_{i_m}\omega \subset A_{\max\{ti_1,k\}}\ld A_{\max\{ti_m, k\}}\omega
\subset A_{t(i_1+\cdots+i_m) +qk+1}\subset A_{tn}.
\]

\subsection{Some basic facts about distortion}\label{sBP} We start by listing few general properties of distortion functions. We denote by $\mathrm{id}$ the identity function on $N$ and $f\cir g$ the usual composition of two functions $f,g:\N\to\N$. Given $t\in\N$ and $f:\N\to\N$ we denote by $t\cdot f$ the scalar multiple of $f$ by $t$. We write $f\le g$ if $f(n)\le g(n)$, for all $n\in\N$. Also, given a homomorphism of algebras $\vp:A\to B$ and a filtration $\alpha=\{ A_n\}$ on $A$, we denote by $\vp(\alpha)$ the \emph{image-filtration} $\{\vp(A_n)\}$ on $\vp(A)$.

\begin{proposition}\label{pGPDF} The following facts are true for the distortion functions of filtrations.
\begin{enumerate}
\item[\emph{(1)}] If $\mu$, $\nu$, $\pi$ are three filtrations on the same algebra $A$ then 

\bee{etrans}
\ds{\mu}{\pi}\le\ds{\nu}{\pi}\cir\ds{\mu}{\nu};
\ene

\item[\emph{(2)}] If $\mu$, $\nu$ are two filtrations on an algebra $A$ and $B$ is a subalgebra of $A$ then
 
\[
\ds{\mu\cap B}{\nu\cap B}\le \ds{\mu}{\nu};
\]

\item[\emph{(3)}] If $C\subset B\subset A$ are algebras with filtrations $\alpha$ on $A$, $\beta$ on $B$ and $\gamma$ on $C$ then
 
\[
\ds{\alpha\cap C}{\gamma}\le \ds{\beta\cap C}{\gamma}\cir\ds{\alpha\cap B}{\beta};
\]

\item[\emph{(4)}] If $\mu$, $\nu$ are two filtrations on an algebra $A$, and $\vp:A\to B$ is a homomorphism of algebras then 

\[
\ds{\vp(\mu)}{\vp(\nu)}\le \ds{\mu}{\nu}.
\]

\item[\emph{(5)}] If $\mu\sim\mu'$ and $\nu\sim\nu'$ are two pairs of pairwise equivalent filtrations on an algebra $A$ then $\dist_{\mu}^{\nu}\sim\dist_{\mu'}^{\nu'}$.

\end{enumerate}
\end{proposition}

\pp (1) Let $\mu=\{ U_n\}$, $\nu=\{ V_n\}$, $\pi=\{ W_n\}$, $f=\ds{\mu}{\pi}$, $g=\ds{\nu}{\pi}$ and $h=\ds{\mu}{\nu}$. Then, by definition of the distortion function, for any $n$ we have $U_n\subset V_{h(n)}\subset W_{g(h(n))}$. Since $f(n)=\min\{ m\,|\,U_n\subset W_m\}$, it follows that $f(n)\le g(h(n))$, proving our first claim.

(2) Let $\mu=\{ U_n\}$, $\nu=\{ V_n\}$, $g=\ds{\mu\cap B}{\nu\cap B}$ and $h=\ds{\mu}{\nu}$. Then, for any $n$, we have $U_n\subset V_{h(n)}$. Taking intersections of both sides of the latter containment with $B$, we obtain  $U_n\cap B\subset V_{h(n)}\cap B$. By definition of the distortion function, $g(n)\le h(n)$, proving our second claim.

(3) This follows by applying (1) and (2). Indeed, let us apply (1) to the following three filtrations on $C$: $\mu=\alpha\cap C$, $\nu=\beta\cap C$  and $\pi=\gamma$. We will have $\ds{\alpha\cap C}{\gamma}\le \ds{\beta\cap C}{\gamma}\cir\ds{\alpha\cap C}{\beta\cap C}$. Now one can rewrite $\ds{\alpha\cap C}{\beta\cap C}=\ds{(\alpha\cap B)\cap C}{\beta\cap C}$. According to (2), $\ds{(\alpha\cap B)\cap C}{\beta\cap C}\le\ds{\alpha\cap B}{\beta}$. Since the distortion functions are non-decreasing, it follows that $\ds{\alpha\cap C}{\gamma}\le \ds{\beta\cap C}{\gamma}\cir\ds{\alpha\cap C}{\beta\cap C}\le\ds{\beta\cap C}{\gamma}\cir\ds{\alpha\cap B}{\beta}$, as claimed.

(4) Let $\mu=\{ U_n\}$, $\nu=\{ V_n\}$, $g=\ds{\vp(\mu)}{\vp(\nu)}$ and $h=\ds{\mu}{\nu}$. Then, for any $n$, we have $U_n\subset V_{h(n)}$. Taking the images of the latter containment under $\vp$, we obtain $\vp(U_n)\subset \vp(V_{h(n)})$. It follows that $g(n)\le h(n)$.

(5) Using (1), we can write 

\[
\dist_{\alpha}^{\beta}\le\dist_{\beta'}^{\beta}\cir\dist_{\alpha'}^{\beta'}\cir\dist_{\alpha}^{\alpha'}.
\]

We may assume that there is $t$ such that both $\dist_{\beta'}^{\beta}\le t\cdot\mathrm{id}$ and $\dist_{\alpha}^{\alpha'}\le t\cdot\mathrm{id}$. It then follows that $\dist_{\alpha}^{\beta}(n)\le t\cdot\dist_{\alpha'}^{\beta'}(tn)$, for all $n\in\N$ and so $\dist_{\alpha}^{\beta}\preceq \dist_{\alpha'}^{\beta'}$. The converse inequality follows in the same manner.
\epf

The just proved proposition  turns out to be instrumental in proving general properties of distortion in algebras. 

\begin{proposition}\label{BOr1}
Let $A$ be an algebra, $B$ and $C$ subalgebras of $A$ such that $C\subset B$. 
\begin{enumerate}
\item[\emph{(i)}] If $A$ is a finitely generated algebra then any two degree filtrations on $A$ are equivalent.
\item[\emph{(ii)}] If $\alpha$ and $\alpha'$ are two filtrations of $A$ and $\alpha\sim \alpha'$ then $\alpha\cap B\sim \alpha'\cap B$. 
\item[\emph{(iii)}] If $A$ is finitely generated then all tame filtrations of $B$ obtained by restriction from the degree filtrations of $A$ are pairwise equivalent. 
\item[\emph{(iv)}] Suppose that $A$ and $B$ are finitely generated and $\alpha$ and $\beta$ are two degree filtrations on $A$ and $B$, respectively. Then $\beta\preceq\alpha\cap B$. Thus $B$ is undistorted in $A$ if and only if there are two degree filtrations $\alpha$ in $A$ and $\beta$ in $B$ such that $\alpha\cap B\preceq\beta$.
\item[\emph{(v)}] Let all $A$, $B$ and $C$ be finitely generated. If $C$ is undistorted in $B$ and $B$ undistorted in $A$, then $C$ is undistorted in $A$.
\item[\emph{(vi)}] Let $A$, $B$, $C$ be as in \emph{(v)}. If $C$ is undistorted in $A$ then $C$ is undistorted in $B$.
\item[\emph{(vii)}] Let $\vp$ be a homomorphism of $A$ which is injective when restricted to $B$. If $A$ and $B$ are \fg and $\vp(B)$ is undistorted in $\vp(A)$ then $B$ is undistorted in $A$. In particular, a retract of an algebra is always an undistorted subalgebra.
\item[\emph{(viii)}] Suppose $A\in\mathfrak{Q}$, where $\mathfrak{Q}$ is a (quasi)variety $\mathfrak{Q}$ of algebras, $\vp:B\to A$ the natural embedding of $B$ in $A$, $\mu:F(X)\to A$, $\nu:F(Y)\to B$ epimorphisms of free algebras of $\mathfrak{Q}$, generated by finite sets $X$ and $Y$. Suppose $B$ is given by some set  $T(Y)\subset F(Y)\times F(Y)$ of defining relations. Let $\pi:F(Y)\to F(X)$ be any homomorphism such that $\mu\cir\pi=\vp\cir\nu$ and $T(X)=(\pi\times\pi)(T(Y))$.  If $B$ is undistorted in $A$ then $B$ is an undistorted subalgebra in an algebra $A'$ defined in terms of generators $X$ and relations $T(X)$. 
\end{enumerate}
\end{proposition}

\pp (i) The proof of this was given in Subsection \ref{BOfua}.

(ii) By definition, $\alpha\preceq\alpha'$ if and only if there is a linear function $f$ such that $\ds{\alpha}{\alpha'}\le f$. Applying Claim (2) of Proposition \ref{pGPDF}, we obtain $\ds{\alpha\cap B}{\alpha'\cap B}\le\ds{\alpha}{\alpha'}\le f$. Thus $\alpha\cap B\preceq\alpha'\cap B$. By symmetry, $\alpha'\cap B\preceq\alpha\cap B$, and our claim follows.

(iii) Follows from (i) and (ii).

(iv) Set $\gamma=\alpha\cap B$.  Choose a generating set $X$ for $A$ which includes the generating set $Y$ for $B$. Let $\alpha'$ and $\beta'$ be respective degree filtrations and $\gamma'=B\cap\alpha'$. Then obviously $\ds{\beta'}{\gamma'}\le \mathrm{id}$. Now $\alpha\sim\alpha'$, $\beta\sim\beta'$ as pairs of degree filtrations on the same algebras and 
$\gamma\sim\gamma'$ by (ii). Thus there is integral $t$ such that $\ds{\beta}{\beta'},\,\ds{\gamma'}{\gamma}\le t\cdot\mathrm{id}$. Applying Claim (1) of Proposition \ref{pGPDF}, we obtain

\[
\ds{\beta}{\gamma}\le\ds{\gamma'}{\gamma}\cir\ds{\beta'}{\gamma'}\cir\ds{\beta}{\beta'}\le (t\cdot\mathrm{id})\cir \mathrm{id}\cir (t\cdot\mathrm{id})= t^2\cdot\mathrm{id}.
\] 

So $\beta\preceq\gamma=\alpha\cap B$, as claimed. Thus if we only assume $\alpha\cap B\preceq\beta$, we will have  $\alpha\cap B\sim\beta$, proving that $B$ is an undistorted subalgebra in $A$. Now if $B$ is undistorted in $A$ then $\alpha\cap B\sim\beta$, and in particular, $\beta\preceq\alpha\cap B$.

(v) If $\alpha$, $\beta$ and $\gamma$ are degree filtrations in $A$, $B$ and $C$, respectively, then by previous claim we only need to establish that $\ds{\alpha\cap C}{\gamma}\le u\cdot\mathrm{id}$, for some integral $u$. By our assumption, we have such constant $t$ for both $\ds{\beta\cap C}{\gamma}$ and $\ds{\alpha\cap B}{\beta}$. Using Claim (3) of Proposition \ref{pGPDF}, we obtain
 
\[
\ds{\alpha\cap C}{\gamma}\le\ds{\beta\cap C}{\gamma}\cir\ds{\alpha\cap B}{\beta}\le (t\cdot\mathrm{id})\cir (t\cdot\mathrm{id})=t^2\cdot\mathrm{id}.
\]

\noindent Selecting $u=t^2$ completes the proof of this claim.

(vi) We apply Claim (1) of Proposition \ref{pGPDF} to filtrations $\mu=\beta\cap C$, $\nu=\alpha\cap C$ and $\pi=\gamma$ on $C$. We will then have $\ds{\beta\cap C}{\gamma}\le\ds{\alpha\cap C}{\gamma}\cir\ds{\beta\cap C}{\alpha\cap C}$. Using Claim (2) of Proposition \ref{pGPDF}, we obtain $\ds{\beta\cap C}{\alpha\cap C}\le\ds{\beta}{\alpha\cap B}$. By (iv) there is $t\in\N$ such that $\ds{\beta}{\alpha\cap B}\le t\cdot\mathrm{id}$. Since $C$ is undistorted in $A$, we may assume that also $\ds{\alpha\cap C}{\gamma}\le t\cdot\mathrm{id}$. Since the distortion functions are non-decreasing, we finally get $\ds{\beta\cap C}{\gamma}\le (t\cdot\mathrm{id})\cir (t\cdot\mathrm{id})=t^2\cdot\mathrm{id}$. Thus $\beta\cap C\preceq \gamma$ and by (iv), $C$ is undistorted in $B$, as claimed.

(vii) Let $\alpha=\{ A_n\}$, $\beta=\{ B_n\}$ be degree filtrations in $A$, respectively, $B$. Then $\vp(\alpha)$ and $\vp(\beta)$ are degree filtrations in $\vp(A)$, respectively, $\vp(B)$. By (iv) we only need to show $\alpha\cap B\preceq \beta$. Since $\vp(\alpha)\cap\vp(B)\preceq \vp(\beta)$, there is $t\in\N$, such that for any $n\in\N$ and any $b\in B\cap A_n$ we have $\vp(b)\in\vp( B_{tn})$. Since the restriction of $\vp$ to $B$ is injective, we have $b\in B_{tn}$. Thus $B\cap A_n\subset B_{tn}$ and so $B$ is undistorted in $A$. 

(viii) From $\mu\cir\pi=\vp\cir \nu$ it follows that 

\begin{eqnarray*}
(\mu\times \mu)(T(X))&=&(\mu\times \mu)(\pi\times \pi ) (T(Y))=(\vp\times \vp) (\nu\times \nu)(T(Y))\\
&\subset & (\vp\times \vp)(\mathrm{diag}(B\times B))\subset \mathrm{diag}(A\times A).
\end{eqnarray*}

Thus $A$ satisfies all relations $T(X)$ and so if $A'$ is defined by generators $X$ and relations $T(X)$ and $\mu':F(X)\to A'$ is the respective canonical epimorphism then there is a natural epimorphism $\ve:A'\to A$ such that $\mu=\ve\cir\mu'$. Finally setting $\vp'(b)=\mu'\cir\pi\cir\nu^{-1}(b)$ correctly defines the desired undistorted embedding of $B$ in $A'$. Indeed, if $f_1,f_2\in F(Y)$ are such that $\nu(f_1)=\nu(f_2)$ then $f_1=f_2$ is a relation in $B$ hence $\pi(f_1)=\pi(f_2)$ is a relation in $A'$ and so $(\mu'\cir\pi)(f_1)=(\mu'\cir\pi)(f_2)$. This shows that $\vp':B\to A'$ is well defined. Now, for any $b\in B$,
 
\[
\ve\cir\vp'(b)=\ve\cir\mu'\cir\pi\cir\nu^{-1}(b)=\mu\cir\pi\cir\nu^{-1}(b)=\vp\cir\nu\cir\nu^{-1}(b)=\vp(b).
\]

Thus, $\vp=\ve\cir\vp'$ and it follows from the injectivity of $\vp$ that $\vp'$ is also injective.

\medskip

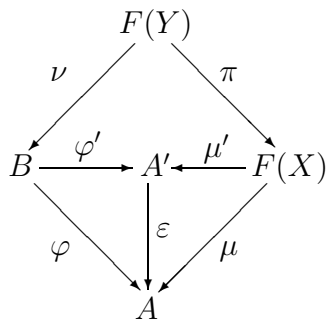
\begin{picture}(0,170)
\put(108,92){\vector(1,-1){40}}
\put(110,98){\vector(1,0){35}}
\put(188,98){\vector(-1,0){28}}
\put(196,92){\vector(-1,-1){41}}
\put(98,95){$B$}
\put(148,95){$A'$}
\put(154,72){$\ve$}
\put(151,92){\vector(0,-1){40}}
\put(190,95){$F(X)$}
\put(146,41){$A$}
\put(140,149){$F(Y)$}
\put(160,145){\vector(1,-1){39}}
\put(146,145){\vector(-1,-1){40}}
\put(114,65){$\vp$}
\put(178,65){$\mu$}
\put(123,103){$\vp'$}
\put(172,102){$\mu'$}
\put(114,131){$\nu$}
\put(178,131){$\pi$}
\end{picture}

\noindent Finally, using our previous Claim (vii), we conclude that $\vp'$ is an undistorted embedding, provided that $\vp$ is undistorted. The proof is now complete.\epf


If $G$ is an $\Omega$-algebra then the linear space $A=F[G]$ with basis $G$  is a linear $\Omega$-algebra. Given a filtration $\alpha=\{ G_n\}$ on $G$, we denote by $F[\alpha]$ the filtration on $A$ whose terms $A_n$ are linear spans of the respective terms of $\alpha$. If $\alpha$ is a finitary filtration then also $F[\alpha]$ is finitary. If $\alpha$ is the degree filtration of $G$ defined by a finite set $X$ then $F[\alpha]$ is the degree filtration of $A$ defined by the same set $X$. Finally, $\dim F[G]_n=\# G_n$. Now suppose $H$ is a subalgebra of $G$. Then the linear span $F[H]$ is a subalgebra of $F[G]$ and $F[G]_n\cap F[H]=\Sp{G_n\cap H}$. It now follows that $\ds{\alpha\cap H}{\beta}=\ds{F[\alpha]\cap F[H]}{F[\beta]}$. This argument works, in particular, in the case of semigroups and semigroup algebras. In the case of groups and group algebras, one can restrict oneself to degree filtrations $\alpha$ defined by symmetric generating sets $X$, that is $X^{-1}=X$.  

\begin{proposition}\label{BOpGA} Let $H$ be a finitely generated subgroup of a finitely generated group $G$, $F[G]$ the group algebra of $G$ over a field $F$, $F[H]$ the group algebra of $H$ naturally embedded in $F[G]$. Then the distortion of the embedding of $H$ in $G$ is the same as the distortion of the embedding of $F[H]$ in $F[G]$. The same claim holds valid in the case of semigroups/monoids and their semigroup/monoid algebras.\epf
\end{proposition}

This latter result leads to a variety of distortions in associative algebras obtained with the help of known results on distortion in groups (see \cite{Gr} in the case of Examples (i) and (ii) and \cite{AP} in the case of Example (iii)). We exhibit here just three.

\begin{example}\label{BOex1} \emph{Let $H$ be cyclic subgroup with generator $c$ in the Heisenberg group $G$ with presentation $G = \langle a,b,c\,|\,[a,b]=c, [a,c]=[b,c]=1\rangle$. Then the distortion of $F[H]$ in $F[G]$ is quadratic.} 
\end{example}
\begin{example}\label{BOex2} \emph{Let $H$ be cyclic subgroup with generator $b$ in the Baumslag - Solitar group $G$ with presentation $G = \langle a,b\,|\,aba^{-1}=b^2\rangle$. Then the distortion of $F[H]$ in $F[G]$ is exponential.}
\end{example}  
\begin{example}\label{BOex3} \emph{Let $g(k)$ be the ``exponential tower with basis $2$ of height $k$'', that is, the function defined by $g(1)=2$ and $g(k)=2^{g(k-1)}$, for $k>1$. Let $H$ be the  subgroup generated by $b$ in the 1-relator Baumslag group $G$ with presentation $G = \langle a,b\,|\,(aba^{-1})b(aba^{-1})^{-1} = b^2\rangle$. Then the distortion of $F[H]$ in $F[G]$ is represented by a function $f$ such that $f(n)=g([\log_2n])$.} 
\end{example}

\subsection{Distortion and Membership Problem}\label{sOPUE}

In this subsection we elaborate on the following observation: knowing that the distortion of a subalgebra $B$ in an algebra $A$ is not too ``ugly'' can help one to build an algorithm deciding whether an element $a\in A$ is actually an element of $B$. Our main result holds valid for any \fg $\Omega$-algebra. However, additional details arise in the case of linear $\Omega$-algebras and so we will start with this kind of algebras. As before, the term ``algebra'' will mean linear ``$\Omega$-algebra''. 

In this subsection we will be considering linear $\Omega$-algebras over constructive fields. A field $F$ is called \emph{constructive} if there is an enumeration of the elements of $F$: $a_1,a_2,\ld$ such that $a_ia_j=a_{f(i,j)}$ and $a_ia_j=a_{g(i,j)}$ where $f,g:\N\to\N$ are recursive functions.

Let $A$ be a linear $\Omega$-algebra over a constructive field $F$, with finite generating set $X$. We say that the Intersection Problem of Finite-dimensional Subspaces (\emph{Intersection Problem}, for short) is decidable in $A$ if there is an algorithm that allows, given two finite subsets $R$ and $S$ of elements of $A$, written as ``polynomials'' in $X$, with linear spans $\Sp{R}$ and $\Sp{S}$ over $F$, to produce a finite subset $T\subset A$, again written as ``polynomials'' in $X$, which is a linear basis of $\Sp{R}\cap\Sp{S}$. An equivalent problem is that of decidability of linear dependence of finite sets. Further, given a \fg subalgebra $B$ of a \fg algebra $A$, we say that the \emph{Generalized  Membership Problem} in $B$ is decidable in $A$ if there is an effective procedure that allows one, for any finite subset $S\subset A$, to effectively produce a finite subset $T\subset A$, such that $T$ is a linear basis of $B\cap\Sp{S}$. Again, $S$ and $T$ are written as ``polynomials'' in $X$ or, more precisely, the elements of the free linear $\Omega$-algebra $\mathcal{F}(X)$ generated by $X$.

\begin{theorem}\label{BOpLDP}
Let $A$ be a finitely generated $\Omega$-algebra where $\Omega$ is a finite set of operations. If $A$ is a linear algebra then we additionally assume that the ground field of coefficients $F$ is constructive and that the Intersection Problem is algorithmically decidable in $A$. Then the following conditions on a finitely generated subalgebra $B$ are equivalent:
\begin{enumerate}
\item[\emph{(1)}] The (Generalized) Membership Problem for $B$ in $A$ is algorithmically decidable;
\item[\emph{(2)}] Any distortion function for $B$ in $A$ is recursive;
\item[\emph{(3)}] A distortion function for $B$ in $A$ is bounded from above by a recursive function.
\end{enumerate}
\end{theorem}

\pp Let $\alpha=\{ A_n\}$ be the degree filtrations  of $A$ and $\beta=\{ B_n\}$ of $B$. We restrict ourselves to the proof in the case of linear $\Omega$-algebras, with few remarks about the ``non-linear'' case at the end of the proof.

To prove Implication (1)$\Rightarrow$(2), assume that $A$ is a linear $\Omega$-algebra and $B$ a subalgebra with Generalized Membership Problem effectively decidable. Given any natural $n$, let us effectively select a finite basis $P_n$ of $A_n$, the $n\th$ term of $\alpha$. Using our hypothesis, we can effectively write the finite set $Q_n$ which is a linear basis of $K_n=B\cap A_n$, the $n\th$ term of the restriction filtration $\alpha\cap B$. After this, starting from $m=n$, we effectively write the linear bases $R_{n,m}$ of the subspaces $K_n\cap B_m$ where each $B_m$ is a \fd space, the $m\th$ term of $\beta$. The least value of $m$ such that $\# R_{n,m}=\# Q_n$ equals the $n\th$ value of the distortion function $\dist_{\alpha\cap B}^{\beta}$. So this function is recursive, proving our first implication.

The implication (2)$\Rightarrow$(3) being trivial, let us prove (3)$\Rightarrow$(1). We assume that $\dist_{\alpha\cap B}^{\beta}$ is bounded from above by a recursive function $f:\N\to\N$, for a fixed choice of degree filtrations $\alpha$ in $A$ and $\beta$ in $B$. Let us prove that the Generalized Membership Problem in $B$ is algorithmically decidable. Indeed, let $S$ be a finite subset of $A$ and $n\in\N$ be a number such that $\Sp{S}\subset A_m$. By our hypothesis, then $B\cap\Sp{S}\subset B_{f(n)}$. Then also $B\cap\Sp{S}\subset B_{f(n)}\cap\Sp{S}$. Since $f(n)$ is computable and $B_{f(n)}$ has an effectively computable basis $Q_{f(n)}$, using the decidability of the Intersection Problem in $A$, we effectively find a basis in $B_{f(n)}\cap\Sp{S}$, hence in $B\cap\Sp{S}$. Thus (3)$\Rightarrow$(1) has been also established, proving our theorem in the case of linear $\Omega$-algebras.
 
 In the case of $\Omega$-algebras which are not linear, the logic of the proof is exactly the same except that we do not need to consider any spans of finite sets of monomials but rather the finite sets themselves, and their intersections, which by necessity are constructive.\epf

As an application, let us consider a (quasi)variety $\mathfrak{Q}$ of linear algebras where there are finitely presented algebras with  algorithmically undecidable equality problem (see \cite{KS}). Additionally we assume that the Intersection Problem is algorithmically decidable in free algebras of finite rank of $\mathfrak{Q}$. As an example, one can consider the varieties of all associative or all Lie algebras, etc.

We will now follow the pattern of a well-known group-theoretic construction of Mikhailova \cite{Mi}, whose relevance to the distortion in groups was first mentioned by Gromov \cite{Gr}.  

\begin{example}\label{BOeMikh} \emph{Let $B$ be a finitely presented linear algebra in a (quasi)variety $\mathfrak{Q}$, with undecidable equality problem. We write $B=\mathcal{F}(X)/I$ where $\mathcal{F}(X)$ is a free algebra with free generating set $X$ and $I$ an ideal in $\mathcal{F}(X)$ generated by a finite set $R$. Let us consider the direct product $\mathcal{G}=\mathcal{F}(X)\times \mathcal{F}(X)$. Suppose $\mathcal{H}$ is a subalgebra in $\mathcal{G}$ generated by the elements of the form $(r,0)$, where $r \in R$,  and also by the set of ``diagonal'' elements $\{ (x,x)\,|\, x\in X\}$. Clearly, $\mathcal{H}$ contains $(I,0)$. Conversely, if $(f,0) \in \mathcal{H}$, then $f$ must be an element of $(I,0)$ which follows because $(I,0)$ is the kernel of the projection of $H$ onto the right factor $\mathcal{F}(X)$ in $\mathcal{G}$.  As a result,  $(f,0)\in  \mathcal{H}$ if and only if $f\in I$. This also tells us that we cannot effectively write the basis of the intersection $\Sp{(f,0)}\cap \mathcal{H}$. Hence the (Generalized) Membership Problem for $\mathcal{H}$ in $\mathcal{G}$ is algorithmically undecidable. According to the above  Theorem \ref{BOpLDP}, the embedding of $\mathcal{H}$ in $\mathcal{F}(X)\times \mathcal{F}(X)$ is distorted, moreover the distortion function is not bounded by any recursive function. }
\end{example}

Replacing $\mathcal{F}(X)\times \mathcal{F}(X)$ by $\mathcal{F}(X)$ is not possible in many important (quasi)va\-rieties. In other words, free algebras of finite rank in many (quasi)varieties have no \fg distorted subalgebras. This is either known (as in the case of groups), or will be shown later in this paper (as in the case of Lie algebras).

\subsection{Multitude of tame filtrations on algebras}\label{sDA}

As in Introduction, any degree filtration of a finitely generated (linear) $\Omega$-algebra is a tame filtration. Indeed, suppose $\# X=r$, $\#\Omega=p$. If an element $a\in A$ is in the $n\th$ term $A_n$ of the degree filtration defined by $X$ on an $\Omega$-algebra $A$ then by induction it easily follows that $a$ can be represented by a word in an $(r+p)$-letter alphabet $X\cup\Omega$ whose length is at most $nq+1$, where $q$ is the maximal arity of operations in $\Omega$. Thus $\# A_n\le C((r+p)^q)^n$, for a certain constant $C$. Therefore, the growth of the sequence of the numbers $\# A_n$ is at most exponential, proving that any degree filtration is a tame filtration. As it turns out, the degree filtrations are just a tip of an iceberg in the vast array of all tame filtrations.
 
\begin{theorem}\label{BOtCDF}
Let $\Omega$ be a finite signature. On any countable $\Omega$-algebra or countably - dimensional linear $\Omega$-algebra there is at least continuum of pairwise non-equivalent tame filtrations. Given any function $g(n)$, one can always choose two tame filtrations $\alpha=\{ A_n\}$ and $\alpha'=\{ A_n^{\prime} \}$ so that for the distortion function $f=\dist_\alpha^{\alpha'}$ one has $f(n)>g(n)$, for all $n>0$.
\end{theorem}

\pp In our proof, we restrict ourselves to the case of $\Omega$-algebras, the case of linear $\Omega$-algebras being quite analogous. Thus, let $A$ be an infinite \fg algebra of finite signature $\Omega$, with finite generating set $X$, $\# X=r$. Let $\{ A_n\}$ be the degree filtration of $A$ determined by $X$. For any $a\in A$ we denote by $d(a)$ its degree with respect to $\{ A_n\}$, that is, $d(a)=n$, if $a\in A_n\setminus A_{n-1}$. As noted previously, $d(a)$ is the minimal number of operations necessary to obtain $a$, starting from $X\cup\Omega_0$. To simplify the notation, we will assume in the future that $\Omega_0\subset X$.

Now let us fix a real number $\lambda$ such that $0<\lambda<1$.  We choose an infinite set $\cw^\lambda$ of elements $w\in A$ with pairwise different values of $[d(w)^{\lambda}]$, for any $w\in\cw^\lambda$. Here $[r]$ is the integral part of a real number $r$. Such set has to exist because $A$ is infinite and each $A_n$ is finite. 

We are going to define a new filtration, $\{ A_n^{\lambda}\}$, in the following way. We choose a set $Y\cup \mv$ so that there is bijection $\vp: Y\cup \mv\to X\cup \mw$ with $\vp(Y)=X$ and $\vp(\mv)=\mw$. We define a function $\da$ on absolutely free algebra $\mathcal{F}$ of signature $\Omega$ with free generating set $Y\cup \mv$ (notice that if $\Omega$ is the group signature then, being nonassociative, an absolutely free $\Omega$-algebra is a much ``larger' object than the free group!). The elements of $\mathcal{F}$, we call them \emph{words}, appear by induction on the number of applications of operations in $\Omega$ starting from the elements of $Y\cup \mv$. We set $d^{\lambda}(y)=0$, for all $y\in Y$, and if $v\in\mv$, then $\vp(v)=w$, for a unique $w\in\mw$, and we set  $d^{\lambda}(v)=[d(w)^{\lambda}+1]$. Using induction, let us assume that we have already assigned the values of $d^{\lambda}$ to some words $u_1,\ld,u_m\in \mathcal{F}$. If $\omega\in\Omega_m$, $1\le m\le q$, and $u=u_{\,1}\ld u_m\omega$ then we set $d^{\lambda}(u)=d^{\lambda}(u_{\,1})+\cdots+d^{\lambda}(u_m)+1$.

Now, for any $n>0$, we define the $n\th$ term of a filtration of $\mathcal{F}$ by setting $\mathcal{F}_n^\lambda=\{u\,|\,d^{\lambda}(u)\le n\}$.  Clearly, $\{ \mathcal{F}_n^\lambda\}$ is an ascending filtration in $\mathcal{F}$. Suppose we have shown that this is also a tame filtration.  Consider a unique homomorphism $\overline{\vp}: \mathcal{F}\to A$ extending the above bijection $\vp: Y\cup\mv\to X\cup\mw$ and set $A_n^{\lambda}=\overline{\vp}(\mathcal{F}_n^{\lambda})$. Then $A_n^{\lambda}$ is a tame filtration in $A$. To provide ourselves with continuum of pairwise inequivalent tame filtrations on $A$, we will later need to prove that  $\{ A_n^{\lambda}\}\not\sim \{ A_n^{\mu}\}$ if $\lambda\ne\mu$.

Let us denote by $S_n^{\lambda}$ the \emph{sphere} of radius $n$ in $\mathcal{F}$, that is, the set of all words $u$ with $d^{\lambda}(u)= n$. Let $\mv_n$ be the set of letters $v$ in $\mv$ with $\da(v)=n$. Then $S_0^{\lambda}=X$ and, using induction on $n\ge 1$, $S_n^{\lambda}=\bS_n^{\lambda}\cup\mv_n$ where $\bS_n^{\lambda}=\bigcup(\bS_n^{\lambda})_{\omega}$ with $(\bS_n^{\lambda})_{\omega}=\bigcup_{i_1+\cd+i_m+1=n}S_{i_1}^{\lambda}\ld S_{i_m}^{\lambda}\omega$, for each $\omega\in\Omega_m$ $(m\ge 1)$. By our hypothesis about $\mv$, $\# S_n^{\lambda}\le \#\bS_n^{\lambda}+1$.

It suffices to prove by induction on $n$ that $\# S_n^{\lambda}\le c^n$, for big enough constant $c$. However, for our argument using induction on $n$, we need to prove a stronger inequality $\# S_n^{\lambda}\le \dfrac{c^n}{(n+1)^2}$, where the constant $c$ is big enough. The choice of $c$ is dependent on the cardinalities of $\Omega$ and $X$. Since the sphere of radius 1 is finite, let $c$ be such that the above inequality holds for $n=1$.

To handle the case of $n>1$, we will need the following. 

\begin{lemma}\label{lBOsec3} For any $m\ge 1$ there is a constant $C(m)>0$, such that 

\[
\sum_{i_1+\cdots+i_m=n} \prod_{k=1}^m  \frac{1}{(i_k+1)^2} \le \frac{C(m)}{(n+2)^2},
\]

\noindent for any $n$.
\end{lemma}

\pp Use induction by $s$ and the convergence of the series $1+\dfrac{1}{2^2}+\dfrac{1}{3^2}+\cdots$.\epf

Assuming our upper bound for the number of words $u$ with $d^{\lambda}(u)\le k$ true for any $k<n$, we will now estimate the number of words of the form $u_1\ld u_m\omega$ where $u_1\in \mathcal{F}_{i_1},\ld,u_m\in \mathcal{F}_{i_m}$, $\omega$ an operation in $\Omega_m$ and $i_1+\cdots +i_m+1=n$.
Using Lemma \ref{lBOsec3} and the induction hypothesis, this number does not exceed

\[
\sum_{i_1+\cdots+i_m+1 = n} \prod_{k=1}^m \frac{c^{i_k}}{(i_k+1)^2} \le c^{i_1+...+i_m} \dfrac{C(m)}{(n+1)^2}
= c^{n-1} \dfrac{C(m)}{(n+1)^2}.
\]

To obtain the upper bound for $\# S_n^{\lambda}$ we have to sum all the expressions just obtained over $\omega\in\Omega\setminus\Omega_0$ and add $1$ in order to take care of possible contributions of letters from $\mv$. This brings us to $1+\left(\dfrac{K}{c}\right) \dfrac{c^{n}}{(n+1)^2}$. Here $K$ is the sum of all $C(m)$ (with possible repetitions) over all $\omega\in\Omega$, of arity $m\ge 1$. It is quite obvious now that $c$ can be chosen so that the following is satisfied for all $n$: 

\[ 1+\left(\dfrac{K}{c}\right) \dfrac{c^n}{(n+1)^2}\le \dfrac{c^n}{(n+1)^2},\mbox{ for all }n>1.
\]

Thus we have shown that $\# S_n^{\lambda}\le \dfrac{c^n}{(n+1)^2}$, for all $m\ge 1$. For this $c$ we will also have $\# S_n^{\lambda}\le c^n$, hence the growth of $\{\# \mathcal{F}_n^{\lambda}\}$ is bounded by an exponential function, as claimed.

Notice that, for any word $u$, one has $d^{\lambda}(u)\ge d(u)^{\lambda}$. This is true for $u\in Y\cup\mv$. Now if $u=u_1\ld u_m\omega$, then, considering $\lambda< 1$ and using induction, we obtain $d^{\lambda}(u) = \sum d^{\lambda}(u_i)+1\ge
\sum d(u_i)^{\lambda} +1 \ge  (\sum d(u_i))^{\lambda}+1 =(d(u)-1)^{\lambda}+1
\ge d(u)^{\lambda}$. In terms of filtrations, if $u\in \mathcal{F}$ is such that $d(u)^{\lambda}>n$ then $u\not\in \mathcal{F}_n$. Since $A_n=\overline{\vp}(\mathcal{F}_n)$, we have that no elements $a\in A$ with $d(a)^{\lambda}>n$ can be contained in $A_n$.

Now we can prove that the filtrations obtained for different values of parameter $\lambda$ are not equivalent. Indeed, suppose $0<\lambda<\mu<1$. Assume that there is natural number $t$ such that the $n\th$ term of the filtration defined by $\lambda$ is a subset of the $(tn)\th$ term of the filtration defined by $\mu$, for all $n\ge 1$. Let us find $a\in \mw$ such that $[d(a)^{\lambda}+1]=n$, for some $n$ such that  $(n-1)^{\mu/\lambda}>tn$, which is possible because $\lambda<\mu$. In particular, $a\in A^{\lambda}_n$. Since $d(a)\ge (n-1)^{1/\lambda}$, it follows that $d(a)^{\mu} \ge ((n-1)^{1/\lambda})^{\mu} > tn$. Hence, $a\notin A_{tn}^{\mu}$, a contradiction. This completes the proof of the first claim of our theorem in the case where $A$ is finitely generated.

To complete the proof of the first claim in the general case, let us assume that $A$ is not finitely generated but has a countable set of generators $\{ a_0, a_1,...\}$ such that for each $i\ge 1$, we have $a_{i}\not\in \alg\{ a_1,\ld, a_{i-1}\}$. Again we introduce free $\Omega$-algebra $\cf$, now with free generating set $\{ y_0, y_1,\ld\}$ and define an epimorphism $\vp:\cf\to A$ extending the bijection $\vp(y_i)=a_i$, for $i=0,1,\ld$.

Now let us choose a real number $\lambda\ge 1$ and define a degree function on $\cf$ by setting $d^{\lambda}(y_i)= [\,i^{\lambda}]$, for all $i\ge 0$. As earlier, $d^{\lambda}$ extends to the whole of $\cf$ and induces a function on $A$ denoted by the same symbol. Using this function, one defines a filtration $\{ A_n^{\lambda}\}$ on $A$ such that $A_0$ consists of the values of $0$-ary operations and element $a_0$. Since $a_i\not\in\alg\{ a_0,\ld,a_{i-1}\}$, it follows that $d^{\lambda}(a_i) = d^{\lambda}(y_i)= [\,i^{\lambda}]$, for  $i=0, 1, \ld$

Now the proof of Lemma \ref{lBOsec3} goes without any changes, including the fact $\# S_n^{\lambda}\le \#\bS_n^{\lambda}+1$, because, for $n$ fixed, there is at most one $i$ such that $d^{\lambda}(y_i)=n$. As before, it follows that $\{ A_n^{\lambda}\}$ is a tame filtration. Finally, if $\lambda<\mu$ then $\{ A_n^{\lambda}\}\not\sim \{ A_n^{\mu}\}$. This follows because 
$d^{\lambda}(a_i)/d^{\mu}(a_i)\to 0$ when $i\to \infty$.

Now let us fix an increasing function $g:\N\to\N$. Again, we start with the case where $A$ is finitely generated. We choose an infinite set $\cw^{\,\prime}$ of elements $w_1,\ld,w_n,\ld\in A$ so that $d(w_n)> g(n)$. Then we consider an absolutely free $\Omega$-algebra $\mathcal{F}$ with free generating set $Y\cup\mathcal{V}^{\,\prime}$ which is in bijection $\vp$ with $X\cup\cw^{\,\prime}$, as before. Let $\overline{\vp}$ be an epimorphism extending the above bijection. In the same manner, as previously, we define the filtration $\mathcal{F}_n^{\,\prime}$, except that  now $d^{\lambda}$ should be replaced by $d^{\,\prime}$, such that $d^{\,\prime}(y)=0$, for all $y\in Y$, and $d^{\,\prime}(v_n)=n$, for all $n=1,2,\ld$, where $v_n\in\mathcal{V}^{\,\prime}$ is such that $w_n=\vp(v_n)$. The proof that $\{ \mathcal{F}_n^{\,\prime}\}$ is a tame filtration is the same, as before. We now set $A_n^{\,\prime}=\vp(\mathcal{F}_n^{\,\prime})$ and consider the distortion function $f=\dist_{\alpha}^{\alpha^{\,\prime}}$ of a tame filtration $\alpha^{\,\prime}=\{ A_n^{\,\prime}\}$ with respect to the degree filtration $\alpha=\{ A_n\}$. By definition, $A_n^{\,\prime}\subset A_{f(n)}$. Since $w_n\in A_n^{\,\prime}\setminus A_{g(n)}$, it follows that $f(n)>g(n)$, for all $n$, as claimed. This completes the proof of the second claim of our theorem in the case where $A$ is finitely generated.

To complete the proof of the second claim in the general case, we may now  
assume that $A$ is not finitely generated but has a countable set of generators $a_0, a_1,...$ such that for each $i\ge 1$, we have $a_{i}\not\in \alg\{ a_1,\ld, a_{i-1}\}$. Again we introduce free $\Omega$-algebra $\cf$ with free generating set $y_0, y_1,\ld$ and define an epimorphism $\vp:\cf\to A$ extending the bijection $\vp(y_i)=a_i$, for $i=0,1,\ld$. Given a function $g$ as above, we  choose $d(y_i)$ arbitrarily with $d(y_i)> \max\{  g(k)\,|\, k\le i\}$. We also set $d(a_i)=d(y_i)$. Now let us consider two filtrations: $\alpha=\{ A^{\lambda}_n\}$ with $\lambda=1$ and $\alpha'=\{ A_n'\}$ defined by $d$, that is, $A_n'=\{ a\,|\, d(a)\le n\}$. In this case,  $a_i \in A^{\lambda}_i\setminus A_{g(i)}'$. Thus, $\ds{\alpha}{\alpha'}>g(n)$. \epf

\begin{corollary}\label{cBOSGLAA} Let $A$ be a countable (semi)group or countably - dimensional associative or Lie algebra over a field. Then $A$ has at least continuum of pairwise inequivalent tame filtrations. Given any function $g: \N\to\N$ one can always find two tame filtrations $\alpha$ and $\alpha'$ in $A$ such that $\dist_{\alpha}^{\alpha'}(n)>g(n)$, for all $n\ge 1$.\epf
\end{corollary}

\begin{remark}\label{hopelast} \emph{We can be more precise about the cardinality of the set of pairwise inequivalent tame filtrations.} 

\emph{First, in the case of countable algebras this cardinality is continuum simply because the cardinality of the countable Cartesian power of finite subsets of a countable is continuum.}

\emph{Second, in the case of countably-dimensional linear $\Omega$-algebras the cardinality of this set can be even greater than continuum, namely as large as the cardinality of the ground field $F$. In Subsection \ref{ssEDF} we will produce such examples already on the polynomial algebra $F[x]$ in one variable.}
\end{remark}

\section{Distortion in classical algebras}

In the previous section we have seen that tame filtrations are abound on virtually any infinite(-dimensional) algebra. As announced in Introduction, we will show in Section 3 that in the classical situation of associative and Lie algebras, every tame filtration of a countably-dimensional algebra $B$ is the restriction of a degree filtration of a suitable finitely generated algebra $A$. But there are important classes of algebras such that if $A$ is a \fg algebra of this class and if a tame filtration of a \fg subalgebra $B$ is the restriction of a degree filtration of $A$, then it is equivalent to a degree filtration of $B$. In other words, all subalgebras of such algebras $A$ are undistorted. As examples, we show that, among others, this property holds for commutative associative algebras and free Lie algebras. Oddly enough, free noncommutative associative algebras may and do have distorted subalgebras. One more topic we study in this section is the connection between the distortion of the embedding of Lie algebras $M\subset L$ and their associative counterparts $U(M)\subset U(L)$. This is interesting when compared with a similar situation with the embedding of groups and their group algebras, see Proposition \ref{BOpGA}.

\subsection{Commutative associative algebras}\label{sNDCA}

For the proof of our main result in this subsection we use simple facts of the theory of Groebner bases in the commutative algebras. One of the many possible references is \cite{KR}.

\begin{theorem}\label{Op1}
The embedding of a finitely generated subalgebra in a \fg commutative associative algebra has no distortion.
\end{theorem}

\pp Let $A$ be a \fg commutative algebra with generators $a_1,\ld,a_k$, $B$ a subalgebra generated by $b_1,\ld,b_{\ell}$, $\alpha=\{ A_n\}$ and $\beta=\{ B_n\}$ respective degree filtrations in $A$ and $B$. We write $A$ as a factor algebra of the polynomial algebra in $k+\ell$ variables $F[x_1,\ld,x_k,y_1,\ld,y_{\ell}] / I$,  where $a_i=x_i+I$, $b_j=y_j+I$, for $i=1,\ld,k$ and $j=1,\ld,\ell$.

We introduce an order on the monomials in $x_1,\ld,x_k,y_1,\ld,y_{\ell}$ in the following way. Let us assume that any $x_i$ is greater than any monomial in $y_1,\ld,y_{\ell}$, but $x_iy_j > x_r$, for any indexes $i,j,r$. Two monomials in $x_1,\ld,x_k$ are ordered by Short-lex, and the same is true for two monomials in $y_1,\ld,y_{\ell}$. In other words, we start comparing two monomials in $x_1,\ld,x_k,y_1,\ld,y_{\ell}$ first by $\deg_X$, where $X=\{ x_1,\ld,x_k\}$. If these degrees are the same then compare by $\deg_Y$, where $Y=\{ y_1,\ld,y_{\ell}\}$. If even these degrees are the same, we begin comparing lexicographically, considering $x_i >y_j$ for any $i,j$. We keep in mind that in the theory of Groebner bases it is allowed to introduce the order in many different ways provided that the set of words becomes a totally ordered monoid with minimum $1$.

Let us consider the Groebner basis of monic polynomials $g_1,\ld,g_s$ in $I$ with respect to this order. We also choose $C>0$ so that $\deg_Y(g_i) \le C$, for all $i=1,\ld,s$.

Being Groebner basis in an ideal $I$ means that $g_1,\ld,g_s\in I$ and that the leading monomials $lm(g_1),\ld,lm(g_s)$ of $g_1,\ld,g_s$ generate an ideal $\ell(I)$ spanned by the leading monomials of all elements in $I$. A linear basis of $F[x_1,\ld,x_k,y_1,$ $\ld,y_{\ell}]$ modulo $I$ is formed by all monomials which are not divisible by any of $lm(g_1), \ld, lm(g_s)$.

Any polynomial $g = g(x_1,\ld,y_\ell)$ can be reduced modulo $I$ to a unique linear combination of such monomials using transformations of the following kind. If $g$ contains a monomial $u = u'u''$ , where $u' = lm(g_i)$, for some $i$, then we replace $u$ by $(u' - g_i)u''$, followed by reduction of like terms.

Notice that under these transformations the weighted ``degree'' $C \deg_X (g) + \deg_Y (g)$  of $g$ cannot grow. Indeed, the ``degree'' $C\deg_X + \deg_Y$ of any monomial in the expression of $g_i$ is at most the ``degree'' of $lm(g_i)$, which follows from the definition of the order and the inequality $\deg_Y (g_i) \le C$. In particular, if we start with a polynomial $g(x_1,\ld,x_k)$ of degree $d$ and end up with a polynomial $f(y_1\ld,y_{\ell})$ of degree $m$ then $m\le Cd$.

Now we can show that $\ds{\alpha\cap B}{\beta}(n)\le Cn$, for any $n=1,2,\ld$. Indeed, let $b\in A_n\cap B$. In this case $b=g(a_1,\ld,a_k)$ where $g(x_1,\ld,x_k)$ is a polynomial of degree $\deg_X g=d\le n$. Let us pass from $g(x_1,\ld,x_k)$ to a polynomial $h(x_1,\ld,y_{\ell})$, which is reduced modulo $I$. By the above, $C\deg_X(h)+\deg_Y(h)\le C\deg_X g\le Cn$. Now since $b\in B$, $b=f(b_1,\ld,b_{\ell})$, for a polynomial $f(y_1,\ld,y_{\ell})$ of some degree $m=\deg_Y f$. Let us assume that $f(y_1,\ld,y_{\ell})$ is of minimal possible degree with respect to $Y$, representing $b$ and depending only on $Y$. We can pass from $f(y_1,\ld,y_{\ell})$ to $h(x_1,\ld,y_{\ell})$ by transformations of the form indicated above. It follows from the definition of the order on the monomials in $X\cup Y$, that the leading monomial $lm(g_i)$ for an element $g_i$ of the Groebner basis for $I$ is a monomial in $Y$ only if $g_i$ is a polynomial in $Y$ and the degree of $lm(g_i)$ is the greatest possible among all the other monomials in $g_i$. By the minimality of the degree of $f(y_1,\ld,y_{\ell})$, it then follows that $h$ is a polynomial in $Y$ and its degree is the same as the degree of $f$. Then by the above, $m\le Cn$, proving our claim.\epf

\begin{example}\label{BOeUUU} \emph{One of possible ways to generalize this result is to consider algebras satisfying nontrivial polynomial identity (PI-algebras) (for some basic facts about algebras with polynomial identities, see \cite{BOSRV}). Unfortunately, a result by Umirbaev \cite{UUU} shows that a free algebra $\cf_{\mathfrak{V}}(Y)$, freely generated by a finite set $Y$, in the variety $\mathfrak{V}$ of associative algebras over computable fields, defined by the following identity}

\begin{equation}\label{BOeIOU} 
[x_1,x_2][x_3,x_4][x_5,x_6][x_7,x_8]\equiv 0,
\end{equation}
 
\emph{contains a \fg subalgebra $B$ with undecidable Membership Pro\-blem. Also, the Intersection Problem is decidable in $\cf_{\mathfrak{V}}$. Indeed, if $\ca(Y)$ is the \fa with free generating set $Y$ (see Subsection \ref{sFSA} for formal definition) then one may write a vector space decomposition $\ca(Y)=\cf_{\mathfrak{V}}(Y)\oplus V(Y)$, where $V(Y)$ is the set of all consequences of (\ref{BOeIOU}) in the $\ca(Y)$. Since (\ref{BOeIOU}) is multilinear, one can write the same equation for each term of the degree filtration of $\ca(Y)$: $\ca(V)_n=\cf_{\mathfrak{V}}(Y)_n\oplus V(Y)_n$. This latter is an equality for finite-dimensional spaces. The space $V(Y)$ is spanned by finitely many monomials $u[v_1,v_2][v_3,v_4][v_5,v_6][v_7,v_8]w$, where the sum of lengths of all words $u,v_1,\ld,w$ is at most $n$. Since $\ca(Y)$ is the \fac we can effectively write the basis of $V(Y)_n$ and hence of $\cf_{\mathfrak{V}}(Y)_n$. Now, given two finite subsets $S$ and $R$ in $\cf_{\mathfrak{V}}(Y)$, we effectively find $n$ such that $S\cup R\subset \cf_{\mathfrak{V}}(Y)_n$, find the basis of $\cf_{\mathfrak{V}}(Y)_n$, as described above, and solve the Intersection Problem for $\Sp{S}\cap\Sp{R}$. Thus Theorem \ref{BOpLDP} applies, showing that none of the distortion functions for $B$ in $\cf_{\mathfrak{V}}(Y)$ is bounded by a recursive function, certainly not by any linear function. Therefore, the embedding of $B$ in $\cf_{\mathfrak{V}}(Y)$ is necessarily distorted.}
\end{example}

It is easy to see that the above argument providing us with an algorithm for the solution of the Intersection Problem works in the case of every variety defined by not only by multilinear but also homogeneous identities within a variety of all algebras, or all associative algebras, or all Lie algebras, etc. In the case of homogeneous rather than multilinear identities one has to complement the set of identities by their partial  linearizations.


\subsection{Free associative algebras and free group algebras}\label{sFSA}

Let $X$ be a nonempty alphabet and $\cw(X)$ the set of all \emph{words} (including the empty word 1) in $X$. One calls $\cw(X)$ the \emph{free monoid} with basis (free generating set) $X$. Let us fix the field $F$ of coefficients and consider the \emph{\fa} $\ca(X)$ over $F$. The elements of $\ca(X)$ are linear combinations of words in $\cw(X)$ with coefficients in $F$. 

The elements of $\ca(X)$ are often called (noncommutative) \emph{polynomials} in $X$. The degree $n=\deg f$ of $f\in \ca(X)$ is the length of the longest word that enters the decomposition of $f$ through $\wax$ with nonzero coefficient. Thus the degree of a monomial $\lambda w$, $0\ne\lambda\in F$, $w\in\wax$ is just the length of $w$.

If we replace $\cw(X)$ by the free group $F(X)$ then the resulting algebra $\cf(X)$ is called a \emph{free group algebra}.

We start by giving a quick proof to an example by U. Umirbaev \cite{UUU1}. Let $\ca=\ca(x,y,z)$ be a \fa with free generators $x,y,z$ and $C=\ca(x,y)$ a (free) subalgebra generated by $x,y$.  Let $I$ be the ideal of $C$ generated by some elements $f_1,\ld,f_k$. We denote by $B$ the subalgebra of $\ca$ generated by the elements $zf_1, \ld,zf_k$, $x$, $y$, $[z,x]=zx-xz$ and $[z,y]=zy-yz$.

\begin{lemma}\label{UOl1}
Let $f$ be an arbitrary element of $C$. Then $f$ is an element of $I$ if and only if $zf$ is an element of $B$. 
\end{lemma}

\pp First we assume $f\in I$. Then $f$ is the sum of expressions of the form $uf_iv$, where $u, v$ are monomials in $x,y$. Let us prove $zf\in B$. It is sufficient to show $zuf_i\in B$, because $v \in B$. We proceed by induction on $\deg u$ with trivial basis when $\deg u=0$. If $u=xu'$ then  $zuf_i = zxu'f_i = [z,x] (u'f_i) + x(zu'f_i)$.  Here $[z,x]\in B$, $u'f_i \in B$ (because this depends only on $x,y$), while the second summand is in $B$ by the induction hypothesis. If $u=yu'$ then we argue similarly but use $[z,y]\in B$.  

Now we assume that $zf \in B$. Let us prove $f \in I$.  For this we notice that $zf$ is homogeneous of degree $1$ in $z$, hence $zf$ is a linear combination of expressions of one of the forms $u[z,x]v$, $u'[z,y]v'$ and  $u''(zf_i) v''$, where $u, v, u', v', u'', v''\in C$. Here the left factors $u, u', u''$ can be removed because in the product $zf$ all monomials start with $z$.  In the expressions $z(f_iv'')$ each factor $f_iv''$ is in $I$ whereas $z(f_iv'')$ itself is in $B$. Therefore, without loss of generality, we can remove these expressions if we replace $f$ by $f - f_iw$. After doing so, we arrive at the equation $[z,x] v + [z,y] v' = zf$, where $v,v'$ depend on the generators $x,y$ only.  A part of this sum, $-xz v - yz v'$ equals zero because all the monomials should start with $z$.  It follows that $v = v' = 0$, which follows because the right $C$-module generated by $xz$ and $yz$ is free. As a result, $[z,x]v+[z,y]v' = zf = 0$, hence $f=0$. Considering that this ``final'' $f$ was obtained by a number of transformations of the form $f\to f-f_iw$, we conclude that, for the ``original'' $f$ we must have $f\in I$. \epf

It is well-known that over any field there exist 2-generator associative algebras with undecidable equality problem \cite{KS}. As a result, by Lemma \ref{UOl1}, we have an example of a \fg subalgebra $B$ in a 3-generator free algebra $\ca$ with undecidable membership problem. It is very easy to obtain an example of a subalgebra in a 2-generator \fac with undecidable membership problem. For this, we simply need to effectively embed $\ca$ in a \fas $\ca'=\ca(u,v)$ with free generators $u,v$, for instance, by setting $x=u^2$, $y=uv$ and $z=v^2$. Then obviously, there is no algorithm to decide the membership problem for the image of $B$ in $\ca'$. 

It is now very easy to notice that also in any free group algebra $\cf(X)$ with $\# X\ge 2$ there is a subalgebra with undecidable membership problem. Indeed, $\ca(X)\subset \cf(X)$, and there is an obvious algorithm to decide whether an element of $\cf(X)$ is an element of $\ca(X)$. Therefore, if $B$ is a subalgebra with undecidable membership problem in $\ca(X)$, it will also be is a subalgebra with undecidable membership problem in $\cf(X)$.

Using the above in conjunction with Theorem \ref{BOpLDP} allows us to make the following conclusion.

\begin{theorem}\label{BOpDSFAS}
\emph{\textbf{(Examples)}} Let  $A$ be a \fa $\ca(X)$ or a free group algebra $\cf(X)$ with $\# X\ge 2$ over a field $F$. Then one can choose a finitely generated distorted subalgebra $B$ in $\ca$ in such a way that none of the distortion functions for $B$ in $\ca$ is bounded from above by any recursive function.
\end{theorem}

\pp The argument in the case of free group algebras being identical to that in the case of free associative algebras, we restrict ourselves to these latter ones.

 Since the above example of a finitely generated subalgebra with undecidable membership problem in a \fa or $\cf(X)$ of rank at least 2 does not depend on the field, we first consider a free associative algebra $\ca_0$ of rank at least 2 over the prime subfield $F_0$ of $F$. Then $F_0$ is constructive and by Theorem \ref{BOpLDP} we have a distorted subalgebra $B_0$ in $\ca_0$. Now let us consider the $F$-algebra $\ca=F\o_{F_0} \ca_0$ and its $F$-subalgebra $B=F\o_{F_0} B_0$. Suppose that $\beta_0=\{ (B_0)_n\}$ and $\alpha_0=\{ (\ca_0)_n\}$ are degree filtrations in $B_0$ and $\ca_0$, respectively, defined by the generating sets $Y$ and $X$, respectively, such that the distortion function $\dist_{\alpha_0\cap B}^{\beta_0}$ is not bounded by a linear function. Then setting $\beta=\{ F\o_{F_0} (B_0)_n\}$ and $\alpha=\{ F\o_{F_0} (\ca_0)_n\}$ gives rise to the degree filtrations on $B$ and $\ca$, respectively, defined by the generating sets $1\o Y$ and $1\o X$, with $\dist_{\alpha\cap B}^{\beta}=\dist_{\alpha_0\cap B}^{\beta_0}$. Therefore, $B$ is a distorted subalgebra of the free associative algebra $\ca$ over an arbitrary field $F$.\epf

\begin{example}\label{BOrDSFJA}
\emph{In Umirbaev's paper \cite{UUU2} there are examples of subalgebras in \emph{free Jordan algebras} similar to those associative in Theorem \ref{BOpLDP}. As we noted at the end of Subsection \ref{sNDCA}, there is no problem with the Intersection Problem in the case of free Jordan algebras of finite rank, since the variety of Jordan algebras is defined within the variety of all linear algebras by homogeneous identical relations. Thus we conclude that \emph{free Jordan algebras of finite rank contain \fg distorted subalgebras}.}
\end{example} 

 \begin{example}\label{BOrSGA}
 \emph{It is well-known (and follows from Schreier rewriting) that any \fg subgroup of a free group of finite rank is undistorted. The same is true in the case of any semigroup $A$ with finite generating set $X$ and \emph{balanced} defining relations, that is relations of the form $u = v$ where $u,v\in\cw(X)$ have equal lengths: $\ell(u)=\ell(v)$. In such a semigroup the length $\ell(g)$ of any element $g$ is well-defined and for any two elements $g,h\in A$ one has $\ell(gh)=\ell(g)+\ell(h)$. If $B$ is a subsemigroup of $A$ generated by a subset $Y$ and $g_1g_2\cdots g_n$ is an element of $B$ then its length with respect to $X$ will be at least $n$, proving that $B$ is undistorted in $A$. Thus it follows from Theorem \ref{BOpDSFAS} that \emph{there exist (semi)group algebras of (semi)groups without distorted sub(semi)groups, which have subalgebras whose distortion function is not even bounded by any recursive function}. As shown in Proposition \ref{Op1}, this cannot happen for commutative (semi)groups.}
\end{example}


\subsection{Free Lie algebras}\label{sFLA}

Let $\mathcal{L}=\mathcal{L}(X)$ be a free Lie algebra over a field $F$ with free generating set $X=\{x_1,\ld,x_m\}$. One can view $\mathcal{L}$ as a subalgebra in \fa $\ca(X)$ generated by $X$, with respect to the bracket operation $[a,b]=ab-ba$. Thus the degree $\deg_X f$ of an element $f\in \mathcal{L}(X)$ is its degree as an element of $\ca(X)$. Each element $f$ has its \emph{leading part} $\mathrm{Lp}\,(f)$, which is the homogeneous component of maximal degree.

\begin{theorem}\label{Op2}
The embedding of a finitely generated subalgebra in a \fg free Lie algebra has no distortion.
\end{theorem}

\pp By Shirshov - Witt's Theorem \cite[Chapter 3]{B} any subalgebra $B$ of a free Lie algebra $\cl$ is itself free. Moreover, one can choose a free generating set $Y$ in $B$ so that the set $\overline{Y}=\{\mathrm{Lp}\,(y)\,|\,y\in Y\}$ is the free generating set of the subalgebra it generates in $\cl$.

Now let us assume that $B$ is a finitely generated subalgebra of $\cl(X) $ and that $Y=\{ f_1,\ld,f_n\}$ is the set of free generators of $B$. Let the degrees of $\mathrm{Lp}\,(f_1)=f'_1,\ld,\mathrm{Lp}\,(f_n)=f'_n$ be $d_1,\ld,d_n$, respectively.
 
Let us consider another free Lie algebra $\mathcal{L}(Z)\subset \ca(Z)$, where $Z=\{ z_1,..,z_n\}$. In addition to the ordinary degree filtration defined by $Z$, with degree function $\deg_Z$, we can endow  $\mathcal{L}(Z)$ with another filtration, with degree function $\dg$, defined on the words in $Z$ by induction on the length as follows. We set $\dg z_1=d_1,\ld,\dg z_k=d_k$, for the letters of the alphabet $Z$. If $w=uv$ is a word of degree $>1$ then by induction $\dg u$ and $\dg v$ are already defined and we set $\dg uv=\dg u+\dg v$. Given an element $h(z_1,\ld,z_n)\in\ca(Z)$ which is a linear  combination of words $u_1,\ld,u_s$, with nonzero coefficients, we set $\dg h(z_1,..,z_n)=\max\{ \dg u_1,\ld,\dg u_s\}$. Finally, we denote by $\mathrm{Lp}'(h)$ the ``leading part'' of $h(z_1,\ld,z_n)$ with respect to the new degree $\dg$, that is, the homogeneous component of the above decomposition of degree $\dg h(z_1,..,z_n)$.
 
  Now, for any nonzero Lie polynomial $h(z_1,...,z_n)$ we consider the equality $g(x_1,\ld,x_m) = h(f_1,\ld,f_n)$. Then the leading part of $g(x_1,\ld,x_m)$ will be $\mathrm{Lp}'(h)(f'_1,$ $\ld,f'_n))$. This polynomial in $X$ is nonzero since $\{ f'_1,\ld,f'_n\}$ is a free generating set.
 
  It follows that $\deg_X g= \dg h \ge \deg_Z h$. Let $\alpha$ be the degree filtration of $\cl(X)$ defined by $X$ and $\beta$ the degree filtation of $B$ defined by $Y=\{ f_1,\ld,f_n\}$. We just showed that $\ds{\alpha\cap B}{\beta}(n)\le n$. Thus $B$ is an undistorted subalgebra in $\cl(X)$.\epf
  
\begin{remark}\label{BOrFLAUA} \emph{Since the universal enveloping algebra $U(\cl(X))$ is isomorphic to the free associative algebra $\ca(X)$ (see, e.g. \cite[Chapter 3]{B}) and by Theorem \ref{BOpDSFAS} free associative algebras of rank $\ge 2$ have distorted subalgebras, it follows that universal enveloping algebras of Lie algebras without distorted subalgebras can have distorted subalgebras. This is quite similar to the case of (semi)groups and their (semi)group algebras (see Remark \ref{BOrSGA} above). However, in both cases just mentioned, the situation is much better when we consider the connection of the distortion of $H$ in $G$ for a \fg subgroup $H$ of a \fg group $G$ or a \fg subalgebra $H$ in \fg Lie algebra $G$ with the distortion of a respective subalgebra $F[H]$ in the group algebra $F[G]$ or $U(H)$ in the universal enveloping algebra $U(G)$. In the case of groups this is described in our simple Proposition \ref{BOpGA} above, where it was shown that $\ds{G}{H}=\ds{F[G]}{F[H]}$. The case of Lie algebras requires greater effort, and this will be done in the next subsection of the paper.}
\end{remark} 
\subsection{Universal enveloping algebras}\label{sDUA}

If $L$ is a Lie algebra with a finitary filtration $\alpha=\{ L_n\}$ then its universal enveloping algebra $U(L)$ is naturally endowed with a finitary filtration $U(\alpha)=\{ U(L)_n\}$ where 

\[
U(L)_n=\sum_{n_1+\cd+n_m\le n}L_{n_1}\cd L_{n_m}.
\] 

If $L$ is finitely generated and $M$ a \fg subalgebra of $L$, then one can speak about the distortions $\dist_L^M$ and $\dist_{U(L)}^{U(M)}$. In the case of groups \emph{vs.} group rings, as shown in Proposition \ref{BOpGA}, similarly defined distortions are simply equal. In this subsection we will show first that we always have $\dist_L^M\le\dist_{U(L)}^{U(M)}$ and later that this inequality can be strict. For sharper results we need the following.

\begin{definition}\label{BOdSAF} Given a function $f:\N\to\N$, we say that $f$ is \emph{superadditive} if for any natural $n$ and $m$ we have $f(n)+f(m)\le f(n+m)$. For any function $f:\N\to\N$ we set $\bf(n) = \max_{n_1+\cd+n_t=n}\{ f(n_1)+\cdots+f(n_t)\}$ and call $\bf$ the \emph{superadditive closure of} $f$.
\end{definition}
Clearly, $\bf$ is the least superadditive function that majorates $f$. It is easy to check that if $f\preceq g$ then $\bf\preceq\bg$.  This allows us to correctly define $\overline{\,[\:\! f\,]\,}$ by setting $\overline{\,[\:\! f\,]\,}=[\:\bf\:\,]$. If $B$ is a \fg subalgebra of a \fg algebra $A$ then the distortion $\ds{A}{B}$ is called \emph{superadditive} if $\ds{A}{B}=\overline{\,\ds{A}{B}\,}$.

The main result of this subsection is as follows.

\begin{theorem}\label{BOtSUPER} Let $M$ be a \fg subalgebra of a \fg Lie algebra $L$. Suppose $U(M)$ is the universal enveloping algebra of $M$ naturally embedded in the universal enveloping algebra of $L$. Then 
\[
\ds{U(L)}{U(M)}=\overline{\,\ds{L}{M}\,}.
\]
\end{theorem}

\pp  Let $L$ be a \fg Lie algebra with \emph{linearly independent} generating set $S =\{ a_1,\ld,a_n\}$,   $M$ a \fg subalgebra, with a finite generating set $T$. One can always assume that $T \subset S$. We denote by $\alpha$ the degree filtration of $L$ defined by $S$, $\beta$ the degree filtration of $M$ defined by $T$ and $f=\ds{\alpha\cap M}{\beta}$. By Claim (5) in Proposition \ref{pGPDF} and Claim (i) in Proposition \ref{BOr1} this particular choice of degree filtrations does not affect our conclusion about the connection between $\ds{U(L)}{U(M)}$ and $\dist_{U(L)}^{U(M)}$.

We have $\Sp{T}\subset\Sp{S}$, for the linear spans of $T$ and $S$.  Let us select linear bases in the terms of $\alpha$ and $\beta$, as follows. Suppose that $\alpha$ has this form: $\Sp{S}= L_1 \subset L_2 \subset\ld$ One can choose the pairs of sets $C_n\subset B_n$ such that $B_n$ is a basis of $L_n$, and $C_n$ a basis of $M\cap L_n$, for all $n=1,2,\ld$, by induction on $n$. Then $B=\cup_{i=1}^\infty B_i$ is a basis of $L$, which contains a basis $C=\cup_{i=1}^\infty C_i$ of $M$ as a subset. If we totally order $B$ in some way then $1$ and the monomials $b_1\cdots b_s$ with $b_1\le\ld\le b_s$, called \emph{PBW-monomials}, form a basis of $U(L)$, called \emph{PBW-basis} (PBW stands for Poincar\' e - Birkhoff - Witt). Those of the above PBW-monomials in which all $b_1,\ld,b_s$ are in $C$ form a PBW-basis for $U(M)$ (see the details in \cite[Chapter 1]{B}). Let us assign \emph{weight} $0$ to $1$ and weight $m=m_1+\cdots+m_s$ to each PBW-monomial $u=b_1\cdots b_s$, where $b_i\in B_{m_i}$, $i=1,\ld,s$. Having done so, one can assign weight to each $w\in U(L)$ as the maximum of weights of PBW-monomials in its unique expression with nonzero coefficients through PBW-basis. 

Let us now consider a product of degree $s$ in $U(L)$ of the form $u=b_{1}\cd b_{s}$ in $U(L)$ where $b_i\in B_{m_i}$, for all $i=1,\ld,s$. We claim that one can write $u$ as the linear combination of PBW-monomials of weight at most $m=m_1+\cd+m_s$ with respect to $B$. This is true by definition if $b_{1}\le\cd\le b_{s}$. But if $b_{j}>b_{j+1}$, for some $j$, then, since $\{ L_n\}$ is a filtration, $[b_{j},b_{j+1}]=\sum_{i=1}^t\xi_ib^{\,\prime}_i$, where all $\xi_i\in F$, for all $i$, and each element $b_i^{\,\prime}\in B$ has weight at most $m_j+m_{j+1}$. We then can write 

\begin{eqnarray*}
&&b_{1}\cd b_{j}b_{j+1}\cd  b_{s}=b_{1}\cd b_{j+1}b_{j}\cd b_{s}+b_{1}\cd [b_{j},b_{j+1}]\cd b_{s}\\&=&
b_{1}\cd b_{j+1}b_{j}\cd b_{s}+\sum_{i=1}^t\xi_ib_{1}\cd b_i^{\,\prime}\cd b_{s}.
\end{eqnarray*}

All the monomials on the right hand side of this equation have weight at most $m_1+\cd+m_s$. If all of them are PBW then we are done. Otherwise, we apply the same transformation to those monomials which are not PBW. After finite number of steps, the process terminates and our claim follows.

Let us now consider an associative word $w=y_1\cd y_d$ of length $d$ and its value  $u=s_1\cd s_d$, where each factor is an element of the generating set $S$ of $L$, $S=B_1$. Then, by the above, $u$ can  be written as a linear combination of PBW-monomials in $B$, each having weight at most $d$. Finally, if $g(y_1,\ld,y_n)$ is an associative polynomial of degree $d$ then it follows that $g(s_1,\ld,s_n)$ can also be written as a linear combination of PBW-monomials in $B$ each of weight at most $d$.

Now suppose $g(s_1,\ld,s_n) \in U(M)$, where $g(y_1,\ld,y_n)$ is an associative polynomial of degree $d$. In this case all PBW-monomials in the expression of $g(s_1,\ld,s_n)$ through PBW-basis $B$ of $U(L)$ are actually PBW-monomials in $C$. Suppose, $c_{1}\cdots c_{t}$ is one of these monomials. Then $c_{1}\le\cdots\le c_{t}$ and the sum of the weights of $c_{1},\ldots, c_{t}$ is at most $d$. If $c_{k}$ has weight $w_k$ then we must have $c_{k}\in L_{w_k}\subset M_{f(w_k)}$ where $f$ is the distortion function for $M$ in $L$. Being a Lie polynomial in the generating set $T$ of $M$, $c_{k}$ is an associative polynomial in the generating set $T$ of $U(M)$, of degree $\le  f(w_k)$. Therefore, $c_{1}\cdots c_{t}$ is the value of a polynomial of degree at most $f(w_1)+\cdots+ f(w_t)\le\bf(d)$ with respect to the generating set $T$ of $U(M)$. We have proved that if the degree of an element of $U(M)$ with respect to the generating set $S$ of $U(L)$ is $n$ then its degree with respect to the generating set $T$ of $U(M)$ is at most $\bf(n)$.
 
Now we need to prove that for each natural $n$ there is an element in $U(M)$ whose degree with respect to $S$ is at most $n$ and with respect to $T$ is at least $\bf(n)$. For a fixed $n  > 0$, let us consider an arbitrary partition of $n$ as the sum of positive integral summands $n = n_1+\cd+n_m$. In each $L_{n_i}$ we choose an element $c_i\in C$, which has degree $f(n_i)$ with respect to $T$. Such element exists by the definition of the distortion function $f$. It is allowed to take $c_j=c_i$ if $n_j=n_i$. We reorder the elements $c_1,\ld,c_m$ to make sure $c_i\le c_j$ if  $i<j$. 

Let us consider the PBW-monomial $c_1\cdots c_m$.  By the choice of $c_i\in L_{n_i}$, its degree in $U(L)$ with respect to $S$ is at most $n=n_1+\cd+n_m$. Proving by contradiction, let us assume that the degree of this monomial with respect to $T$ is less than $d =f(n_1)+\cd+f(n_m)$, that is, $c_1\cdots c_s$ is a linear combination of monomials $h(s_1,\ld,s_t)$ in $T$ such that $\deg h(y_1,\ld,y_t) < d$. Let us write $h(s_1,\ld,s_t)$  as a linear combination of PBW-monomials in $C$. As before, the weight of each monomial after rewriting in terms of PBW-basis does not increase, if compared with the sum of the weights of the factors in $h(s_1,\ld,s_t)$. Therefore, the resulting weight will be less than $d$ because the original weight was less than $d$. But the final result of our reductions is our PBW-monomial $c_1\cdots c_m$, whose weight is $f(n_1)+...+f(n_m) = d$, a contradiction. Thus, for a given $n$ we have found an element in $U(M)$ of degree at most $n$ with respect to $S$ whose degree with respect to $T$ is at least $f(n_1)+\cdots+f(n_k)$. Because the partition of $n$ was arbitrary, we have obtained that the distortion function of $U(M)$ in $U(L)$ is at least $\bf(n)$, as claimed.\epf

\subsubsection{Examples of subalgebras with non-superadditive distortion}\label{sSF}

In this excerpt we would like to complement Theorem \ref{BOtSUPER} by exhibiting an example of a Lie algebra $L$ and its subalgebra $M$ such that $\ds{L}{M}\ne\ds{U(L)}{U(M)}$. By Theorem   \ref{BOtSUPER}  we only need to provide an example of a subalgebra $M$ of a Lie algebra $L$ such that the distortion $\ds{L}{M}$ is not superadditive. Actually, such examples can be given not only in the case of Lie algebras but also in the case of associative or Jordan algebras, etc.

Let $\vp,\psi$ be two maps on the set $S=\{ v(0),v(1),\ld\}$,  given by $\vp(v(i))=v(i+1)$, for all $i$, and $\psi(v(j-1))= v(j^2)$ if $j = 2^{2^i}$, for $i =0,1\ld$, otherwise, $\psi(v(j-1)) = v(0)$.  These settings define the action of the free monoid $\cw(a,b)$ on $S$: $a$ acts as $\vp$ and $b$ as $\psi$. Let $S_n  \subset  S$ be the ball of radius $n$ with center $v(0)$ for this action, that is, the set of all elements $w(v(0))$, where $w$ is of length at most $n$ in $\cw(a,b)$. Let $Q_n$ be the ball of radius $n$ with the same center, for the action of the free submonoid $\cw(a)$. The distortion function for the action of $\cw(a)$ with respect to the action of $\cw(a,b)$ by definition equals $f(n) =  \min \{ m\, |\, S_n\subset Q_m \}$. Now $S_{2^{2^i}}$ contains $\psi(v(2^{2^i}-1))= v(2^{2^{i+1}})$, and $v(2^{2^{i+1}})\in Q_{2^{2^{i+1}}}\backslash Q_{2^{2^{i+1}}-1}$. It follows that $f(n) \ge n^2$ for $n = 2^{2^i}$. But if $m  \in [n, n^2-1]$, and $v(m)$ results from $v(0)$ without transformations of the form $v(2^{2^i}-1) \to v(2^{2^{i+1}})$, then it is easy to estimate that we will need to apply $\vp$ and $\psi$ more times compared to the case where this transformation has been applied. (Using this transformation saves $2^{2^{i+1}} -2^{2^i}$ applications of $\vp$ and all other possible applications of $\psi$ save less than $2^{2^i}$. Here we take into account that all $v(i)$ where $i\in[n, n^2-1)$ are mapped by $\psi$ to $v(0)$.) It follows that if $n = 2^{2^i}$ and $m\in [n, n^2-1]$ then the vector with the greatest label in $S_m$ is $v(f(n) + (m-n))$ and so $f(m) = f(n)+m-n\le f(n)+m$.

Now let us assume that $g:\N\to\N$ is a superadditive function and $f$ is equivalent to $g$, that is, $f(n) \le tg(tn)$ and $g(n) \le tf(tn)$ for some positive integer $t$.  We choose $n=2^{2^i}>t^4(t^2+1)$. If $m=t^2(t^2+1)n$ then $n\le m <n^2$, and we can write 
\bea
&&\hspace{-.75cm}(t^2+1)f(n) \le t(t^2+1)g(tn) \le tg((t^2+1)tn)\le t^2f(t^2(t^2+1)n)\\&&=t^2f(m)\le t^2(f(n)+m)
\eqa
and so $f(n)\le t^2m=t^4(t^2+1)n$. Since $f(n)\ge n^2$, we immediately arrive at a contradictory inequality $n\le t^4(t^2+1)$. Thus the distortion function $f(n)$ is not equivalent to a superadditive function.

\begin{example}\label{BOcNSD}
In each of the cases of associative, Lie or Jordan algebras, there exists a \fg algebra $L$ with a \fg subalgebra $M$ such that $\ds{L}{M}$ is not superadditive.
\end{example}

\pp Let us consider a linear space $V$ over a field $F$ with basis $S$ from the above example and extend by linearity the maps $\vp$ and $\psi$ to linear transformations of $V$. Then $V$ naturally becomes a left module over the \fa $\ca=\ca(x,y)$ over $F$ with $x$ acting as $\vp$ and $y$ as $\psi$. The linear space $T=\ca\oplus V$ becomes an associative algebra if we keep the operations of $\ca$, the left action of $\ca$ on $V$ and additionally set $v_1v_2=0$ and $vg=0$, for all and $v_1,v_2,v\in V$ and $g\in\ca$, with free term zero.  One can then make $T$ into a Lie algebra, respectively, Jordan algebra, if one sets $[a,b]=ab-ba$, respectively, $a\cir b=ab+ba$. Since the argument that follows is similar for all three cases, let us consider Lie algebras. Let $L=\cl(x,y)\oplus V$ be a Lie algebra generated by $x,y,v(0)$ and $M=\Sp{x}\oplus V$ a (Lie) subalgebra of $L$ generated by $x,v(0)$ (it is natural to call Lie algebras like $M$ \emph{cyclic}, see also Subsection \ref{ssEDF}). The $m\th$ term of the degree filtration $\beta$ of $M$ defined by $\{ x,v(0)\}$ is $\Sp{x}\oplus \Sp{Q_m}$. The intersection of the $n\th$ term of the degree filtration $\alpha$ of $L$ defined by $\{ x,y,v(0)\}$ with $M$ is $\Sp{x}\oplus \Sp{S_n}$. Therefore, $\ds{\alpha\cap M}{\beta}=f$, the function from our previous example. We know that $f$ is not equivalent to any superadditive function. \epf

This result and Theorem \ref{BOtSUPER} allow one to produce the following.

\begin{example}\label{BOcLAUO} There exist a \fg Lie algebra $L$ and its \fg subalgebra $M$ such that $\ds{U(L)}{U(M)}\ne\ds{L}{M}$.
\end{example}

\pp This follows because one can choose $L$ and $M$ as in Example \ref{BOcNSD} so that $\ds{L}{M}$ is not superadditive and at the same time, by Theorem \ref{BOtSUPER}, $\ds{U(L)}{U(M)}$ is always superadditive.\epf

\section{Realizing tame filtrations as degree filtrations under embeddings}\label{sRDD}

In this section we obtain our central results concerning tame filtrations in both associative and Lie algebras. We show that any tame filtration of a countably dimensional algebra $B$ can be obtained by restriction from a degree filtration of a finitely generated algebra $A$, where $B$ is embedded as a subalgebra. Then we prove that actually $A$ can be chosen simple. 

One general remark about the statements and proofs of our theorems in this and the next section refers to ``if and only if'' claims. As we noted in Introduction, the restriction of a tame filtration of an algebra to a subalgebra is always a tame filtration. Therefore, when we claim that a filtration $\beta$ of an algebra $B$ is a tame filtration if and only if $\beta$ is equivalent (or equal) to the restriction $\alpha\cap B$ of a degree filtration of an algebra $A$ where $B$ is embedded as a subalgebra we actually do not need to proof the ``if'' part. So in each such theorem we will be proving the ``only if'' part, without further comment. 

\subsection{Composition lemmas for associative and Lie algebras}\label{sFMA}

Let $\ca(X)$ be the free associative algebra with free generating set $X$. Given a nonzero polynomial $f\in\ca(X)$, the highest word, with respect to Shortlex, which enters $f$ with nonzero coefficient is called the \emph{leading term} or the \emph{leading word} of $f$ and is denoted by $\bf$. We call $f$ \emph{monic} if the coefficient of $\bf$ in $f$ is 1.

In the definitions and lemma that follow, the total order on $\mathcal{W}(X)$ does not need to be Shortlex, however it must be a semigroup order, that is, $u\le v$ must imply $uw\le vw$ and $wu\le wv$, for any word $u,v,w\in \mathcal{W}(X)$.  The remainder of this subsection can be traced back to \cite{BC} (see also \cite{BK}).

\begin{Definition}\label{dBCA} Let $\mathcal{W}(X)$ be given a semigroup total ordering, $f$ and $g$ two monic polynomials in $\ca(X)$, with leading words $\bf$ and $\bg$ . Then, there are two kinds of \emph{compositions} of $f$ and $g$:
\begin{enumerate}
\item[\emph{(i)}] If $w$ is a word such that $w = \bf a = b\bg$ for some $a, b\in \mathcal{W}(X)$ with $\deg(\bf)+\deg(\bg) >\deg(w)$,
then the polynomial $(f, g)_w = fa- bg$, is called the \emph{intersection composition} of $f$ and $g$, with respect to $w$.
\item[\emph{(ii)}] If $w = \bf = a\bg b$ for some $a, b \in \mathcal{W}(X)$, then the polynomial $(f, g)_w = f - agb$ is called
the \emph{inclusion composition} of $f$ and $g$, with respect to $w$.
\end{enumerate}
\end{Definition}
Notice that the leading word of $(f, g)_w$ is strictly less than the leading word of $fa$ or $gb$ in the first case and of $f$ in the second.
\begin{Definition}\label{dCC} A set of monic polynomials $S$ is \emph{closed under compositions} or is a \emph{Groebner-Shirshov set}
if any composition $(f, g)_w$ of elements $f,g\in S$ can be written in $\ca(X)$ as
$(f, g)_w =\sum \alpha_ia_if_ib_i$, where $f_i\in S$, $\alpha_i\in F$, and the leading word $a_i\bf_ib_i$ of each $a_if_ib_i$ is strictly less than  $w$.
\end{Definition}
\begin{cla}\label{CLA} Let $S$ be a subset of the free associative algebra $\ca(X)$ over a field $F$, consisting of monic polynomials, $I$ a two-sided ideal generated by $S$, $A = \ca(X)/I$. Then
the following conditions are equivalent:
\begin{enumerate}
\item[\emph{(i)}] $S$ is a Groebner-Shirshov set.
\item[\emph{(ii)}] for any element $f\in I$ we have $\bf = a\bs b$ for some $s \in S$ and $a, b \in \mathcal{W}(X)$.
\item[\emph{(iii)}] The set of cosets represented by reduced words $\{ u+I\,|\, u\in \mathcal{W}(X),\, u\ne a\bs b, s\in S, a,b\in \mathcal{W}(X)\}$ is a basis of $A$.\epf
\end{enumerate}
\end{cla}

In the case of Lie algebras the notion of a set of Lie polynomials closed under composition is defined as follows. We always view a free Lie algebra $\mathcal{L}(X)$ as a Lie subalgebra in $\ca(X)$ generated by $X$ under the bracket operation $[a,b]=ab-ba$. The elements of $\lax$ are called \emph{Lie polynomials} in $X$.  The leading word $\bf$  of a Lie polynomial $f\in\mathcal{L}(X)$ is defined because $f$ is an element of $\ca(X)$. Also we call a Lie polynomial \emph{monic} if it is monic as an associative polynomial. 

The elements of $\ca(X)$ resulting from $X$ by repeated application of the bracket operation only are often called (higher) \emph{commutators}. Forgetting brackets on a commutator $c$ produces a unique word $w$, which is called the (associative) \emph{carrier} of $c$. 

If $X$ is totally ordered then there is so called \emph{Shirshov order} $\preceq$ on $\wax$ which is lexicographic unless we compare a word with its proper prefix; in that case the prefix is greater than the word. This order induces an order on the set of commutators: $[u]\preceq[v]$ if and only if $u\preceq v$. With each word $w\in\wax$ one can associate the set $C(w)$ of its cyclic permutations. If $w$ is greater than any other word in $C(w)$, we call $w$ \emph{regular}. If $f\in\cl(X)$ then $\bf$ is always regular. In particular, if $c$ is a commutator then $\bc$ is regular. 

A special set of commutators forms a basis of the vector space $\cl(X)$ over $F$. Its elements, called \emph{basic commutators}, are constructed using induction by degree. The elements of $X$ are basic commutators of degree 1. Suppose we have defined basic commutators of degree $<n$. Let $c,d$ be any basic commutators of degrees $k,n-k<n$, satisfying
\begin{enumerate} 
\item[(i)] $d\prec c$
\item[(ii)] if $c=[c_1, c_2]$ then $c_2\preceq d$. 
\end{enumerate}
In this case, we declare $[c,d]$ a basic commutator of degree $n$. For the basic commutators, the leading word and the associative carrier are the same thing, and the map $c\to\bc$ is the bijection between the set of basic commutators and the set of regular words. In particular, if $[c,d]$ is a basic commutator then $\overline{[c,d]}=\bc\bd$. The existence of inverse to the above map means that on any regular word one can set brackets so that we obtain a basic commutator.

To formulate an important Special Bracketing Lemma by Shirshov (see \cite[Lemma 3.10]{BC}), we need the following notation. Let $u$ and $v$ be regular associative words $u = avb$, $a,b\in \mathcal{W}(X)$. Suppose we have turned $u$ in a basic commutator by setting brackets in such a way that one pair of matching brackets embraces a subword $vc$ of $u$ where $b=cd$. Symbolically, we write  $[u] = [a[vc]d]$. Let us write $c=c_1\cdots c_m$ where all $c_i$ are regular and $c_1\preceq\ldots\preceq c_m$. Afterwards, we set brackets in a unique way on $v$ and each $c_i$ to obtain basic commutators $[v],[c_1],\ld,[c_m]$. Next we form the (left-normed) commutator $w=[\ld[[v],[c_1]],\ldots,[c_m]]$ and replace $[vc]$ by $w$ in $[u]$. The resulting commutator will be denoted by $\kappa(a,v,b)$. If now $g$ is a monic Lie polynomial with leading word $v$ included as above in a regular associative word $u$ then we proceed exactly as before but when we form $w$ we replace $[v]$ by $g$. The resulting Lie polynomial will be denoted by $\sigma(a,g,b)$. 

If we use this notation then the following is true.

\begin{cbl}\label{lS} Let $u$ and $v$ be two regular associative words such that $u = avb$, $a,b\in \mathcal{W}(X)$.
 \begin{enumerate}
\item[\emph{(i)}] In the unique setting of brackets on $u$ which makes it into a basic commutator two pairs of matching brackets are set as follows: $[u] = [a[vc]d]$, where $b = cd$, $c, d \in \mathcal{W}(X)$.
\item[\emph{(ii)}] If we use the bracketing of \emph{(i)} to form the commutator $\sigma(a,[v],b)$, as described before this lemma, then
$\overline{\sigma(a,[v],b)} = u$.
\epf
\end{enumerate}
\end{cbl}

We will use the above notation in the definitions and lemma that follow.

\begin{Definition}\label{dLC} Let $f$ and $g$ be two monic Lie polynomials in $\mathcal{L}(X)\subset\ca(X)$. Then there are two kinds of Lie compositions:
\begin{enumerate}
\item[\emph{(i)}] If $w = \bf = a\bg b$ for some $a, b \in \mathcal{W}(X)$, then the polynomial $\langle f, g\rangle_w = f -\sigma(a,g,b)$ is called
the \emph{composition of inclusion} of $f$ and $g$ with respect to $w$.
\item[\emph{(ii)}] If $w$ is a word such that $w = \bf b = a\bg$ (then $w$ has to be regular!), for some $a, b \in \mathcal{W}(X)$, with $\deg(\bf)+\deg(\bg) >\deg(w)$,
then the polynomial $\langle f, g\rangle_w = \sigma(1,f,b) -\sigma(a,g,1)$ is called the \emph{composition of intersection}
of $f$ and $g$ with respect to $w$.
\end{enumerate}
\end{Definition}

Notice that the leading word of $\langle f, g\rangle_w$ is strictly less than that of $f$ in the first case and each of $f b$ or $a g$ in the second.

\begin{Definition}\label{dLCC} A nonempty set $S\subset \mathcal{L}(X)$ of monic Lie polynomials is \emph{closed under composition} or is a \emph{Groebner-Shirshov set} if any composition $h=\langle f, g\rangle_w$ of $f,g\in S$ with respect to $w$ can be written as $h=\sum_i\alpha_i\sigma(a_i,s_i,b_i)$ where $s_i\in S$, $\alpha_i\in F$, $a_i,b_i\in \mathcal{W}(X)$ and $a_i\bs_ib_i<w$, for all $i$.
\end{Definition}
\begin{Proposition}\label{pBM} A set $S\subset \mathcal{L}(X)\subset\ca(X)$ is Groebner-Shirshov in $\mathcal{L}(X)$ if and only if it is such in $\ca(X)$.
\end{Proposition} 

\begin{cll}\label{tSh} Let $S \subset \mathcal{L}(X) \subset\ca(X)$ be a nonempty set of Lie polynomials in $X$. Let $I$ be the Lie ideal  generated by $S$ in $\mathcal{L}(X)$ and $J$ the two-sided associative ideal  generated by $S$ in $\ca(X)$. Then the following conditions are equivalent.
\begin{enumerate}
\item[\emph{(i)}] $S$ is a Groebner-Shirshov set in $\mathcal{L}(X)$.
\item[\emph{(ii)}] for any elements $f\in J$ we have $\bf = a\bs b$ for some $s \in S$ and $a, b \in \mathcal{W}(X)$.
\item[\emph{(iii)}] The cosets  $[u]+I$, such that $u$ is a regular word in $\mathcal{W}(X)$, without subwords $\bs$ where $s\in S$, form a basis in factor-algebra $\mathcal{L}(X)/I$. 
\item[\emph{(iv)}] The cosets  $u+J$, $u$ a word in $\mathcal{W}(X)$ without subwords $\bs$ where $s\in S$, form a basis in factor-algebra $\mathcal{A}(X)/J$.\epf
\end{enumerate}
\end{cll}

\subsection{Tame filtrations in associative algebras}\label{sCA}

Let $\cb$ be a subalgebra in a \fa $\ca=\ca(X)$, $\# X>1$, generated by a set $\mathcal{M}$ of words such that no nonempty suffix of any word in $\mathcal{M}$ is a prefix of another word in $\mathcal{M}$ and also none of the words in $\mathcal{M}$ is a subword of another word in $\mathcal{M}$. We call this condition ``non-overlapping''. Such sets exist and, moreover, one can choose $\mathcal{M}$ so that the growth of the number of elements in the set $\cm_n$ of words of degree $\le n$ in $\cm$ is a function that majorates an exponential function $c^n$, with $c>1$, for all sufficiently large $n$ (``exponential sets''). As an example of such set in the case where there are only two variables, one can consider the set of all words $x^3ywxy^3$, where $ywx$ has no subwords $x^3$ or $y^3$. But if we do not restrict the number of variables, then, given any natural $c$, we can produce a non-overlapping set $\cm$ with $\#\cm_n\ge c^n$, for all $n\ge 1$ if we proceed in the following way. We choose $X$ with $\#X=c+2$, select two letter $x$ and $y$ and consider the set $\cm$ of all words $xwy$ such that $w$ does not depend  on $x,y$. For this set $\cm$, we would have $\#\cm_n\ge c^n$, for any $n>2$. But since we also  want to have $\#\cm_1\ge c$ and $\#\cm_2\ge c^2$ we set $d=c^2$, and add to $X$ new variables $z_1,\ld,z_c$, $z_1^{(1)},\ld,z_{d}^{(1)}$, $u_1^{(2)},\ld,u_{d}^{(2)}$ and to $\cm$ new words of length one: $z_1,\ld,z_c$ and of length two: $z_1^{(1)}z_1^{(2)},\ld,z_{d}^{(1)}z_{d}^{(2)}$.

Notice that from the ``non-overlapping'' property of $\cm$, it is immediate that the set $\cu$ of products of words in $\cm$ (including the empty word 1) is a free submonoid $\cu\cong\cw(\cm)$ and the linear span $\cb$ of $\cu$ is the free unital associative algebra $\cb\cong\mathcal{A}(\cm)$. Two easy properties of $\cu$ are as follows. If $u\in\cu$  and $u = avb$, where also $v\in\cu$, then $a,b\in\cu$. Also if $u,v,w\in\cw(X)$, $w\ne 1$ and $vw,wu\in\cu$ then $u, v, w\in\cu$. Another remark is that any total order on $\cw(X)$ induces a total order on $\cu$.
 
  Now notice that given a total semigroup ordering of the words in a free associative algebra, any ideal $I$ has a Groebner - Shirshov basis. This is simply any basis $\cs$ of $I$ consisting of monic polynomials and such that the leading words of different elements of the basis are different. So if $I$ is an ideal in $\cb=\ca(\cm)$, with induced order from $\cw$, then there is Groebner - Shirshov basis $\cs$ of $I$, as an ideal of $\cb$. 
  
\begin{lemma}\label{BOlC} Let $\cm$ be a non-overlapping set in the free associative algebra $\ca(X)$ of rank $\ge 2$, $\cb=\alg\cm$, $I$ an ideal of $\cb$, $\cs$ a vector space basis of $I$ consisting of monic polynomials such that the leading words of different elements are different. Consider the ideal $J=\mathrm{id}_{\ca} I$ of $\ca$ generated by $I$. Then the leading word of each element of $J$ contains the leading word of an element in $\cs$, as a subword.
\end{lemma}

\pp  We would easily derive this by Composition Lemma for Associative Algebras (Subsection \ref{sFMA}) provided that we have checked that the set $\cs$ is closed under composition. 
For the intersection composition, we need to take two polynomials $f, g\in \cs$ and assume that there are words $a,b\in\cw(X)$ such that $w = \bf a = b\bg$. Since $\bf,\bg\in\cu$, by one of the above mentioned properties of $\mathcal{M}$, it follows that $a,b\in\cu\subset\cb$. Then both $fa$ and $bg$ are in $I$. Each can be written as a linear combination of elements of $\cs$ with leading word $w$ in both cases. As a result, $(f,g)_w=\sum \alpha_if_i$ where each $f_i$ is in $\cs$ and $f_i<w$, for all $i$. Therefore, $\cs$ is closed under the intersection composition. For the inclusion composition, suppose that $w=\bf=a\bg b$, for $a,b\in\cw(X)$. Then, as before, $a,b\in\cu\subset\cb$ and so both $f$ and $agb$ are in $I$. The same argument, as just before, shows that $\cs$ is closed also under the inclusion composition $(f,g)_w=f-agb$. Now the conditions of Composition Lemma for Associative Algebras are satisfied. As a result, the leading word $\bf$ of every polynomial $f\in J$ contains as a subword a leading word $\bg$ of a polynomial $g \in  \cs\subset I$. \epf 

Adopting the notation of the previous lemma, let us assume that $f \in  J\cap\cb$. In this case, $\bf$  is the product of words in $\mathcal{M}$ and also by this lemma, $\bf  = a \bar g b$, $g\in\cs$. Since $\bg$ is also a product of words in $\mathcal{M}$, it follows that $a$ and $b$ are products of words in $\mathcal{M}$. Finally, $agb \in I$ and then $f-agb\in J\cap\cb$ with a lesser leading word. Using induction by leading word shows that $f \in I$. We now can conclude that $J \cap \cb = I$.  Actually, this result is a particular case (when $f\in \cb\cap J$) of the following.
 
\begin{lemma}\label{lcep} Let $f$ be a polynomial in $\cb$. Then in the coset $f+I$ there is a polynomial $f_0$, such that $\overline{f_0} \le \bh$, for any polynomial $h$ in the coset $f+J$. In particular, $J \cap \cb = I$.
\end{lemma}
 
\pp Let $f_0$ be a polynomial in $f+I$ with leading word minimal possible in $\ca$. Proving by contradiction, let us assume that there is $h\in f + J$, such that $\bh < \overline{f_0}$. Let us apply Lemma \ref{BOlC} and Composition Lemma for Associative Algebras to $f_0 - h$, which is a nonzero element of $J$. We have $\overline{f_0 - h} =a \bg b$, where $g \in I$.  But $\overline{f_0 - h} = \overline{f_0} \in \cb$, hence, as earlier, $a, b$ are products of elements in $\mathcal{M}$ and then $agb \in I$.  It follows that $f_0$ can be replaced by a difference $f_0 - agb$ whose leading word is lower, in contradiction with the choice of $f_0$. \epf

One of the simplest non-overlapping sets of words in the variables $\{ x,y,z\}$ is the set $\cm=\{ xy^iz\,|\, i=0,1,2,\ld\}$. This set can be used  to prove that \emph{any countably-dimensional algebra $B$ can be embedded in a finitely generated algebra} (A.I.Malcev's result, see \cite{AIM}). Indeed, let $b_0,b_1,b_2,\ld$ be a subset of $B$ that generates $B$. Consider $\cb=\alg{\cm}\subset \ca(x,y,z)$. As noted above, $\cb$ is a free associative algebra with free generating set $\cm$. Then the map $xy^iz\mapsto b_i$, $i=0,1,2,\ld$, extends to an epimorphism $\nu:\cb\to B$. Let $I$ be the kernel of $\nu$, $J$ the two-sided ideal of $\ca$ generated by $I$, $A=\ca/J$ and $\mu$ the natural epimorphism $\mu:\ca\to\ca/J$. By Lemma \ref{lcep}, we have $J\cap\cb=I$ and so the well-defined map $\vp:B\to A$ given by $\vp(b)=\mu(\nu^{-1}(b))$ is the desired embedding of $B$ in a 3-generator algebra $A$. 

For our argument in Section \ref{sHEUA}, we need a modification of this general result, as follows.

\begin{proposition}\label{BOrABC} Suppose $B$ is an arbitrary unital countable $F$-algebra with a finite generating set $S$  and $b_1,b_2,\ld$
is an arbitrary enumeration of all of its elements (each element may occur infinitely many times) such that $\deg_Sb_i\le i$, for all $i=1,2,\ld$ Then there exists a unital \fg $F$-algebra $A$ containing
$B$ as an undistorted subalgebra and elements $a,b,c\in A$, such that $1=ac$ and $b_i=ab^ic\;(i=1,2,\ld)$. 
\end{proposition}

\pp The argument preceding the statement of this proposition applies, with $b_0=1$. If we select in $A=\ca/J$ three elements $a=\mu(x)$, $b=\mu(y)$, $c=\mu(z)$ and identify $B$ with its image under emebedding $\vp$ then we will have $ac=1$ and $ab^ic=b_i\in A$, for all $i=1,2,\ld$, as needed.  

To prove the unistortedness of this embedding, we recall the generating set $S$ in $B$ and choose the generating set $T=\{ a,b,c\}$ for $A$. Let $\beta=\{ B_n\}$ be the degree filtration of $B$ defined by $S$ and $\alpha=\{ A_n\}$ the degree filtration of $A$ defined by $T$. We select $u\in B$ and assume $u\in A_n$. In this case, there is $f(x,y,z)\in \ca$ of degree $n$ such that $u=f(a,b,c)$. By Lemma \ref{lcep} then there is $f_0(x,y,z)\in\cb$ whose degree with respect to $\{ x,y,z\}$ is at most $n$ and such that $f_0(x,y,z)+I=u$. In this case we can rewrite $f_0(x,y,z)$ in terms of the free generating set $\cm$ of $\cb$: $f_0(x,y,z)=g(xz,xyz,\ld,xy^kz)$. Each monomial of $g$ has the form $xy^{i_1}z\cdots xy^{i_m}z$ and $i_1+\cdots+i_m+2m\le n$. Each such monomial will map under $\vp$ to $b_{i_1}\cdots b_{i_m}$ and this element has degree at most 
$i_1+\cdots+i_m\le n$ with respect to the generating set $S$ of $B$. Hence $u=\vp(g)\in B_n$. Thus, $\ds{\alpha\cap B}{\beta}(n)\le n$, and the embedding has no distortion. \epf

We now prove the main result of this subsection.

\begin{theorem}\label{BOpDdeg} 
Let $B$ be a unital associative algebra over a field $F$.
\begin{enumerate} 
\item[$\mathrm(1)$] A filtration $\beta$ on $B$ is tame if and only if $\beta\sim\alpha\cap B$ where $\alpha$ is a degree filtration on a unital 2-generator associative algebra $A$ where $B$ is embedded as a subalgebra.
\item[$\mathrm(2)$] A filtration $\beta$ on $B$ is tame if and only if $\beta=\alpha\cap B$ where $\alpha$ is a degree filtration on a unital \fg associative algebra $A$ where $B$ is embedded as a subalgebra.
\end{enumerate}
\end{theorem}

\pp If $B$ is a subalgebra in a \fg algebra $A$, and $\beta\sim\alpha\cap B$ where $\alpha$ is a degree filtration of $A$ then, as noted in Introduction just after Definition \ref{dDFil}, $\beta$ is tame. This proves the ``if'' parts in Claims (1) and (2).

Let us now prove the ``only if'' part in Claim (1). Suppose that $\beta=\{ B_n\}$ is a tame filtration on a unital associative algebra $B$. We will use the notation of the first two paragraphs of the current subsection. Notice that the subalgebra $\cb$ defined earlier is a free graded subalgebra in a graded algebra $\ca$ and that $\mathcal{M}$ is a free graded generating set of $\cb$. We write $\ca=\bigoplus_{m=0}^\infty \ca^{(m)}$, where $\ca^{(m)}$ is the linear span of words of length $m$ in $X$. Since $\cb$ is generated by words, we have $\cb=\bigoplus_{m=0}^\infty \cb^{(m)}$ where $\cb^{(m)}=\cb\cap \ca^{(m)}$. These gradings induce filtrations: the degree filtration with terms $\ca_n=\bigoplus_{m=0}^n \ca^{(m)}$ on $\ca$ and a tame filtration with terms $\cb_n=\bigoplus_{m=0}^n \cb^{(m)}$ on $\cb$. We also set $\mathcal{M}^{(m)}=\mathcal{M}\cap \ca^{(m)}$. 

Since the growth of the sequence $\{ \# \mathcal{M}^{(n)}\}$ majorates an exponential function, there is a positive integral constant $C$ such that $\dim B_n \le \#\mathcal{M}^{(Cn)} \le \dim  \cb_{Cn}$, for any $n$. This allows us to define an epimorphism $\nu: \cb\to B$ so that for each $n=1,2,\ld,$ $\Sp{\nu(\mathcal{M}^{(Cn)})}=B_n$.  We also may assume that $\nu(\mathcal{M}^{(k)})=0$, for the values of $k$ not divisible by $C$. Since $\mathcal{M}^{(Cn)}\subset\cb_{Cn}$ we have $B_n\subset\nu(\cb_{Cn})$. The converse inclusion is also easy. Indeed,  
\bea
\nu(\cb_{Cn})&=&\Sp{\sum_{l_1+\ld+l_s=Cn}\nu(\mathcal{M}^{(l_1)})\cdots\nu(\mathcal{M}^{(l_s)})}\\
&=& \Sp{\sum_{k_1+\ld+k_t=n}\nu(\mathcal{M}^{(Ck_1)})\cdots
\nu(\mathcal{M}^{(Ck_t)})}\\
&\subset&\sum_{k_1+\ld+k_t=n}B_{k_1}\cdots B_{k_t}\subset B_n.
\eqa
So we have $B_n\,=\,\nu(\cb_{Cn})$.

Now let $I$ be the kernel of $\nu$. If $J$ is the ideal of $\ca$ generated by $I$ then by Lemma \ref{lcep}, $\cb\cap J=I$. Let us set $A=\ca/J$ and suppose that $\mu$ is the canonical epimorphism from $\ca$ to $A$. The degree filtration of $\ca$ induces the degree filtration on $A$ given by $A_n=(\ca_n+J)/J$, for $n=0,1,2,\ld$. The natural embedding $\vp:B\to A$ is given by $\vp(u)=\mu(\nu^{-1}(u))$. We identify $B$ with its image in $A$, using $\vp$. To complete the proof of our theorem, it is sufficient to show that $B_n = B\cap A_{Cn}$. 

Since $\cb_{Cn}\subset \ca_{Cn}$ and for any $u\in B_n$ we have $u=\nu(b)=b+I$, where $b\in \cb_{Cn}\subset \ca_{Cn}$, we have $\vp(u) =\widetilde{\nu}(\nu^{-1}(u) = \mu(b) = b+J \in A_{Cn}$. Thus  $B_n\subset B\cap A_{Cn}$. 

To prove the inverse inclusion, we pick $\mu(a)\in B\cap A_{Cn}$. Then $a$ can be chosen in $\ca_{Cn}$ and there is $u\in B$ such that $\mu(a)= \vp(u) =  \mu(\nu^{-1}(u)$. If $u=\nu(b)$, for some $b\in \cb$, then $a+J=b+J$. Using Lemma \ref{lcep}, we can find $b_0\in b+I$ such that $b_0\in \cb_{Cn}$. Then $\mu(a)=\mu(b_0) =\nu(b_0)\in \nu(\cb_{Cn})\subset B_n$, and so $B\cap A_{Cn}\subset B_n$, as claimed. Now the proof of Claim (1) is complete.

To prove the ``only if'' part in Claim (2), notice that the previous proof works for any number of elements in the set $X$, $\#X\ge 2$. But if  $\# X$ is big enough, as specified in the first paragraph of this subsection, then the growth of the set $\cm$ can be made faster than any exponential function. Consequently, the constant $C$ that was used to cover $B_n$ by $\cb_{Cn}$ can be take equal 1. In this case we get $B_n=B\cap A_n$ and so the original tame filtration on $B$ is simply the restriction of a degree filtration on an appropriate \fg algebra $A$. This complete the proof of Claim (2), hence of the whole theorem. \epf

An argument very similar to the one used in the proof Theorem \ref{BOpDdeg}, allows one to derive the following result about the monoids.

\begin{theorem}\label{pUEM} Let $N$ be a countable monoid.
\begin{enumerate} 
\item[\emph{(i)}] There exist a 3-generator monoid $M$ where $N$ is embedded as a submonoid. 
\item[\emph{(ii)}] If $N$ is  finitely generated then the embedding of $N$ in a 3-generator monoid $M$ can be done without distortion.
\item[\emph{(iii)}] A filtration $\beta$ on $N$ is a tame filtration if and only if there is a \fg monoid $M$ with a degree filtration $\alpha$ such that $\beta\sim\alpha\cap N$. 
\end{enumerate}
\end{theorem}

\pp All claims follow if we replace the occurrences of the word ``monoid'' by the words ``monoid algebra''. Notice that \fa $\ca=F[\cw(X)]$, its free subalgebra $\cb=F[\cw(\cm)]$ are monoid algebras. 

To prove Claim (i) we have to use the argument preceding Proposition \ref{BOpGA}. If the elements $\{b_0,b_1,\ld\}$ are the elements in $N$ then $\nu:\cb\to B$ to $\cw(\cm)$ is an epimorphism of monoids $\cw(\cm)$ to $N$. The kernel $I$ of $\nu$ is then generated by the elements of the form $u-v$ where $u,v$ are some elements of $\cm$. The same elements generate $J$ and thus $A=\ca/J$ is a monoid algebra of the 3-generator monoid $M=\mu(\cw(X))$. The restriction of $\psi$ to $N$ is an embedding of monoids from $N$ to $M$, proving Claim (i). 

To prove Claim (ii) we have to apply Proposition \ref{BOrABC} in conjunction with Proposition \ref{BOpGA}, which relates the distortion of the embedding of monoids to that of respective monoid algebras.

To prove the Claim (iii), we only need to note that tame filtrations on $N$ come from a tame filtrations on $B=F[N]$ and the degree filtration on $M$ comes from the degree filtration on $\mathcal{W}(X)$ defined by $X$. Thus, these filtrations are merely the restrictions of the tame filtration of $B$ to $N$ and the degree filtration of $A$ to $M$, and the number of elements in the $n\th$ term of each of these filtrations for $N$ or $M$ is the dimension of the the respective filtrations for $B$ and $A$. As a results, our claim about tame filtration follows from Theorem \ref{BOpDdeg}. \epf

   
 
\subsection{Tame filtrations in Lie algebras}\label{sCLA}

Let $X$ be a totally ordered finite set and $\ca(X)$ the free associative algebra with $X$ as the set of free generators. Let $x,y\in X$ and assume that $X$ is ordered in such a way that $x$ is greater than any other letter in $X$. This ordering expands to one of the total orderings on $\cw(X)$: Shortlex or Shirshov ordering described in Subsection \ref{sFMA}. As in the previous subsection, we will be using here an exponentially growing non-overlapping set $\cm$  in one of two forms: $x^3ywxy^3$, where $ywx$ has no subwords $x^3$ or $y^3$, and also, provided that $\# X>3$, the set of all word $xuy$, where $u\in\cw(X)$ is a word without letters $x$ and $y$. Notice that in either case $\cm$ consists of regular words in the sense of Subsection \ref{sFMA}. Some of the properties of the words in the free semigroup $\cw(\cm)$ have been indicated in the first paragraphs of the previous subsection.

Now let $\mathcal{N}$ be the set of basic commutators obtained by setting brackets on the elements of $\mathcal{M}$, $\cb$ the subalgebra in the \fa $\ca$ generated by $\mathcal{M}$, and  $\cc$ the Lie subalgebra in the free Lie algebra $\cl=\lax$ generated by $\mathcal{N}$. As in the previous subsection, we have free monoid $\cu\cong\cw(\cm)$, the free associative algebra $\cb\cong\ca(\cm)$ and now the free Lie algebra $\cc\cong\cl(\cn)$.

\begin{lemma}\label{BOlMBN} Let $u\in \cw(\cm)$ be a regular associative word in $\cw(X)$. If we set brackets on $u$ in a unique way so that the resulting word is a basic commutator $[u]$ in $\cl(X)$ then $[u]\in\cl(\cn)$.
\end{lemma}

\pp Let us write $w=w_1\cdots w_m$, where each $w_i$ is an element of $\cm$. In the case where $\cm$ is the set of the first kind, the argument is as follows. We apply Shirshov's method of bracketing (see \cite[Lemma 3.10]{BC}). It consists of subsequent iterations during each of which we find the smallest word $v$ among those already bracketed (as noted before Special Bracketing Lemma in Subsection \ref{sFMA}, the words of length 1 are already bracketed!). For each already bracketed subword $u\ne v$ we look at the maximal number $l$ of bracketed subwords equal $v$, which follow $u$, and introduce the bracketed subwords of the next level by setting brackets on $uv^l$ in a left-normed way to obtain bracketing $[[u,\underbrace{v],v,\ldots v}_{l}]$. Suppose that after $t$ iterations of the process the brackets respect the decomposition $w=w_1\cd w_m$. Let us assume to the contrary that after the next iteration, the brackets have been set to give the commutator $[[u,\underbrace{v],v,\ldots v}_{l}]$, whose associative carrier contains a suffix of $w_k$ and a prefix of $w_{k+1}$. Since the maximal letter $x$ cannot be the last letter of the regular word $v$, the subword $v\cd v$ must contain the prefix $x^3$ of $w_{k+1}$. Since $uv\cd v$ is regular, it follows that $x^3$ is a prefix of $u$. Now since $v\cdots v$ contains $x^3$ but does not end by $x$, it follows that $v$ itself must contain $x^3$. But $v$ is regular and so $x^3$ must be the prefix of $v$. As a result,  $u=w_k$ and $v$ is a prefix of $w_{k+1}$. But Shirshov's bracketing assumes that $v$ is the smallest of the words already bracketed in course of the first $t$ iterations (but before the $(t+1)^{\mathrm{st}}$ begins!). This means that after $t$ steps, our word is of the form $[v_1]\cd[v_t]$, where each $v_i$ is not less than $v$, hence $x^3$ is the prefix of $v_i$. But then $v_1=w_1, v_2=w_2,\ld$, which is what we want.

In the case of $\cm$ of the second kind the argument goes through if we replace $x^3$ by $x$. \epf

This Lemma and its proof allow one to modify Special Bracketing Lemma (Subsection \ref{sFMA}).

\begin{corollary}\label{BOlMBS} Let $u,v\in \cw(\cm)$ be regular associative words, $u=avb$ where $a,b \in \cw(X)$. Suppose $g\in\cl(\cn)$ has leading word $v$. Then $\sigma(a,g,b)$ is an element of the Lie ideal of $\cl(\cn)$ generated by $g$.
\end{corollary}

\pp First of all, by the properties of the set $\cm$ we have that $a,b\in\cm$. Using the same method of choosing the smallest letter, etc., as in the previous lemma, we can see that throughout the process of  forming $\kappa(a,v,b)=[a[[v],[c_1],\ld,[c_n]]d]$ all commutators arising are in $\cl(\cn)$ and then when we replace $[[v],[c_1],\ld,[c_n]]$ by $g$ we obtain the desired property of $\sigma(a,g,b)$.\epf 

We now denote by $I$ an ideal of Lie algebra $\cc$, $J$ the two-sided (associative) ideal of $\ca$ generated by $I$, and $K$ the ideal of Lie algebra $\cl$, generated by $I$. Using Zorn's Lemma, one can choose  a  linear basis $\mS$ of $I$ in such a way that different elements of $\mS$ have different leading words. Each element in $\mS$ is a linear combination of basic commutators and  so these leading words are regular associative words.

\begin{lemma}\label{lLT}  Let $I$, $J$ and $K$ be ideals in $\cc$, $\ca$ and $\cl$, respectively, as just defined, and $f\in J$ or $f\in K$, $f\ne 0$. Then the leading word $\bf$ has a subword equal to the leading word of an element in $\cs$.
\end{lemma} 

\pp Since $K\subset J$, we only need to consider the case where $f\in J$. Then, by Composition Lemma for Associative Algebras (Subsection \ref{sFMA}), we need to check that $\cs$  is closed under associative compositions. In view of Proposition GS (Subsection \ref{sFMA}), it is sufficient to check this claim for Lie compositions of two elements in $\cs$.

First, let us consider the inclusion composition $\langle f,g\rangle_w$ for two elements $f,g\in \cs$ where $w=\bf=u\bg v$, for some $u,v\in\cw(X)$. Since $\bf\in\cw(\cm)$, it follows that in this case also $u,v\in\cw(\cm)$. By Corollary \ref{BOlMBS}, $\sigma(u,g, v)$ is an element of $I$. As such, $\sigma(u,g, v)$ can be written as a linear combination of elements of $\cs$, with a leading basic commutator $c$. The leading basic word of $\sigma(u,g, v)$ is then $\bc$ and hence by Special Bracketing Lemma (Subsection \ref{sFMA}), $\bc=u\bg v=\bf$. It follows then that the leading basic commutator in the expression of $f$ is also $c$ and then $\langle f,g\rangle_w=\sigma(u,g, v)-f$ is a linear combination of elements in $\cs$, strictly less than $f$.   


Second, let us consider the composition of intersection $\langle f,g\rangle_w$ for two elements $f,g\in I$ where $w=\bf u=v\bg$, for some $u,v\in\cw(X)$. Again, using the same argument, $u,v\in\cw(\cm)$. By Corollary \ref{BOlMBS}, $\sigma(1,f, u)$ is in the ideal of $\cl(\cn)$ generated by $f$ and $\sigma(v,g,1)$ is in the ideal of $\cl(\cn)$ generated by $g$. Each of those ideals is in $I$. As in the previous case, then each of $\sigma(1,f, u)$, $\sigma(v,g,1)$ is a linear combination of elements of $\cs$. By Special Bracketing Lemma (Subsection \ref{sFMA}) the leading words of both Lie polynomials are the same. When we form the composition $\langle f,g\rangle_w=\sigma(1,f, u)-\sigma(v,g,1)$, we observe that this is a linear combination of elements of $\cs$ whose leading words are strictly smaller than $w$. 

To comply with Definition \ref{dLCC} (Subsection \ref{sFMA}), let us notice that for any $s\in \cs$ we always have $\sigma(1,s,1)=s$. In that case, for both types of composition, we have 

\[
\langle f,g\rangle_w=\sum\alpha_is_i=\sum\alpha_i\sigma(1,s_i,1)
\]

\noindent where each $s_i\in\cs$ satisfies $\overline{s_i}<\bs$. In this case also $1\cdot \overline{s_i}\cdot 1<\bs$ and  by Definition \ref{dLCC}, we have that $\cs$ is closed under composition.

Now we conclude that $\cs$ is a Groebner-Shirshov basis for $J$ and by Composition Lemma for Associative Algebras (Subsection \ref{sFMA}) our claim follows. \epf

The next lemma is an analogue of Lemma \ref{lcep}. We will need it to prove a Lie analogue of Theorem \ref{BOpDdeg}.

\begin{lemma}\label{lLiecep} Let $f \in  \cc$ and $f_0$ be an element of minimal degree in $f+ I$. Then $\deg f_0 \le \deg h$, for any $h\in f+J$, hence $\deg f_0 \le \deg g$ for any $g\in f+K$.
\end{lemma}

\pp Arguing by contradiction, we consider $f_0 - h \in J$. By Lemma \ref{lLT}, $\overline{f_0} = \overline{f_0-h}$ contains as a subword the leading word $\bf$ of an element $f\in I$. We have $\overline{f_0}=u\bf v$.
Since $f_0 \in  \cc$, we have that $\bf\in\cw(\cm)$. As before, then also $u,v\in\cw(\cm)$. In this case again by Corollary \ref{BOlMBS}, the special bracketing on $w=u\bf v$ produces an element $\sigma(u,f,v)\in I$ whose leading word is the same as that of $f_0$. Subtracting $\sigma(u,f,v)$ from $f_0$ we find an element of $f_0+I$ whose leading word is strictly smaller than that of $f_0$, which is a contradiction.\epf

As in the case of associative algebras, an immediate corollary is the following. Recall that $I$ is an ideal of $\cc=\cl(\cn)$, $K$ the Lie ideal of $\cl(X)$ generated by $I$ and $J$ the two-sided associative ideal of $\ca(X)$ generated by $I$.

\begin{lemma}\label{BOlLCEP}
   $J \cap  \cc = I = K\cap \cc$.\epf
   \end{lemma}

The next result and its proof are Lie algebra analogues of Theorem \ref{BOpDdeg}. This is the best result about general tame filtrations on Lie algebras in this paper. Notice that if a Lie algebra $H$ is a Lie subalgebra in an associative algebra $R$ with tame filtration $\{ R_n\}$ then $\{ H\cap R_n\}$ is always a tame filtration in $H$.

\begin{theorem}\label{BOpLieDdeg} 
\begin{enumerate} 
\item[$\mathrm(1)$] A filtration $\chi$ on a Lie algebra $H$ is tame if and only if $\chi\sim\gamma\cap H$ where $\gamma$ is the degree filtration on a 2-generator Lie algebra $G$ where $H$ is embedded as a subalgebra, if and only if $\chi\sim\rho\cap H$ where $\rho$ is the degree filtration on a 2-generator associative algebra $R$ where $H$ is embedded as a Lie subalgebra.
\item[$\mathrm(2)$] A filtration $\chi$ on a Lie algebra $H$ is tame if and only if $\chi=\gamma\cap H$ where $\gamma$ is the degree filtration on a \fg Lie algebra $G$ where $H$ is embedded as a subalgebra, if and only if $\chi=\rho\cap H$ where $\rho$ is the degree filtration on a \fg associative algebra $R$ where $H$ is embedded as a Lie subalgebra.
\end{enumerate}
\end{theorem}

\pp As in all similar theorems (see notes after Definition \ref{dDFil} and just before this theorem), it suffices to prove the ''only if'' claim in both (1) and (2). Let us start with Claim (1). Suppose $\chi=\{ H_n\}$ is a tame filtration of $H$. Similarly to the approach of Theorem \ref{BOpDdeg}, we consider the free associative algebra $\ca(X)$, $\# X=2$, the non-overlapping set $\cm$ whose growth majorates an exponential function, the set $\cn$ of basic commutators obtained by setting brackets of words from $\cm$, a free associative algebra $\cb\cong\ca(\cm)$ and a free Lie algebra $\cc\cong\cl(\cn)$. The natural degree filtrations $\{ \ca_n\}$ of $\ca(X)$ and $\{ \cl_n\}$ of $\cl(X)$ by restriction induce tame filtrations $\{ \cb_n\}$ on $\cb$ and $\{ \cc_n\}$ on $\cc$.

Using the same constant $C$ and the argument that follows in Theorem \ref{BOpDdeg}, we define an epimorphism $\theta:\cc\to H$ such that $H_n=\theta(\cc_{Cn})$. Let $I$ be the kernel of this homomorphism, $K$ the Lie ideal of $\cl(X)$ generated by $I$ and $J$ the two-sided associative ideal generated by $I$. We set $G=\cl(X)/K$ and $R=\ca(X)/J$. Both these algebras have degree filtrations $\gamma=\{ G_n\}$ and $\rho=\{ R_n\}$ defined by the natural images of $X$. By Lemma \ref{BOlLCEP}, $H$ naturally embeds in both $G$ and $R$ and using Lemmas \ref{lLiecep} and \ref{BOlLCEP} in place of Lemma \ref{lcep} allows us to conclude that $\chi\sim\gamma\cap H$ and $\chi\sim\rho\cap H$. This takes care of Claim (1).

Now since constant $C$ in our present proof is the same as in the proof of Theorem \ref{BOpDdeg}, we can easily pass from present Claim (1) to present Claim (2) in the same manner as we did in the case of associative algebras. \epf



\subsection{Undistorted embeddings in simple algebras}\label{sUESA}

In this subsection we will discuss the possibility of undistorted embedding of an algebra as a subalgebra in a simple algebra. Note the many important results about embedding of algebras in simple algebras have been obtained by L.A.Bokut', starting with \cite{Bsim}, and his coauthors. One of the most up-to-date sources, which also contains an extensive list of references, is \cite{BCM} (the behavior of filtrations is not among the questions studied in those papers).

\begin{theorem}\label{BOpSA}
Any finitely generated associative, respectively, Lie algebra can be embedded without distortion in a 2-generator simple associative, respectively, Lie algebra.
\end{theorem}

\pp  First let $B$ be a \fg \emph{associative} algebra with generators $a_1,\ld,a_m$ and $\beta=\{ B_n\}$ the respective degree filtration. Let us add to the above generators two more: $x, y$, and consider an algebra $A$ given by the set of generators $a_1,\ld,a_m,x,y$, and a set of defining relations which is the union of the set $\mathcal{S}$ of ALL relations of $B$ and a set $\mathcal{R}$ of additional relations which we are going to define next. By $\alpha=\{ A_n\}$ we denote the degree filtration of $A$ defined by the generating set $\{ a_1,\ld,a_m,x,y\}$

To start with, we consider an auxiliary set $\cm$ in alphabet $\{ x,y\}$ satisfying the same conditions as the set $\mathcal{M}$ in Section \ref{sCA}: no word in $\cm$ is a proper subword of another word in $\cm$ and no proper prefix of a word in $\cm$ can be a proper suffix of another word in $\cm$. Next we introduce a well-order on the words in free monoid $\cw(a_1,\ld,a_m, x,y)$ as follows. First we set  $a_1<a_2<\ld<a_m<x<y$. For words of arbitrary length, we write $u<u'$ if either the length of $u$ is less than the length of $u'$ and in the case of equality, if $u$ is lesser than $u'$ lexicographically (ShortLex).

Now let $\{f_1,f_2,\ld\}$ be the list of all monomials in $a_1,\ld,a_m$, $x,y$. We will introduce defining relations, two at a time, for each $i=1,2,\ld$, provided that $f_i$ has no subwords equal to the leading words of relations in $\mathcal{S}$ and previously introduced relations of the set $\mathcal{R}$ under construction. If this condition is met,  we first introduce a new relation $u_if_iv_i=1$, where the ``wings''$u_i\ne v_i$ are arbitrary words in $\cm$ such that the length of each of the  $u_i$ and $v_i$ is at least twice the degree of $f_i$ and of any of previously introduced relations, with the exception of the relations from $\mathcal{S}$. Then we introduce $u_i'f_iv_i'=0$ with the same conditions on the new ``wings'' $u_i'\ne v_i'$ and additionally we will require that the length of each of the new wings $u_i',v_i'$ is at least two lengths of each of the old ones.

Notice that $\mathcal{S} \cup \mathcal{R}$ is complete under composition (see Definition \ref{dBCA}). Indeed, the set of ``old'' relations $\mathcal{S}$ was complete from start, as the set of ALL relations of $B$. Any ``new'' relation, that is, from $\mathcal{R}$, has the form $ufv - g$, where $g$ is one of polynomials $0$ or $1$. The leading word of such relation is $ufv$, where $u,v$ are monomials in $x,y$. Therefore no partial overlapping with the leading words of the ``old'' relations is possible. It could be possible that a leading word of an ``old'' relation is a subword in $f$. But we deliberately excluded such monomials $f$ in the process of construction of $\mathcal{R}$. Finally, the leading words of ``new'' relations cannot overlap or be subwords of each other, thanks to the choice of the set $\cm$ of words in the variables $x,y$. Thus, the closure of $\mathcal{S}\cup \mathcal{R}$ under composition is itself, as needed. 

To proceed further, first notice that $A$ is nonzero. Indeed, since $1$ does not contain any leading words of the relations in $\mathcal{S}\cup \mathcal{R}$, by Composition Lemma For Associative Algebras (Subsection \ref{sFMA}) it follows that $1\ne 0$ and so $A\ne\{ 0\}$. To show that $A$ is simple, it is sufficient to find, for each nonzero in $A$ polynomial $f$, two monomials $u,v$ such that $ufv=1$. Let us write $f$ as a linear combination of monomials $f(i)$ not containing leading words of relations in $\mathcal{S}\cup \mathcal{R}$. If the number of summands in this expression is 1 then by our construction, there are $u,v\in\mathcal{U}$ such that $\mathcal{S}\cup \mathcal{R}$ contains $ufv=1$. Then the ideal $I$ of $A$ containing $f$ must contain 1 and thus be equal to the whole of $A$. 

Let us now assume that $f$ is a linear combination of monomials $f(1)$, $f(2)$, $\ld$, $f(n)$ of the form $f=\lambda_1f(1)+\lambda_{2}f(2)+\cd+\lambda_nf(n)$ where all coefficients $\lambda_1,\ld \lambda_n$ are nonzero and $f(1)>f(2)>\ld>f(n)$, in the sense of the ordering introduced in the second paragraph of the present proof. One of the relations introduced by us was $u(1)f(1)v(1)=0$. Then $u(1)fv(1)$ will be the linear combination

\[
u(1)fv(1)=\lambda_2u(1)f(2)v(1)+\cd+\lambda_nu(1)f(n)v(1)\in I
\] 

of less than $n$ summands. We need to show that none of the monomials $u(1)f(i)v(1)$, $i=2,\ld,n$, in the previous equation contains a leading word of a relation in $\mathcal{S}\cup \mathcal{R}$, as a subword. Proving by contradiction, let us assume, say, that $u(1)f(2)v(1)$ contains such leading word $w$. If $w$ has no letters $x,y$ then $w$ is a subword in $f(2)$, which is impossible. Otherwise,  $w=uf_iv$, for some $i$, where the ``wings'' $u$ and $v$ are in the set of non-overlapping words $\cm$. As just above, we cannot have $w$ completely inside $f(2)$ and so either $u$ has overlapping with $u(1)$ (hence, equal to $u(1)$) or $v$ has overlapping $v(1)$(hence, equal to $u(1)$), or both. In either case, by the nature of defining relations in $\mathcal{R}$, we must have $f_i=f(1)$ and  $w=u(1)f(1)v(1)$. Since this is a subword of $u(1)f(2)v(1)$, the length of $f(2)$ must be at least the length of $f(1)$. Since $f(1)>f(2)$, $f(2)$ cannot be longer than $f(1)$. So their lengths are the same and hence $f(1)=f(2)$, a contradiction.

As a result, our ideal $I$ generated by $f$ contains a linear combination of $n-1$ nonzero monomials with nonzero coefficients which allows us to apply induction and conclude that $I=A$. 

To prove that the embedding of $B$ in $A$ is undistorted it is sufficient to show that if $f(a_1,\ld,a_m)\in B_k\setminus B_{k-1}$ for some $k$ and $f(a_1,\ld,a_m)=g(a_1,\ld,a_m,x,y)$ then $g(a_1,\ld,a_m,x,y)\not\in A_{k-1}$. In our proof by contradiction, we will additionally assume without loss of generality that $f(a_1,\ld,a_m)$ is different in $B$ from a polynomial with lesser leading word. Since the difference $h = g(a_1,\ld,x,y) - f(a_1,\ld,a_m)$ equals $0$ in $A$, we can apply Composition Lemma For Associative Algebras (Subsection \ref{sFMA}), and then the leading word of $h$ must contain as a subword the leading word of one of defining relations. Now by our assumption, $g(a_1,\ld,a_m,x,y)\in A_{k-1}$. Hence the leading word of $h$ equals the leading word of $f$ and then $f$ equals in $B$ to a polynomial with lesser leading word, which is a contradiction. 

Similar construction applies also in the case of Lie algebras. Again, given a Lie algebra $M$ with generators $a_1,\ld,	a_m$, $m\ge 2$, and ALL defining relations $\cs$, we add two new variables $x,y$ and define a Lie algebra $L$ by generators $a_1,\ld,a_m,x,y$ and a set of defining relations $\cs\cup\mathcal{R}$, where $\mathcal{R}$ is defined as follows. 
We call here $a_1,\ldots,	a_m,x,y$ suitable commutators; by induction 
$[u,v]$ is called suitable, if $u, v$ are suitable  and  $\bar u>\bar v$.
Notice that all basic commutators are suitable.
We enumerate all suitable commutators :  $g_1, g_2,\ldots$. The associative support of each suitable commutator is an associative regular word. Then we consider the set $\cm=\{ x^t(xy)^ty^2\,|\,t=2,3,4,\ld\}$, similar to the language $P$ in  \cite[Theorem 1]{BOAMG}, satisfying the non-overlapping condition. Next we totally order the free monoid $\cw(a_1,\ld,a_m,x,y)$ so that $x$ is the largest variable. We can easily observe that $\cm$ consists of regular words. Afterwards, we impose brackets on (regular) words of the set $\cm$ to produce the set $\cn$ of basic commutators $[u]$, for $u\in\cm$. 
	
Now for each $g_i$ such that $\bg_i$ has no subwords equal to the leading words of relations in $\cs$ and already introduced relations of $\mathcal{R}$, we add to $\mathcal{R}$ the total number of $m+3$ relations as follows. We consider $2(m+3)$ words $u_i^{(j)}, v_i^{(j)}\in\cm$, $j=0,1,\ld,m+2$, satisfying the following conditions. The length of each of $u_i^{(0)}, v_i^{(0)}$ is at least twice the length of $g_i$ or of any of previously introduced relations except those in $\cs$. Also the length of each of $u_i^{(j)}$ and $v_i^{(j)}$ is at least twice the length of each of  $u_i^{(j-1)}$ or $v_i^{(j-1)}$, for $j=1,\ld,m+2$. Finally, we require that $\bg_i\preceq v_i^{(j)}\prec u_i^{(j)}$, for all $j$. With these conditions in place, all the commutators $[u_i^{(j)}]$, $[v_i^{(j)}]$ obtained by setting brackets on these $2(m+3)$ words are suitable. Also, by our definition of suitable commutators, each commutator $[[[u_i^{(j)}],g_i],[v_i^{(j)}]]$ is suitable. Now the relations of $\mathcal{R}$ added on the $i\th$ step are  $[[[u_i^{(0)}],g_i],[v_i^{(0)}]]=0$, $[[[u_i^{(1)}],g_i],[v_i^{(1)}]]=a_1$, $\ld$, $[[[u_i^{(m+2)}],g_i],[v_i^{(m+2)}]]=y$. 

Notice that $\mathcal{S} \cup \mathcal{R}$ is closed under Lie composition (see Definition \ref{dLCC}). Indeed, the set of ``old'' relations $\mathcal{S}$ was complete from start, as the set of ALL relations of $M$. Any ``new'' relation, that is from $\mathcal{R}$, has the form $[[u,g],v]-h$, where $h$ is one of polynomials $0$ or $a_1,\ld,y$. The leading word of this relation is $u\bg v$, where $u,v\in\cm$. Therefore no overlapping with the leading words of the ``old'' relations is possible. It could be possible that a leading word of an ``old'' relation were a subword in $\bg$. But we deliberately excluded such commutators $g$ in the process of construction of $\mathcal{R}$. Finally, the leading words of ``new'' relations cannot overlap or be subwords of each other, thanks to the choice of the set $\cm$ of words in the variables $x,	y$. Thus, the closure of $\mathcal{S}\cup \mathcal{R}$ under composition is itself, as needed.

To proceed further, first notice that $L$ is nonzero. Indeed, since $x$ does not contain any leading words of the relations in $\mathcal{S}\cup \mathcal{R}$, by Composition Lemma for Lie Algebras (Subsection \ref{sFMA}), it follows that $x\ne 0$ and so $L\ne\{ 0\}$. To show that $L$ is simple, it is sufficient, for each nonzero Lie polynomial $f$,  to find $m+2$ monomials $u^{(1)},v^{(1)},\ld,u^{(m+2)},v^{(m+2)}\in\cm$ such that $[[[u^{(1)}],f], v^{(1)}]=a_1,\ld,[[[u^{(m+2)}],f], v^{(m+2)}]=y$. Notice that thanks to Composition Lemma for Lie Algebras (Subsection \ref{sFMA}), each elements of $L$ is a linear combination of basic commutators $g(i)$ such that $\overline{g(i)}$ has no subwords equal to the leading words of relations in $\mathcal{S}\cup \mathcal{R}$. If $f$ is a suitable commutator itself, then by our construction, such $u^{(1)},v^{(1)},\ld,u^{(m+2)},v^{(m+2)}\in\cm$ exist and then the ideal $I$ of $L$ containing $f$ must contain all generators of $L$, hence be equal to the whole of $L$. 

Any element $f\in L$ can be written as a \emph{reduced} linear combination $f=\lambda_1g(1)+\lambda_2g(2)+\cd+\lambda_ng(n)$  of suitable (even basic) commutators $g(1),\ld,g(n)$ with nonzero coefficients which additionally satisfies  $\overline{g(1)}> \overline{g(2)}>\ld>\overline{g(n)}$ and such that none of $\overline{g(1)},\overline{g(2)},\ld,\overline{g(n)}$ contains a leading word of a relation in $\cs\cup\mathcal{R}$, as a subword. Any reduced
combination of length $\ge 1$ is different from 0 in $L$ according to the Composition Lemma for Lie Algebras. We need to prove that if a reduced linear combination as above is an element of an ideal $I$ of $L$ then $I=L$. Let us use induction by $n$. We already handled the case $n=1$. Now suppose we already sorted out the case of linear combinations of length at most $n-1$ and we now deal with a linear combination $f$ of length $n>1$, as before. 

Now one of the relations introduced by us was $[[[u(1)],g(1)],[v(1)]]=0$, for appropriate $u(1),v(1)\in\cm$. Let us set $f'=[[[u(1)],f],[v(1)]]$. Then

\[
f'=\lambda_2[[[u(1)],g(2)],[v(1)]]+\cd+\lambda_n[[[u(1)],g(n)],[v(1)]]\in I.
\] 

Since by construction, $\overline{g(2)}<\overline{g(1)}\preceq v(1)\prec u(1)$, all commutators on the right hand side of the latter equation are suitable. The leading words of these commutators are regular words $u(1)\overline{g(2)}v(1)$, \ld, $u(1)\overline{g(n)}v(1)$. As in the case of associative algebras, none of these leading words contains the leading word of a defining relation in $\cs\cup\mathcal{R}$,  as a subword. We also have 
$\overline{[[[u(1)],g(2)],[v(1)]]}>\ld>\overline{[[[u(1)],g(n)],[v(1)]]}$.
As a result, any ideal $I$ containing $f$ contains also $f'$ which is a reduced linear combination of length $n-1$, and so $f'\ne 0$. This allows us to apply induction hypothesis and conclude that $I=L$.
 
The undistortedness of the embedding of $M$ in $L$ follows by exactly the same argument, as in the ``associative portion'', if we apply Composition Lemma for Lie Algebras (Subsection \ref{sFMA}) in place of its associative counterpart.
\epf

In combination with our previous results, the theorem just proved yields the following.

\begin{corollary}\label{BOtDdegSA}
A filtration $\beta$ on a unital associative, respectively, Lie algebra $B$ over a field $F$ is tame if and only if $\beta\sim\alpha\cap B$ where $\alpha$ is a degree filtration on a \fg \emph{simple} unital associative, respectively, Lie algebra $A$ over $F$.
\end{corollary}

\pp As usual, no need to prove the ``if'' claim of the theorem. Using Theorem \ref{BOpDdeg} in the case of associative algebras and Theorem \ref{BOpLieDdeg} in the case of Lie algebras, we embed $B$ in a finitely generated (associative, Lie) algebra $C$ with a degree filtration $\gamma=\{ C_n\}$ so that $\beta\sim\gamma\cap B$. This would mean that $\ds{\gamma\cap B}{\beta} \le t\cdot\mathrm{id}$, for some integer $t$. Afterwards, using Theorem \ref{BOpSA}, we embed $C$ in a simple (associative, Lie) algebra $A$ with degree filtration $\alpha$ so that $\gamma\sim\alpha \cap C$. Then $\ds{\alpha\cap C}{\gamma} \le u\cdot\mathrm{id}$, for some other integer $u$. It remains to apply Claim (3) of Proposition \ref{pGPDF} to see that $\ds{\alpha\cap B}{\beta} \le (tu)\cdot\mathrm{id}$. Thus $\beta\sim \alpha\cap B$, as claimed.\epf

\begin{remark}\label{BOrUIP} \emph{A quick look at our proof of Theorem \ref{BOpSA} reveals that we can replace the word ``simple'' in the statement by ``divisible'' which we understand as follows. A unital algebra $A$ is ``divisible'' if for any nonzero $a$ there exist $p,q\in A$ such that $paq=1$. Such algebras form a subclass of the class of simple algebras defined within first order logic. It is known that simple algebras do not form such a class. Note that the embeddings suggested in \cite{Bsim} also enjoy this property.}
\end{remark}


\section{Undistorted embeddings in finitely presented algebras}\label{sHEUA}

As we noted in Introduction, not every tame filtration $\beta$ of an infinite - dimensional \fg algebra $B$ is equivalent to the intersection $\alpha\cap B$ where $\alpha$ is a degree filtration of a finitely presented algebra $A$ where $B$ is embedded as a subalgebra. This explains the necessity of restrictions we impose on algebras and their filtrations in this section. 

For example, if a \fg algebra $B$ is embedded as a subalgebra in a finitely presented algebra $A$, both over a constructive field $F$, then it is well-known that $B$ can be defined by a recursively enumerable set of defining relations. So to cover all important cases we have to impose the condition of recursive enumerability of defining relations of $B$. Then, by the above remark, there are two many tame filtration on a given infinite(-dimensional) algebra, so we have to select a reasonably narrow countable subset of filtrations $\beta$ for which we could expect the desired results true (see constructive filtrations below). 

Finally, a very strong practical reason is that, in general, the embedding theorems of algebras in finitely presented algebras are known only in few cases, like groups (G. Higman's theorem  \cite{Hig}), semigroups (Murskii's theorem \cite{Mu-u-u}) and associative algebras (Belyaev's  theorem \cite{VB}). Unfortunately, Lie algebras are not on this list (see this story in \cite{KS}).

In the case of an algebra $B$ without the structure of a vector space (such as semigroups or groups, etc) the \emph{constructivity} of a finitary filtration $\beta=\{ B_n\}$ means that there exists an algorithm that lists out the pairs $(b_i,\deg_{\beta}b_i)$, containing all elements of $B$. In the case of linear algebras, this list must contain all pairs $(b_i,\deg_{\beta}b_i)$ where $b_i$'s are all elements of a \emph{$\beta$-basis} $\cb$. Here by $\beta$-basis $\cb$ of $B$ we mean a union $\cb=\cup_n\cb_n$ where $\cb_n$ is a basis of $B_n$ containing $\cb_{n-1}$.

We have the following ``constructive'' modifications of our results in the previous section. They serve as a tool in the proof of Theorems \ref{BOtBODtM} and \ref{BOtM}, the latter being our central result in this section.

\begin{lemma}\label{BOpDdegC} 
Let $B$ be an associative algebra over a field $F$ with a recursively enumerable set of defining relations. A constructive filtration $\beta$ on $B$ is a tame filtration if and only if $\beta=\alpha\cap B$ where $\alpha$ is a degree filtration on a \fg associative algebra $A$ with recursively enumerable set of defining relations, where $B$ is embedded as a subalgebra.
\end{lemma}

\pp We simply have to adapt the proof of Theorem \ref{BOpDdeg} but notice that the homomorphism $\vp$ in the proof follows the algorithm of enumerating the elements of the $\beta$-basis $\cb$ (with their $\beta$-degrees). Now the defining relations of $A$ are the defining relations of $B$ in which the elements of the free generating set $\cm$ of $\cb$ are replaced by their expressions as words in $X$. Since by construction we have an algorithm for enumeration of the preimages under $\vp$ of the elements of $B$, and an algorithm for enumeration of the defining relations for $B$, we have an algorithm for enumeration of defining relations for $A$.\epf

Similarly, in the case of monoids, the adaptation of the proof of Theorem \ref{pUEM}, gives the following.

\begin{lemma}\label{pUEMC} Let $N$ be a recursively presented monoid. A constructive filtration $\beta$ on $N$ is a tame filtration if and only $\beta\sim\alpha\cap N$ where $\alpha$ is the degree filtration in a \fg monoid $M$ with recursively enumerable set of defining relations, where $B$ is embedded as a submonoid.\epf
\end{lemma}


We start our treatment of the possibility of writing tame filtration in terms of degree filtrations of finitely presented algebras with the case of monoids.

\subsection{Monoids}\label{ssMu-u-u}

Suppose a monoid  $M$ is given by a set of generators $X$ and a recursively enumerable set of relations. In the paper \cite{Mu-u-u} it is proven that $M$ embeds in a finitely presented  monoid  $M'$ with a set of generators $X'$ that includes $X$. Actually, a more precise result is true.

\begin{theorem}\label{BOtMur}
 Any \fg monoid $M$ with a recursively enumerable set of defining relations can be embedded in a finitely presented monoid $M'$ as an undistorted submonoid.
\end{theorem}
 
\pp We will use the argument from \cite{Mu-u-u}. However, since we need additional information about the lengths $|\;\;  |_X$ and $|\;\;|_{X'}$ of the words in monoids $M$ and $M'$ arising in Murskii's proof with respect to their generating sets $X$ and $X'$, respectively, we need few comments about the course of Murskii's proof.

The alphabet $X'$ includes, in particular, a copy $\wX$ of alphabet $X$, hence each word $p$ in alphabet $X$ has a copy $\wp$. The next claim is contained in \cite[Lemmas 3.3 and 3.1]{Mu-u-u}. In monoid $M'$, if a word $p$ in alphabet $X$ is equal to a word $w$ in alphabet $X'$ then $w=p_0u_1p_1u_2\cdots  p_l$ so that:
\begin{enumerate}
\item[(1)] each $p_i$ is a word in the alphabet $X$;
\item[(2)] after deleting from each $u_i$ some subwords we obtain a word
$u'_i$, containing a subword of the form $\wq_i$, where $q_i$ is a word in alphabet $X$ such that the following holds:
\item[(3)] The word $p_0q_1p_1\cdots q_lp_l$  equals $p$ in monoid $M$.
\end{enumerate}

Using (1) - (3), one can compare the length $|p|_X$ of an element of $M$, represented by a word $p$ and the length $|w|_{X'}$,  where $w$ represents the same element in $M'$. We have:

\bea
&&|p|_X \le  |p_0|_X + |q_1|_X +\ld+ |p_l|_X \le |p_0|_{X'}+|q_1|_{X'}+...+|p_l|_{X'}\\
  &\le& |p_0|_{X'} + |u'_1|_{X'} +\cd + |p_l|_{X'}\le |p_0|_{X'}+|u_1|_{X'}+\cd+|p_l|_{X'} = |w|_{X'}.
\eqa

Since $w$ is an arbitrary word equal $p$ in monoid $M'$, we have an inequality $|p|_X  \le  |p|_{X'}$, that is, what we wanted to achieve.\epf

If we combine Theorem \ref{BOtMur} with Lemma \ref{pUEMC}, and use Claim (3) of Proposition \ref{pGPDF}, then we obtain the following.  

 \begin{theorem}\label{BOtBODtM}
A constructive filtration $\beta$ on a monoid $N$ is tame if and only if $\beta\sim\alpha\cap N$ where $\alpha$ is a degree filtration on a finitely presented monoid $M$ where $N$ is embedded as a submonoid.\epf
 \end{theorem}
 
 We now switch to to unital associative algebras.
 
\subsection{Unital associative algebras}\label{ssBack}

In this subsection we consider unital associative algebras with recursively enumerable set of defining relations over fields finitely generated over prime subfield.

It was proven in a paper by Belyaev \cite{VB} that every associative  algebra $A$ with a recursively enumerable set of defining relations, over a unital commutative and associative ring $K$ or a field $K$ that is finitely generated over its prime subfield, can be isomorphically embedded in a finitely presented algebra $B$ over $K$. Belyaev's proof does not guarantee that the embedding is undistorted, and neither is it unital. In what follows, we base on Belyaev's proof to produce a new proof ensuring that both properties are satisfied. 

\begin{theorem}\label{BT} Let $B$ be an arbitrary finitely generated unital associative algebra with a recursively enumerable set of defining relations, over a field $F$ which is finitely generated over prime subfield. Then there exists a finitely
presented unital associative $F$-algebra $A$ in which $B$ is contained
as an undistorted unital subalgebra.
\end{theorem}

To prove this theorem we will need to review and modify several Belyaev's preliminary results. 

\begin{lemma}\label{Bl1} Let $X$ be an alphabet that includes letters $a_1,\ld,a_n,b_1,\ld,b_n$. Suppose $I$ is a two-sided ideal of a free unital associative algebra $\ca(X)$ over a field $F$ with free generating set $X$. Let $\vp$ be a homomorphism of unital  algebras $\vp: \alg{\{ a_1,\ld,a_n\}}\to\alg{\{ b_1,\ld,b_n\}}$ such that  $\vp(a_i)=b_i$, for $i=1,\ld,n$, and $\vp(I\cap\,\! \alg{\{ a_1,\ld,a_n\}}) \subset  I$. Consider $X_1=X\cup\{x,y,z,\beta_1,\ld,\beta_n\}$ and naturally embed $\ca(X)$ in $\ca(X_1)$. Then there is a two-sided ideal $I_1$ in $\ca(X_1)$ containing all elements 
\begin{enumerate}
\item[\emph{(i)}] $xzy-1$;
\item[\emph{(ii)}] $xa_izy-b_i\,(i=1,\ld,n)$;
\item[\emph{(iii)}] $a_iz-z\beta_i\,(i=1,\ld,n)$;
\item[\emph{(iv)}] $xz\beta_i\beta_j-b_ixz\beta_j\; (i,j=1,\ld,n)$,
\end{enumerate}
such that $I_1\cap\ca(X)=I$. In addition, if $f\in\ca(X)$ and $f_0\in f+I$ is such that $\deg_{X}f_0\le \deg_{X}g$, for all $g\in f+I$, then $\deg_{X_1}f_0\le \deg_{X_1}h$, for all $h\in f+I_1$. 
\end{lemma}

\pp It will be convenient, for our further argument, to extend $\vp$ to a linear transformation of the vector space $\ca(X)$. Given a word in $\beta$ in the free monoid $\mathcal{W}(\beta_1,\ld,\beta_n)$, we denote by $\beta(\ba)$ the result of replacing each $\beta_i$ by
$a_i$ ($i=1,\ld,n$) in $\beta$. We define $I_1$ in $\ca(X_1)$ as the ideal generated by $I$, all elements in (iii) and (iv) in the statement of Lemma and also all elements $xwz\beta y-\vp(w\beta(\ba))$ where $w$, respectively, $\beta$ run through all words in $\wax$, respectively, $\mathcal{W}(\beta_1,\ld,\beta_n)$. Thus, the element of the latter kind include (i) and (ii), for obvious choices of $w$ and $\beta$. 

We prove our lemma arguing by contradiction. Let us choose $f\in\ca(X)$ and $f_0\in f+I$ of the least degree $d$ possible in $f+I$. Assume that there is an element $h\in\ca(X_1)$ of strictly lesser degree representing $f_0+I_1$. Let us consider the difference $v=f_0-h$. Then $v\in I_1\setminus \{ 0\}$ and so $v$ can be written as a linear combination, with coefficients in $F$, of the elements of the form

\bee{Be1}
v_1uv_2, w_1(xwz\beta y-\vp(w\beta(\ba)))w_2, w_1(a_iz-z\beta_i)w_2,w_1(xz\beta_i\beta_j-b_ixzb_j)w_2
\ene

where $u\in I$, $v_1,v_2,w_1,w_2\in\mathcal{W}(X_1)$. We may assume that $v_1$ has no suffix, and that $v_2$ has no prefix which is a letter of $X$. 

Let us view $x$, respectively, $y$ as a left, respectively, right parenthesis. In a word $w\in \mathcal{W}(X_1)$ with
``properly'' arranged parentheses, these are naturally divided into pairs $(x,y)$, a left parenthesis
$x$ and its corresponding right parenthesis $y$. The \emph{depth} of a fixed pair $(x,y)$ is the difference between the number of occurrences of $x$ and $y$ to the left of $x$ in $(x,y)$.
If $w$ has a pair of parentheses with depth $s$ but no pair with depth 
$s+1$ , then the number $s$ is called the \emph{rank} of $w$. If $w$ has no parentheses, its rank is zero.
It is easy to observe that if in the expression for $v$ we group together the monomials with
``properly'' arranged parentheses, then we again obtain a linear combination of the elements (\ref{Be1}). Since in the expression for $v=f_0-h$ there are no monomials of degree $\ge d$ with improperly arranged brackets, if we write $v$ as a linear combination of elements of the form (\ref{Be1}), then all monomials of degree $\ge d$ with such ``parentheses'' cancel, while the terms of degree $<d$ can be included in $h$ (remember that we proceed in our proof ``by contradiction''!), which changes $v$. As a result, from the very beginning, we may assume that in our linear combination of the elements (\ref{Be1}), we do not have monomials with ``improperly'' arranged parentheses.

Let $s$ be the largest number such that the expression for $v$ involves words of rank $s$. If $s=0$ then an argument similar to the one just given allows us to assume that $z$ is not present in our expression (those monomials of degree $\ge d$ with $z$ cancel out and those of degree $<d$ can be included in $h$). Thus $v$ is a linear combination of the elements $v_1uv_2$, where $v_1,v_2\in \cw(X\cup\{\beta_1,\ld,\beta_n\})$ (this case was omitted in \cite{VB}). However, all $v_1uv_2$ where $v_1v_2\not\in\cw(X)$ of degree $d$ must cancel because we do not have such monomials in $v=f_0-h$. Those with degree $<d$ can be moved to $h$. Now we have $v\in I$ and then $h= f_0 - v  = f_0\mod I$. Since $\deg h <d$, we obtain a contradiction with the choice of $f_0$.

Now suppose $s>0$. We will show that $v$ has a presentation of the
same form, in which all words have rank less than $s$. The proof of our lemma will then be complete by induction on $s$. 

A word $w\in \mathcal{W}(X_1)$ with properly arranged parentheses is called \emph{good} if its rank is
either less than $s$, or else equal to $s$ but for any pair $(x,y)$ of depth $s-1$ the subword starting at $x$ and ending at $y$, for these $x,y$, equals $xuz\beta y$ with $u\in\wax$ and $\beta\in \mathcal{W}(\beta_1,\ld,\beta_n)$. Otherwise we call $w$ \emph{bad}.

If we examine the elements of the last three types in (\ref{Be1}), we quickly observe that, in each type, both summands are good or bad words at the same time, which makes some of elements in (\ref{Be1}) (including the first type) good and some bad. Since $f_0\in\ca(X)$, the linear combination $b$ of all bad summands has degree strictly less than $d$, so replacing $h$ by $h-b$ removes all bad summands from our expression for $v$. So in our argument by contradiction we can assume that we do not have bad summands in the expression for $v$ through (\ref{Be1}).

In each word $w\in\cw(X_1)$ of rank $s$ and degree $\ge d$ in this expression of $v$ for each pair $(x,y)$ of depth
$s-1$ we replace the subwords $xuz\beta y$ for these $x,y$ by $\vp(u\beta(\ba))$. Since these subwords cancel while reducing to $f_0-h$, they should cancel also if we replace the respective subwords  $xuz\beta y$ by any symbol ? not in $X_1$. Therefore, the linear combination $v$ will not change if we replace everywhere these subwords by $\vp(u\beta(\ba))$. (Again, those $w$ with degree $<d$ will be moved to $h$.) 

We again obtain an expression for $v$ and now we are going to show that this, as before, is a linear combination of elements (\ref{Be1}). Since the ranks of words in the new expression are less than $s$, this
will complete the proof of the lemma.

Consider a summand $v_1uv_2$ where $u\in I$.  Suppose $u=\sum_i\alpha_iw_i$ where $\alpha_i\in F$ and $w_i\in\cw(X)$. Clearly,
the ranks of all of the words $v_1w_iv_2$ are the same. If their common rank is equal to $s$,
then under the replacement described above the words $w_i\in \mathcal{W}(X)$ are affected only when
$v_1= v_1'x, v_2= z\beta yv_2'$. Now after the replacement, we have
\bea \sum_i\alpha_iv_1^{''}\vp(w_i\beta(\ba))v_2^{''}=v_1^{''}\left(\sum_i\alpha_i\vp(w_i\beta(\ba))\right)v_2^{''}.
\eqa
But $\sum_i\alpha_iw_i\beta(\ba)\in I$, and by our hypothesis about $\vp$,  $\sum_i \alpha_i \vp(w_i\beta(\ba))\in I$.

Consider a summand $w_1(xwz\beta y-\vp(w\beta(\ba)))w_2$. If in the word $w_1\vp(w\beta(\ba))w_2$ the subword $\vp(w\beta(\ba))$
occurred within a pair of depth $s-1$, then clearly in the word $w_1xwz\beta yw_2$ the
pair $(x,y)$ would have depth $s$, which is impossible. Therefore, obviously, after the replacement
the expression under consideration either vanishes or keeps the same form.

Consider a summand $w_1(a_iz-z\beta_i)w_2$. It suffices to look at the case where $w_1= w_1'xw$, $w_2= \beta yw_2'$ and this pair $(x,y)$ has depth $s-1$. If this is the case, after the replacement, we arrive at 

\bea
w_1^{''}\vp(wa_i\beta(\ba))w_2^{''}-w_1^{''}\vp(wa_i(\beta_i\beta)(\ba))w_2^{''}.
\eqa
Since $(\beta_i\beta)(\ba)= a_i\beta(\ba)$, this expression vanishes.

Finally, consider a summand $w_1(xz\beta_i\beta_j-b_ixz\beta_j)w_2$. Again, it suffices to look at the
case where $w_2= \beta yw_2'$ and the pair considered $(x,y)$ has depth $s-1$. After the replacement
we obtain

\bea
w_1^{''}\vp((\beta_i\beta_j\beta)(\ba))w_2^{''}-w_1^{''}b_i\vp((\beta_j\beta)(\ba))w_2^{''}.
\eqa
This expression also vanishes, since by the choice of $\vp(w)$ we have

\bea
\vp((\beta_i\beta_j\beta)(\ba))= b_ib_jb_{i_1}\cdots b_{i_k}= b_i\vp((\beta_j\beta)(\ba)),
\eqa
if $\beta = \beta_{i_1}\cdots\beta_{i_l}$.
Now the proof is complete.\epf

It will be useful to restate Lemma \ref{Bl1} without involving free algebras.

\begin{lemma}\label{Br1}\emph{(1)} Let $A$ be an algebra with two subalgebras $\alg\{a_1,\ld,a_n\}$ and $\alg\{b_1,\ld,b_n\}$ such that there is an algebra homomorphism 

\[
\vp: \alg\{a_1,\ld,a_n\}\to\alg\{b_1,\ld,b_n\}
\]

satisfying $\vp(a_i)=b_i$, $i=1\ld,n$. Then $A$ can be embedded without distortion in an algebra $A_1$, whose generators are those of $A$ and some $x,y,z,\beta_1,\ld,\beta_n$, and whose defining relations of are those of $A$ and a finite number of relations \emph{(i)} through \emph{(iv)} from Lemma \emph{\ref{Bl1}}.
 
\emph{(2)} If $A_1$ is any algebra with elements $a_1,\ld, b_n$, $x,y,z,\ld,\beta_n$ satisfying \emph{(i)} through \emph{(iv)} from Lemma \emph{\ref{Bl1}}, and $f(a_1,\ld,a_n)=0$ is a relation in $A_1$ then also $f(b_1,...,b_n)=0$ is a relation of $A_1$.
\end{lemma}
\pp Claim (1) is just a restatement of Lemma \ref{Bl1}. 
To prove Claim (2), it is sufficient to check that the mapping  $\vp:\alg\{a_1,\ld,a_n\}\to\alg\{b_1,\ld,b_n\}$ given by $\vp (a)=xayz$ is a homomorphism of algebras. Since $\vp$ is linear, we only need to check the claim when $a=a_{i_1}\cd a_{i_k}$: 
\bea
&&\vp(a_{i_1}a_{i_2}\cdots a_{i_k})
=xa_{i_1}a_{i_2}\cdots a_{i_k}zy
=xz\beta_{i_1}\beta_{i_2}\cdots \beta_{i_k}y\\&=&b_{i_1}b_{i_2}\cdots b_{i_{k-1}}xa_{i_k}zy=b_{i_1}b_{i_2}\cdots b_{i_k}=xa_{i_1}zyxa_{i_2}zy\cdots xa_{i_k}zy\\&=&\vp(a_{i_1})\vp(a_{i_2})\cdots\vp(a_{i_k}):
\eqa
Also, applying relation $xzy = 1$, we easily derive $\vp(1)=1$.
\epf

The next result is a slight modification of Belyaev's lemma devoted to the development of the technique which allows one to switch from ``additive'' defining relations of an associative algebra to ``multiplicative'' ones.

\begin{lemma}\label{Bl2} Let $f(i,j)$ be a recursive function defined for $i,j=1,2,\ld,i\ne j$ 
such that $f(i,j)=f(j,i)$. Suppose that $Y\subseteq \N^2$ is a recursively enumerable set and $Q$ a unital algebra over a field $F$ with generators $x,y,z$ and  defining relations

\[
\{ xy^iz+xy^jz=xy^{f(i,j)}z\,|\, i\ne j;\,i,j=1,2,\ld\}\cup\{xy^iz=xy^jz\,|\,(i,j)\in Y\}\cup\{xz = 1\}.
\]

Then there exists an algebra $Q_1$  over a field $F$ with the following properties:
\begin{enumerate}
\item[\rm(1)]  $Q_1$ has a finite number of generators and a recursively enumerable set of defining
relations, one of which has the form $\alpha+\beta=\gamma$ and the others are equalities of words in an  appropriate alphabet;
\item[\rm(2)] $Q$ is a subalgebra of $Q_1$ and the degree of any element in $Q$, with respect to the generating system of $Q$, cannot decrease when we consider this element as an element of $Q_1$, with respect to a generating system of $Q_1$.
\end{enumerate}
\end{lemma}
\pp Let us define an increasing sequence of natural numbers \bea
n_2(1),n_3(1),n_3(2),n_4(1),n_4(2),n_4(3),n_5(1),n_5(2),n_5(3),n_5(4),\ld
\eqa
 by setting $n_2(1)=f(1,2)$ and if $n_j(i)$ immediately follows after $n_k(l)$ then we set $n_j(i)=\max\{ n_k(l)+1,\,f(i,j)\}$. We also define $n_j(i)$ for $j\le i$ by setting $n_j(i)=n_{i+1}(j)$. Now for $i<j$ we set $s(i,j)=s(j,i)=n_j(i)$.

Thus, for any $n\in\N$ there
exists at most one pair $(i,j)$, $1\le i< j$, such that $n=s(i,j)=s(j,i)$.  With this choice we would automatically have $n_i(j) \ge  i-1$, for all $j$. 

Let us define the desired algebra $Q_1$ by the set of generators $\{ x,y,z,u,\alpha,\beta,\gamma\}$ and
the set of defining relations 
\bea
&&\{ xy^iz=xy^jz\,|\, (i,j)\in Y\}\cup \{ xy^iz=xu^{n_i(j)}\ve_{ij}z\,|\, i,j=1,2,\ld\}\\
&&\cup\, \{ xy^{f(i,j)}z=xu^{s(i,j)}\gamma z\,|\, i\ne j=1,2,\ld\}\cup\{xz=1\}\cup \{\alpha+\beta=\gamma\}.
\eqa
Here $\ve_{ij}$ equals $\alpha$ if $i+j$ is even and $\beta$ if $i+j$ is odd.

First we note that the defining relations of $Q_1$ imply those of $Q$. Indeed, if $i<j$,
then it follows from the relations of $Q_1$ that
\bea xy^iz+xy^jz&=&xu^{n_i(j-1)}\ve_{i,j-1}z+xu^{n_j(i)}\ve_{ij}z=xu^{s(i,j)}(\ve_{ij}+\ve_{i,j-1})z\\&=&xu^{s(i,j)}\gamma z=xy^{f(i,j)}z.
\eqa 

Now consider the ideal $I_1$ of the free $F$-algebra $\fx{x,y,z,u,\alpha,\beta,\gamma}$ generated by
the relations of $Q_1$. Any element of this ideal can be written as a linear combination of the elements of the form

\beq{Be1p411}
&&w_1(xy^iz-xy^jz)w_2\,((i,j)\in Y),\;
w_1(xy^iz-xu^{n_i(j)}\ve_{ij}z)w_2,\\ &&w_1(xu^{f(i,j)}z-xy^{s(i,j)}\gamma z)w_2\,(i\ne j),\; w_1(xz-1)w_2,\; w_1(\alpha+\beta-\gamma)w_2,\nonumber
\eqe
where $w_1,w_2\in \mathcal{W}(x,y,z,u,\alpha,\beta,\gamma)$.

By a reduction of a word $w\in \mathcal{W}(x,y,z,u,\alpha,\beta,\gamma)$ we will mean the simultaneous replacement in $w$ all subwords of the following forms:
\begin{itemize}
\item[] $xu^n\alpha z$ by $xy^iz$, where $i$ is such that $n=n_i(j)$ and $i+j$ is even or
\item[] $xu^n\beta z$ by $xy^iz$, where $i$ is such that $n=n_i(j)$ and $i+j$ is odd or
\item[] $xu^n\gamma z$ by $xy^{f(i,j)}z$, where $i\ne j$ is such that $n=s(i,j)$ or 
\item[] $xz$ by 1.
\end{itemize}

It should be stressed, that the kind of replacement we apply is fully defined, that is, knowing the word uniquely defines by what word it has to be replaced. For instance, $xu^n\alpha z$ should be replaced using the reduction in the first line, etc.

As in Lemma \ref{Bl1}, we proceed by contradiction. Again, suppose we have $v=f_0 - h$, with the same conditions on $f_0$ and $h$.  The reductions of the four kinds just described should be applied to the monomials of degree $\ge d$. The monomials of the form (\ref{Be1p411}) of degree $<d$ should be attributed to $h$. It is important here that when we apply reduction, the degree of monomials does not increase. This follows from the choice of numbers $s(i,j) = n \ge \max\{f(i,j),\,  i-1\}$. 
 
Having completed all reductions, we arrive at $v$ which is the linear combinations of the words of the form $w_1(\alpha+\beta-\gamma)w_2$. It follows that the leading word of the polynomial $f$ (which does not depend on $\alpha, \beta, \gamma$) cannot cancel, which is a contradiction. Thus the proof is complete.\epf

Now we proceed to the proof of Theorem \ref{BT}.

\pp Suppose $B$ is an arbitrary unital countable $F$-algebra, $S$ a finite generating set for $B$ and $b_1,b_2,\ld$ an arbitrary enumeration of all of its elements (each element occurs at least twice). It is easy to see that after a possible renumeration we can have $\deg_Sb_i\le i$, for all $i=1,2,\ld$. By Corollary \ref{BOrABC}, there exists a unital $F$-algebra $C$ with three generators $a,b,c$, with a recursively enumerable set of defining relations, such that $B\subset C$, $1=ac$, $b_i=ab^ic$, for $i=1,2,\ld$, and this embedding has no distortion. 

Let $Y= \{(i,j)\,|\,ab^ic=ab^jc\}$. Then $Y\subset\N^2$ is a recursively enumerable set. Let $f(i,j)$
be a recursive function, defined for all $i,j\in\N$ with $i\ne j$, such that $f(i,j)= f(j,i)$ and

\bea
ab^ic+ab^jc=ab^{f(i,j)}c\mbox{ in } C.
\eqa
Consider the $F$-algebra $Q$ with generators $x,y,z$ and defining relations

\bea
\hspace{-.2in}\{ xy^iz+xy^jz=xy^{f(i,j)}z\,|\, i\ne j;\,i,j=1,2,\ld\}\cup\{xy^iz=xy^jz\,|\,(i,j)\in Y\}.
\eqa

By Lemma \ref{Bl2}, $Q$ is contained as an undistorted subalgebra in an $F$-algebra $Q_1$ which has a recursively enumerable set of defining relations

\bee{Be1p413}
\Sigma(x,y,z,u,\alpha,\beta,\gamma)\cup\{\alpha+\beta=\gamma\},
\ene
where $\Sigma$ contains only word equalities in $x,y,z,u,\alpha,\beta,\gamma$.

Next we consider the tensor product $Q_1\o C$, where $Q_1$ is embedded as the set of elements $\{ q_1\o 1\,|\, q_1\in Q_1\}$ and $C$ as the set of elements $\{ 1\o c\,|\, c\in C\}$. Both embeddings are easily seen to be undistorted.

There exists a homomorphism of algebras $\vp:Q\to C$. By Claim (1) of Lemma \ref{Bl1}, $Q_1\o C$ can be embedded as undistorted subalgebra in an $F$-algebra $Q_2$ with additional finite set of defining relations.

Now let us consider the monoid $G$ with generators $x',y',z',u',\alpha',\beta',\gamma'$ and set of defining relations $\Sigma(x',y',z',u',\alpha',\beta',\gamma')$. By Theorem \ref{BOtMur}, there exists a \fg monoid
$G_1$ containing $G$ as an undistorted submonoid, with a finite set of defining relations

\bee{Be3p413}
\Sigma_1(x',y',z',u',\alpha',\beta',\gamma').
\ene

Let $F[G]$ denote the monoid algebra $G$. By Proposition \ref{BOpGA}, $F[G]$ is a unital undistorted subalgebra of $F[G_1]$. By Claim (1) of Lemma \ref{Bl1},  $Q_2\o F[G_1]$ can be embedded as an undistorted subalgebra in an $F$-algebra $Q_3$ with additional finite set of defining relations.

At this stage of the proof, we have obtained a chain of undistorted embeddings 

\[
B\to C\to C\o Q_1\to Q_2\to Q_2\o F[G_1]\to Q_3.
\] 

We know that all additive relations of $B$, denoted by $R(B,\mathrm{ad})$, follow from the set of all additive relations $ab^ic+ab^jc=ab^{f(i,j)}c$ (the addition table) of $C$, denoted by $R(C, \mathrm{ad, left})$. Next we identify $C$ with $C\o 1$ in $C\o Q_1$, and denote by $R(C, \mathrm{ad, right})$ the relations of the set $R(C, \mathrm{ad, left})$ written in terms of generators $x,y,z$. The relations in $R(C, \mathrm{ad, right})$ follow from from some set $\Sigma \cup \{ r\}$ of relations of $Q_1$ identified with $1\o Q_1$,where $\Sigma$ consists only of some equalities of words while  $\{ r\}$ is the set of just one singular relation. By Claim (2) of Lemma \ref{Bl1} all relations in $R(C, \mathrm{ad, left})$ follow from relations of $Q_2$ of the form $\Sigma \cup R_2$ where $R_2$ is a finite set. 

 Next, all relations $\Sigma(x',\ld)$ of $F[G]$ follow from a finite number of relations of $F[G_1]$ and hence, again by Claim (2) of Lemma \ref{Bl1}, all relations $\Sigma(x,\ld)$ of $Q_2$ follow from a finite number of relations of $Q_3$. As a result, all relations in $R(B,\mathrm{ad})$ follow from  a finite number of relations of $Q_3$. 
 
A similar, actually, even simpler, argument works for the multiplication table $R(B,\mathrm{mul})$ of $B$. Indeed, let $H$ be a monoid given by relations $\bar a \bar b^i \bar c \bar a \bar b^j\bar c = \bar a \bar b ^{g(i,j)}\bar c$, where $g$ is a recursive function, whose existence follows from the recursive enumerability of the set of relations $b_ib_j=b_k$ in $B$. We embed $Q_3$ in $Q_3\o F[H]$. By Theorem \ref{BOtMur} and Proposition \ref{BOpGA} $F[H]$ is embedded without distortion in a finitely presented monoid algebra $F[H_1]$, and $Q_3$ is  undistorted in $Q_3 \o F[H_1]$. By Claims (1) and (2) of Lemma \ref{Bl1} this tensor product can be embedded without distortion in an algebra $Q_4$ in such a way that all relations in $R(B,\mathrm{mul})$ follow from a finite number of relations of $Q_4$.

Finally, let us notice that all relations of $B$ follow from its addition and multiplication tables and finitely many relations of the form $\alpha x = x'$, where $\alpha$ is a generator of $F$, $x$ a generator of $B$ (and $x' \in B$). Thus, all relations of $B$ follow from a finite number of relations of an algebra $Q_4$, where $B$ is contained as subalgebra. By Claim (viii) of Proposition \ref{BOr1}, $B$ is an undistorted subalgebra in an algebra $A$ given by these relations of algebra $Q_4$.
\epf

This theorem allows us to prove our final result characterizing tame filtrations in associative algebras.

\begin{theorem}\label{BOtM} Let $R$ be a unital associative algebra over a field $F$, which is finitely generated over prime subfield. A constructive filtration $\rho$ on $R$ is a tame filtration if and only if $\beta\sim\alpha\cap R$ where $\alpha$ is a degree filtration on a finitely presented unital algebra $A$ in which $R$ is embedded as a unital subalgebra.
\end{theorem}

\pp 
As in several cases before, we can restrict ourselves to the proof of ``only if'' portion of the statement. Using Theorem \ref{BOpDdegC}, we can embed $R$ without distortion in a finitely generated algebra $B$ presented by a recursively enumerable set of relations in such a way that $\rho\sim\beta\cap R$ where $\beta$ is a degree filtration on $B$. By definition then $\ds{\beta\cap R}{\rho}\le t\cdot\mathrm{id}$, for some $t\in\N$. Then we use Theorem \ref{BT} and embed $B$ without distortion in a finitely presented algebra $A$ with a degree filtration $\alpha$. This allows us to write $\ds{\alpha\cap B}{\beta}\le t\cdot\mathrm{id}$. Now by Claim (3) in Proposition \ref{pGPDF}, we can write $\ds{\alpha\cap R}{\rho}\le \ds{\alpha\cap B}{\beta}\cir\ds{\beta\cap R}{\rho}\le t^2\cdot\mathrm{id}$ or that $\alpha\cap R\preceq \rho$. Using  Claim (iv) in Proposition \ref{BOr1}, we obtain $\rho\sim\alpha\cap R$. The proof is now complete. \epf

\subsection{More Examples: Distortion of ``cyclic'' subalgebras}\label{ssEDF}

We start this subsection by exhibiting the variety of tame degrees on associative and Lie algebras of particularly simple form. Using the results on the distortion in groups \cite{AO}, one can show that for any real
number $\theta$ such that $0 < \theta< 1$ there is a tame degree on $F[x]$, which is equivalent
to $n^{\theta}$ (where $n$ is the ordinary degree of a polynomial). As it turns out, these
degrees come from the tame degrees on the free monoid $\cw(x)$. In the remainder of the subsection, we show that actually, there are tame filtrations (hence tame degrees) on $F[x]$ that are not induced by tame filtrations on $\cw(X)$.

Let $B=F[x]$ be the polynomial algebra in one variable. Consider any function $\vp$ on $\{ 0, 1,\ld\}$ which is positive on $\{ 1,2,\ld\}$. Let us define on $B$ a ``degree-like'' function  $d^\vp(f) = \vp(\deg(f))$. Setting $B_n^\vp=\{ f\in B\,|\, d^{\vp}(f)\le n\}$ defines a filtration on $B$ provided that $\vp$ is \emph{subadditive}, that is, $\vp (a+b) \le \vp(a) +\vp(b)$. 

If the ratio $\dfrac{\vp(n)}{\log n}$ cannot be separated from 0 when $n$ grows indefinitely then for any $k>0$ there is $n=n(k)$ such that $B_n$ contains all polynomials $f$ with $\deg f > \exp (kn)$; then $\{ B_n\}$ is not a tame filtration.

But if $\log n = O(\vp(n))$, for a subadditive function $\vp$, as above, then $\{ B_n\}$ is a tame filtration.

By Theorem \ref{BOpDdeg}, $B$ can be embedded in a 2-generator algebra $A^{\vp}$ with a degree filtration $\{ A_n^\vp\}$ so that $\{ B\cap A_n^\vp\}\sim\{B_n^\vp\}$. In this case, the degree of a polynomial $f \in F[x]$ with respect to the generators of $A$ will be a function equivalent to $\vp(\deg (f))$. For example, we may choose $\vp=(\log n )^{\mu} n ^{\nu}$, where $0<\nu<1$.


If $\vp$ is computable then by Theorem \ref{BOtM} such degree-like functions can be achieved in finitely presented algebras. For example, we can get the degree function for $F[x]$ in the form $(\log n )^{\mu} n ^{\nu}$, where $0<\nu<1$ in a finitely presented algebra provided that $\mu$ and $\nu$ are constructive real numbers (say, $\pi-e$, etc.). 
   
 A similar example can be obtained in the case of Lie algebras. Let $M$ be an abelian Lie algebra with basis $\{ a_0,a_1,a_2,\ld\}$. We form a semidirect product $L=Fx\oplus M$ of this algebra with a one dimensional algebra $Fx$ so that $[x,a_i]=a_{i+1}$ (similar algebras have been introduced in subsubsection \ref{sSF}). However, in this case, using Theorem \ref{BOpLieDdeg}, we can only guarantee an embedding in a finitely generated Lie algebra, which does not need to be finitely presented.
 
It is obvious that any subadditive ``degree-like'' function $d^{\vp}$ defines filtration not only on $F[x]$ but also on monoid $\mathcal{W}(x)$ and, conversely, is defined by its restriction to this monoid. But one should not think that even in a ``simple'' case like $F[x]$, all tame filtrations are equivalent to  ``degree-like'' filtrations, that is, filtrations extending the filtrations of $\mathcal{W}(x)$.

\begin{proposition}\label{BOpNDLF}
There exist a family of pairwise inequivalent tame filtrations on $F[x]$, labeled by nonzero elements of $F$, none being an extension of a tame filtration on $\cw(x)$, naturally embedded in $F[x]$.
\end{proposition}

\pp We start with an infinite increasing sequence of natural numbers $d_1=1, d_2, d_3,\ld$ such that $\dfrac{d_{n+1}}{d_n}\ra\infty$. There is $n_0$ such that $\max\left\{ \dfrac{d_m}{m},1\right\} < \dfrac{d_n}{n}$ for all $n\ge n_0$ and all $m<n$.

Let us consider a polynomial ring $F[y_1,\ld,y_s,\ld]$ in the set $\{ y_1,\ld,y_s,\ld\}$ of \emph{weighted} variables such that the weight of each $y_s$ is $s$. This assignment of the weight to the variables leads to a well-defined assignment of the weight to each polynomial $G(y_1,\ld,y_s)\in F[y_1,\ld,y_s,\ld]$. Let us fix $\lambda\in F$ and consider polynomials $f_n^\lambda(x) = x^{d_n}+\lambda x^{d_n-1}$, for all $n=1,2,\ld$.  Since $f_1 = x+\lambda$, any polynomial $f=f(x)\in F[x]$ can be written as $f(x)=G(f_1^\lambda(x),\ld,f_s^\lambda(x))$, where $G(y_1,\ld,y_s)\in F[y_1,\ld,y_s,\ld]$.  Let $D(f)$ be the minimum of the weights of all polynomials  $G(y_1,\ld,y_s)$ such that $f(x)=G(f_1^\lambda(x),\ld,f_s^\lambda(x))$. Clearly, setting $B_n^\lambda=\{f\,|\, D(f)\le n\}$ defines a filtration $\beta^\lambda=\{ B_n^\lambda\}$ on $B=F[x]$. Let us denote by $p_n$ the number of monomials of weight $n$ in $F[y_1,\ld,y_s,\ld]$. The set of such monomials is the disjoint union $S_1\cup\ld\cup S_n$ where the monomials in $S_1$ contain $y_1$, the monomials in $S_2$ do not contain $y_1$ but contain $y_2$, etc. Thus, $p_n\le p_{n-1}+ p_{n-2}+\cdots+p_0$ and using induction by $n$, we easily obtain $p_n \le 2^n$. It follows that $\beta^\lambda$ is a tame filtration of $B=F[x]$.

Now let us estimate from below the $\beta^\lambda$-degree (we will be simply saying $\lambda$-degree and write $\deg^\lambda$) of $f=f_n^\mu$ $(=x^{d_n}+\mu x^{d_n-1})$, where $n>n_0$ and $\mu\ne\lambda$. Suppose that $f=G(f_1^\lambda(x),\ld,f_s^\lambda(x))$ where the weight of $G(y_1,\ld,y_s)$ is minimal possible, as above. We need to consider three cases. 

\emph{Case} 1: $s<n$. When we replace in $G(y_1,\ld,y_s)$ each $y_m$ by $f_m^\lambda$, we obtain a polynomial whose degree does not exceed $\max\left\{\dfrac{d_1}{1},\ldots, \dfrac{d_s}{s}\right\}$ times the weight of $G(y_1,\ld,y_s)$. Since $\dfrac{d_m}{m} \le \dfrac{d_{n-1}}{n-1}$, we conclude that the weight of $G(y_1,\ld,y_s)$ cannot be less than $c_n = \dfrac{d_n}{d_{n-1}/(n-1)}$.  It follows that the $\lambda$-degree of $f_n^\mu$ cannot be less than $c_n$. Since $\dfrac{d_{n}}{d_{n-1}}\ra\infty$, we have $\dfrac{c_n}{n}\ra \infty$, for such $n$.

\emph{Case} 2: $s=n$. In this case, $G(y_1,\ld,y_n) =y_n^t g_t +  y_n^{t-1} g_{t-1}+\ld$, where $t>0$, each $g_j$ is a polynomial in $y_1,\ld,y_{n-1}$  and $g_t$ is nonzero. Then, as in the previous case, we replace all $y_1,\ld,y_n$ by $f_1^\lambda,\ld,f_n^\lambda$ and obtain $f_n^\mu$. If $g_t$ is not a constant or if $t>1$ then the leading term (as a polynomial in $x$) resulting from the first summand, has degree greater than $d_n$. This term can get canceled with a term arising in another $g_m(f_1^\lambda,\ld,f_{n-1}^\lambda)$ only when the degree of  $g_m(f_1^\lambda,\ld,f_{n-1}^\lambda)$ is at least $d_n$ and then, as in the previous case, the weight of $g_m$ and $G$ cannot be less than $c_n$. Now if $t=1$ and $g_1$ is a constant $c\ne 0$, then the first summand produces $cf_n^{\lambda}$ and in order to obtain $f_n^{\mu}$ on the left hand side, the weight of $g_0$ (hence of $G(y_1,\ld,y_n)$) must be greater or equal to $\dfrac{d_n-1}{d_{n-1}/(n-1)}= c_n-o(1)$.

\emph{Case} 3: $s>n$. Then $G =y_s^t g_t +  y_s^{t-1} g_{t-1}+\ld$, where $t>0$,  each $g_j$ is a polynomial in $y_1,\ld,y_{s-1}$ and $g_t$ is nonzero. In this case, the same argument as in Case 2 shows that the weight of $G(y_1,\ld,y_s)$ is at least $c_s$ and $\dfrac{c_s}{s} > \dfrac{c_n}{n}$. 

As a result, in any case, $\deg^{\lambda}f_n^\mu$ grows faster than $n$ (in the sense that $\dfrac{\deg^{\lambda}f_n^\mu}{n}\ra\infty$), while $\deg^{\mu} f_n^\mu\le n$. It follows that $\beta^\lambda\not\sim\beta^\mu$ if $\lambda\ne\mu$, as claimed.

In a particular case where $\mu=0$, we also have the following.
If we restrict any $\beta^\lambda$ with $\lambda\ne 0$ to $\mathcal{W}(x)$, we obtain a filtration $\beta$ whose extension to $F[x]$, viewed as monoid algebra, produces a filtration $\beta'=F[\beta]$, for which $\deg_{\beta'} f^{\lambda}_i \ge \deg_{\beta'} x^{d_i} = \deg_{\beta} f^{0}_i = \deg^{\lambda} f^0_i$ grows faster than $\deg^{\lambda}f_i^\lambda$. Since the restriction  of $\beta'$ to $\mathcal{W}(x)$ coincides with $\beta$, it follows that $\beta^\lambda$ is not equivalent to any filtration, which is extended from $\mathcal{W}(x)$. In other words, $\beta^\lambda$ is not a ``degree-like'' filtration, for any $\lambda\ne 0$.\epf

Our last remark is as follows. 

\begin{remark} \emph{By Theorem \ref{BOpDdeg}, if we fix $\lambda\in F$ then every filtration $\beta^\lambda$ on $B=F[x]$ is a restriction of a degree filtration $\alpha^\lambda$ of certain finitely generated algebra $A^\lambda$ where $B$ is contained as a subalgebra. However, none of $A^\lambda$ can be chosen commutative. Indeed, by Theorem \ref{Op1} about the embeddings in the commutative case, the embedding $B\subset A^\lambda$ is undistorted. If $\beta$ is the standard degree filtration on $B=F[x]$ then by Claim (i) of Proposition \ref{BOr1}, $\alpha^\lambda\cap B\sim\beta$. It follows that $\beta^\lambda\sim\beta$, which was proven impossible because $\beta$ is the extension of the standard degree filtration of $\cw(x)$ while $\deg^{\lambda} f^{\lambda}_n \le n$ and $\deg f^{\lambda}_n=d_n$ where  $\dfrac{d_n}{n}\to \infty$ as $n\to\infty$. }
\end{remark}
\begin{center}
\textbf{Acknowledgment}
\end{center}

The authors would like to thank the following mathematicians with whom they discussed various matters related to this paper: L. A. Bokut', V.  Drensky, R. I. Grigorchuk, M. V. Sapir, U. U. Umirbaev. We also thank Qiuhui Mo who pointed out an error in the Lie algebra portion of the original proof of Theorem \ref{BOpSA}.








\end{document}